\documentclass{amsart}

\usepackage{hyperref}
\usepackage{tikz-cd}
\usepackage{amsfonts,amsmath,amssymb,mathrsfs, stmaryrd, amsthm}

\usepackage{enumitem}
\usepackage{cite}
\usepackage{amsrefs}
\usepackage{bm}

\title{On the algebraic Mackey deformation}
\author{Yvann Gaudillot-Estrada}
\address{Université de Lorraine, CNRS, IECL, F-57000 Metz, France}
\email{yvann.gaudillot-estrada@univ-lorraine.fr}

\newtheorem{theorem}{Theorem}[section]
\newtheorem{theoremintro}{Theorem}
\newtheorem{proposition}[theorem]{Proposition}
\newtheorem{lemma}[theorem]{Lemma}
\newtheorem{corollary}[theorem]{Corollary}
\newtheorem{corollaryintro}[theoremintro]{Corollary}

\theoremstyle{definition}
\newtheorem{definition}[theorem]{Definition}
\newtheorem{remark}[theorem]{Remark}

\newtheorem{convention}{Convention}

\newcommand{\hatcpt}[1]{
  \vphantom{#1}\smash[t]{\widehat{#1}}
}
\begin{document}

\begin{abstract}
    This paper provides a conceptual understanding of the Mackey bijection from the admissible dual of a real reductive group $G$ into the admissible dual of its Cartan motion group $G_0$, by studying the primitive spectrum of a Hecke algebra over $\mathbb{C}[t]$ arising from the contraction of $G$ to $G_0$. If $K$ denotes a maximal compact subgroup of $G$, our approach leads to the following two results: the admissible Mackey bijection is a piecewise homeomorphism, each piece being defined by the length of the lowest $K$-types; the Hecke algebras of $G$ and $G_0$ have, up to $K$-types not exceeding a certain length, the same periodic cyclic homology.
\end{abstract}

\maketitle

\section*{Introduction}

On the basis of conjectural ideas of Mackey \cite{MackeyAnalogy} and Higson's previous work \cite{MAhigson} on complex groups, Afgoustidis constructed an explicit bijection \cite{Afgoustidis1}, the \textit{Mackey bijection}, between the tempered dual of any real reductive group $G$ and the unitary dual of its Cartan motion group $G_0$. This correspondence between two different unitary representation theories can be strengthened in several ways, especially by considering the reduced group C*-algebras \cites{Afgoustidis3,Membedding, AfgoustidisClare}, but not only \cites{AfgoustidisAubert, Yu}

However, concerning such a correspondence at the level of admissible representations, nothing is known except that the Mackey bijection has a natural extension to admissible duals \cite{Afgoustidis1}*{§5}. If $K$ denotes a maximal compact subgroup of $G$, the Mackey bijection can thus be seen a one-to-one correspondence between isomorphic classes of simple $(\mathfrak{g},K)$-modules and isomorphism classes of simple $(\mathfrak{g}_0,K)$-modules. The goal of this paper is to bring a conceptual understanding of the latter by studying an algebraic deformation of the Hecke algebra \cite{KnappVogan}*{Ch.I} of the pair $(\mathfrak{g},K)$ onto the Hecke algebra of the pair $(\mathfrak{g}_0,K)$. Our approach is partly inspired by the algebraic formalism introduced by Bernstein--Higson--Subag in \cite{BHS} for deformation families of Harish-Chandra pairs.

Let us denote by $U$ and $U^0$ the Hecke algebras associated to the pairs $(\mathfrak{g},K)$ and $(\mathfrak{g}_0,K)$, respectively. Let us recall that they are precisely defined so that their nondegenerate modules correspond to $(\mathfrak{g},K)$-modules and $(\mathfrak{g}_0,K)$-modules, respectively. Consider the algebraic family of Harish-Chandra pairs $\{(\mathfrak{g}_t,K)\}_{t \in \mathbb{C}}$ corresponding to the deformation to the normal cone of the group inclusion $K\subset G$, see \cite{BHS}*{§2.1.2, §2.3.1}. We refer to this algebraic family as the \textit{algebraic Mackey deformation}. One can associate to it a Hecke $\mathbb{C}[t]$-algebra $\bm{U}$, which specializes at any $t=z \in \mathbb{C}^\times$ to $U$ and to $U^0$ at $t=0$. The nondegenerate $\bm{U}$-modules then have a natural interpretation in terms of algebraic families of Harish-Chandra modules \cite{BHS}*{§2.4} along the Mackey deformation.

The purpose of this paper is to prove that, at the level of the primitive spectrum of $\bm{U}$, the algebraic Mackey deformation is piecewise constant over $\mathbb{C}$. This will have two main consequences: the admissible Mackey bijection is a piecewise homeomorphism with respect to Jacobson topologies; when restricted to $K$-types not exceeding a certain length, $U$ and $U^0$ have the same periodic cyclic homology.

Let us now be more precise. Let $(l_n)_{n \in \mathbb{N}}$ be the collection of all possible lengths of $K$-types, listed in ascending order. For every $n \in \mathbb{N}$, let $\chi_n$ be the projection onto the isotypical components of the $K$-types having length less than $l_n$ and let $\bm{I}_n$ be the ideal of $\bm{U}$ generated by $\chi_n$. We prove the following for every $n\in \mathbb{N}$.

\begin{theoremintro}\label{maintheorem}
    There exists a $\mathbb{C}$-algebra $B_n$ and a morphism of $\mathbb{C}[t]$-algebras
    \[\bm{\pi}_n : \chi_n (\bm{I}_n/\bm{I}_{n-1})\chi_n \longrightarrow  B_n[t]\]
    such that the map $J \mapsto \bm{\pi}^{-1}_n(J)$ is a homeomorphism of the primitive spectrum of $B_n[t]$ onto the primitive spectrum of $\chi_n (\bm{I}_n/\bm{I}_{n-1})\chi_n$, the two spaces being equipped with the Jacobson topology.
\end{theoremintro}

Concretely, the morphism $\bm{\pi}_n$ is built from the representation of $\bm{U}$ on certain algebraic deformations of the principal series associated to the cuspidal data belonging to $\mathcal{C}_n$. The construction of a parabolic induction for families is actually a significant part of the first half of the paper (Subsections \ref{changeofring} and \ref{DeformationPI}). Theorem~\ref{maintheorem} also relies in an essential way on Delorme's work on generalized Harish-Chandra homomorphisms \cites{Delorme84}. Let us also mention that $\pi_n$ is an isomorphism for $G = \mathrm{SL}_2(\mathbb{R})$; we do not know whether this is true in general.

The algebra $B_n$ is well known; it implicitly appears in the afore-mentioned paper of Delorme. We will show that it is Morita equivalent to the commuting direct sum of skew group rings
    \[\bigoplus_{\sigma \in \mathcal{C}_n} R_\sigma \ltimes\mathbb{C}[\mathfrak{a}^*_\sigma /W_\sigma'],\]
where $\mathfrak{a}_\sigma^*$ and $R_\sigma \ltimes W'_\sigma$ denote for each $\sigma \in \mathcal{C}_n$, the complex vector space parameterizing the generalized principal series and the Weyl group associated to $\sigma$, respectively. In its form, our Theorem can thus be compared to the Wassermann Formula \cite{Wassermann}.
    
From Theorem \ref{maintheorem}, we can derive the piecewise bicontinuity of the admissible Mackey bijection. Let us state this properly. For any ring $S$, let us write $\mathrm{Prim}(S)$ for its primitive spectrum, which we endow with the Jacobson topology. If $S$ is a $\mathbb{C}[t]$-algebra having countable dimension over $\mathbb{C}$, the map associating to any $J \in \mathrm{Prim}(S)$ the element $z\in \mathbb{C}$ such that $J\subset (t-z)S$ is well defined and equips $\mathrm{Prim}(S)$ with a bundle structure over $\mathbb{C}$, the latter being endowed with the Zariski topology; the fiber of $\mathrm{Prim}(S)$ over each $z \in \mathbb{C}$ is then the primitive spectrum of the specialization of $S$ at $z$. For any element $S$ of $\{U,U^0, \bm{U}, U[t]\}$, the natural map from the set of isomorphism classes of simple $S$-modules to the primitive spectrum of $S$ turns out to be injective (Proposition \ref{criterion_point_primitive}). In particular, the Mackey bijection can be seen as a bijection $\mathrm{Prim}(U) \to \mathrm{Prim}(U^0)$. For each $S \in \{U,U^0, \bm{U},U[t]\}$ and $n \in \mathbb{N}$, let $\mathrm{Prim}(S)_n$ be the set of annihilators of simple $S$-modules having lowest $K$-types of length $l_n$. As a consequence of Theorem \ref{maintheorem}, we obtain the following for every $n \in \mathbb{N}$.

\begin{corollaryintro}\label{intro_piecewise_homeo}
     The Mackey bijection restricts to a homeomorphism \[\mathrm{Prim}(U)_n \to \mathrm{Prim}(U^0)_n.\]
     In fact, there is an isomorphism $\mathrm{Prim}(U[t])_n \to \mathrm{Prim}(\bm{U})_n$ of continuous bundles over $\mathbb{C}$ (explicitly constructed from $\bm{\pi}_n$) whose restriction to the fibers over $1$ and $0$ are the identity map and the Mackey bijection, respectively.
\end{corollaryintro}

In the above statement, $U[t]$ should be viewed as the Hecke algebra of the constant family of Harish-Chandra pairs over $\mathbb{C}$ with fiber $(\mathfrak{g},K)$. This way, the Mackey deformation is seen to be piecewise spectrally constant. Note that if we restrict to the tempered dual $G$ and instead consider the Fell topology, the Mackey bijection is known to be a piecewise homeomorphism \cite{AfgoustidisClare}*{Theorem B}.

An application of Corollary \ref{intro_piecewise_homeo} is that for any $K$-type $\mu$, the set of primitive ideals of $\bm{U}$ corresponding to simple $\bm{U}$-modules containing $\mu$ as a lowest $K$-type is isomorphic, as a continuous bundle over $\mathbb{C}$, to a complex affine variety of the form $\mathbb{C}\times\mathfrak{a}_\sigma^*/W_\sigma$ for some cuspidal datum $\sigma$.

The second consequence of Theorem \ref{maintheorem} concerns the periodic cyclic homology of corner subalgebras of $\bm{U}$, $U^0$ and $U$. For every $S \in \{\bm{U}, U^0, U\}$, let $S_n = \chi_n S \chi_n$. Building on Baum--Nistor's result \cite{BaumNistor}*{Theorem 7} stating that spectrum preserving morphisms (Definition \ref{definition_spectrum_preserving}) of finite type $\mathbb{C}$-algebras induce isomorphisms in periodic cyclic homology, we prove the following for every $n \in \mathbb{N}$.

\begin{corollaryintro}
    The specializations of $\bm{U}$ at $t=1$ and at $t=0$ induce isomorphisms
    \((\mathrm{sp}_1)_* : \mathrm{HP}_*(\bm{U}_n) \to \mathrm{HP}_*(U_n)\) and \((\mathrm{sp}_0)_* : \mathrm{HP}_*(\bm{U}_n) \to \mathrm{HP}_*(U^0_n),\) respectively. In particular, we have an isomorphism $(\mathrm{sp}_1)_* \circ (\mathrm{sp}_0)_*^{-1} : \mathrm{HP}_*(U^0_n) \to \mathrm{HP}_*(U_n)$.
\end{corollaryintro}

This corollary is of course reminiscent of the Connes-Kasparov isomorphism, see \cite{BaumConnesHigson}*{p.263} and \cites{MAhigson, Afgoustidis3}; see also \cite{Bradd}, where corner subalgebras associated to finite sets of $K$-types appear too. Moreover, let us mention that, in order to apply \cite{BaumNistor}*{Theorem 7}, finite generation conditions are needed. To check that the latter are satisfied, we consider an algebraic family $\mathcal{Z}(\bm{\mathfrak{g}})$ of central subalgebras of $\{U(\mathfrak{g}_t)\}_{t \in \mathbb{C}}$, see Subsection \ref{Center}, and prove that $\bm{U}_n$ is a finitely generated $\mathcal{Z}(\bm{\mathfrak{g}})$-module.

The paper is organized as follows. In Section \ref{sec1}, in the spirit of \cite{BHS}, we introduce Harish-Chandra pairs over arbitrary unital $\mathbb{C}$-algebras, their module categories and their Hecke algebras, which we rather call enveloping algebras. In Section \ref{sec2}, we define the Harish-Chandra pair over $\mathbb{C}[t]$ underlying the algebraic Mackey deformation and construct an algebraic deformation of the parabolic induction. Finally, Section \ref{sec3} is devoted to the proof of the results stated above. The paper also contains two appendices: Appendix \ref{AUA} is about approximately unital algebras and Appendix \ref{relInd} compares various induction processes.

\section{Harish-Chandra pairs and their enveloping algebras}\label{sec1}

A significant part of this paper is about extending usual representation theoretic tools to Harish-Chandra pairs over more general rings than $\mathbb{C}$ and make sure that they still have the same properties. In this first section, we thus introduce Harish-Chandra pairs over arbitrary $\mathbb{C}$-algebras and their module category (Subsection \ref{gKmods}), we define and briefly study their enveloping algebras (Subsection \ref{UgK}) and we discuss change of ring functors (Subsection \ref{changeofring}). In subsection \ref{compactgroups}, we fix some convention about the representation theory of compact Lie groups.

Let us also mention that several approximately unital algebras will be introduced in this section. We refer to Appendix \ref{AUA} for an overview of what is needed for our purpose on this topic. Moreover, for the rest of the paper, we shall use the following convention.

\begin{convention}\label{NotationComplexLie}
    If $G$ is a Lie group, its Lie algebra is denoted by $\mathfrak{g}_\mathbb{R}$ and its complexification by $\mathfrak{g}$. The latter is equipped with the conjugation, which is by definition the unique antilinear map $\mathfrak{g} \to \mathfrak{g}$ fixing every element of $\mathfrak{g}_\mathbb{R}$. In this context, if $\mathfrak{h}$ is a conjugation-invariant $\mathbb{C}$-linear subspace of $\mathfrak{g}$, we write $\mathfrak{h}_\mathbb{R} = \mathfrak{h} \cap \mathfrak{g}_\mathbb{R}$. Conversely, the complexification of every linear subspace $\mathfrak{h}_\mathbb{R}$ of $\mathfrak{g}_\mathbb{R}$ is denoted by $\mathfrak{h}$.
\end{convention}

\subsection{Algebras and representations of compact Lie groups}\label{compactgroups} We start by reviewing the representation theory of compact Lie groups.

Let $K$ denote any compact Lie group. The representations of $K$ are always assumed to be defined over $\mathbb{C}$ and locally $K$ finite, see for instance \cite{KnappVogan}*{p.45}. A morphism of representations of $K$ is by definition a linear $K$-intertwining map. For every pair $(\sigma,\tau)$ of representations of $K$, we denote by $\mathrm{Hom}_K(\sigma,\tau)$ the set of morphisms $\sigma \to \tau$. The Hopf algebra of $K$-finite (for the action by translation) continuous functions on $K$ is denoted by $\mathbb{C}[K]$. We recall that a representation of $K$ is the same as a right $\mathbb{C}[K]$-comodule.

The set of equivalence classes of irreducible representations of $K$ is denoted by $\widehat{K}$ and its elements are called \textbf{$K$-types}. For a given $\mu \in \widehat{K}$ and representation $V$ of $K$, the isotypic component of $V$ of type $\mu$ is denoted by $V^\mu$. If $V^\mu \neq 0$, then $\mu$ is said to be a $K$-type of $V$. Moreover, the quotient of $\mathrm{dim} \,V^\mu$ by the vector space dimension of $\mu$ is called the multiplicity of $\mu$ in $V$. A representation of $K$ is said to be \textbf{admissible} if each of its $K$-type multiplicities is finite.

Since $\mathbb{C}[K]$ is a coalgebra, its linear dual $\mathbb{C}'[K]$ inherits a unital algebra structure. If we view $\mathbb{C}'[K]$ as a space of generalized distributions, then its algebra structure corresponds to the convolution product. For every $k \in K$, the dirac mass at $k$ is an element of $\mathbb{C}'[K]$, so that we can identify $K$ with a subgroup of $\mathbb{C}'[K]^\times$. The union of all finite dimensional ideals of $\mathbb{C}'[K]$ is denoted by $U(K)$, we call it \textbf{the Hecke algebra of $K$}, see \cite{KnappVogan}*{Section I.2}. By construction, $U(K)$ is approximately unital (Appendix \ref{AUA}).

Each representation of $K$ inherits a nondegenerate (Appendix \ref{AUA}) $U(K)$-module structure: if $V$ is a representation of $K$ and $c : V \to  V\otimes_\mathbb{C}\mathbb{C}[K]$ is the associated comodule map, every element $x$ of $\mathbb{C}'[K]\supset U(K)$ acts on $V$ via $(x \otimes \mathrm{id}_V)\circ c$. Conversely, if $V$ is a nondegenerate $U(K)$-module then $V = U(K) \otimes_{U(K)} V$ and the representation of $K$ on $U(K)$ by left multiplication induces a representation of $K$ on $V$. This correspondence between representations of $K$ and nondegenerate $U(K)$-modules is a category equivalence. For more information, see the proof of \cite{KnappVogan}*{Theorem 1.57}.

One can also give a concrete description of $U(K)$. Let $(E_\pi)_{\pi \in \widehat{K}}$ be a family of representatives of $\hatcpt{K}$. The coalgebra decomposition (from the Peter-Weyl theorem)
\[
    \mathbb{C}[K] = \bigoplus_{\pi \in \widehat{K}}\, E_\pi^* \otimes E_\pi
\]
induces by duality an algebra isomorphism
\[ \mathbb{C}'[K] \to \prod_{\pi \in \widehat{K}} \, \mathrm{End}_\mathbb{C}(E_\pi),\qquad x \mapsto (x \curvearrowright E_\pi)_{\pi \in \widehat{K}},\]
which identifies $U(K)$ with $\bigoplus_{\pi \in \widehat{K}} \mathrm{End}_\mathbb{C}(E_\pi)$.

From the above, it follows that for every finite subset $F$ of $\hatcpt{K}$, there is a unique idempotent $\chi_F$ of $U(K)$ whose action on every $E_\pi$ is given by
\begin{equation}\chi_F \curvearrowright E_\pi = \begin{cases}
    \mathrm{id}_{E_\pi} & \text{if $\pi \in F$}\\
    0 & \text{otherwise.}
\end{cases}\end{equation}
If $\mathcal{F}(\hatcpt{K})$ denotes the set of all finite subset of $\hatcpt{K}$, ordered by inclusion, then $(\chi_{F})_{F \in \mathcal{F}(\widehat{K})}$ is an approximate unit (Appendix \ref{AUA}) of $U(K)$.

\subsection{Harish-Chandra pairs and their modules}\label{gKmods} Our goal now is to introduce Harish-Chandra pairs over arbitrary unital commutative $\mathbb{C}$-algebras, as well as their module categories. Let us mention that the objects we define in this subsection are very similar to the algebraic families introduced in \cite{BHS}. For the rest of the subsection, let $A$ denote any unital commutative $\mathbb{C}$-algebra.

\begin{definition}
    A \textbf{Harish-Chandra pair over $A$} is a pair $(\mathfrak{g},K)$, where $\mathfrak{g}$ is a Lie algebra over $A$ and $K$ is a compact Lie group, together with a Lie $\mathbb{C}$-algebra homomorphism $\iota : \mathfrak{k} \to \mathfrak{g}$ and a representation $\kappa$ of $K$ on $\mathfrak{g}$ by Lie $A$-algebra automorphisms such that $\kappa(k) \iota(X) = \iota(\mathrm{Ad}(k)X)$ for every $k \in K$ and $X\in \mathfrak{k}$.
\end{definition}

In the following, the letters $\iota$ and $\kappa$ will always be used in the above context. Occasionally, in order to avoid ambiguity, we will write $\iota = \iota_{\mathfrak{g},K}$ and $\kappa = \kappa_{\mathfrak{g},K}$.

\begin{definition}
    A \textbf{morphism $(\mathfrak{h}, L) \to (\mathfrak{g}, K)$ of Harish-Chandra pairs} over $A$ is defined as a pair $(\phi, \Phi)$, where $\phi : \mathfrak{h} \to \mathfrak{g}$ is a Lie $A$-algebra homomorphism and $\Phi : L \to K$ is a Lie group homomorphism, such that $\phi \circ \iota_{\mathfrak{h}, L} = \iota_{\mathfrak{g}, K}\circ d\Phi$, where $d\Phi$ denotes the differential of $\Phi$, and $\phi \circ\kappa_{\mathfrak{h}, L}(l) = \kappa_{\mathfrak{g}, K}(\Phi(l)) \circ \phi$ for every $l \in L$.
\end{definition}

Let $(\mathfrak{g},K)$ be a Harish-Chandra pair over $A$. If $L$ is a subgroup $K$ and $\mathfrak{h}$ is a Lie $A$-subalgebra of $\mathfrak{g}$ stable by $L$, then $(\mathfrak{h},L)$ defines a Harish-Chandra pair over $A$, equipped with a natural morphism $(\mathfrak{h},L) \to (\mathfrak{g},K)$. Such a Harish-Chandra pair is said to be a \textbf{subpair} of $(\mathfrak{g},K)$.

In the following, we call $\mathfrak{g}$-module any $A$-module $M$ equipped with a Lie $A$-algebra homomorphism $\mathfrak{g} \to \mathrm{End}_A(M)$.

\begin{definition} \label{gKmodule}
     A \textbf{$(\mathfrak{g},K)$-module} is a $\mathfrak{g}$-module $M$ equipped with an $A$-linear representation of $K$ satisfying the following compatibility conditions:
    \begin{enumerate}[label = (\roman*)]
        \item the action of $\mathfrak{k}$ on $M$ via $\iota$ coincides with the differentiation of the representation of $K$,
        \item for every $k \in K, X\in \mathfrak{g}, m\in M$, we have $k(Xm) = (\kappa(k)X)(km)$.
    \end{enumerate}

    A \textbf{morphism of $(\mathfrak{g},K)$-modules} is an $A$-linear map intertwining the actions of $\mathfrak{g}$ and $K$. For every pair $M, N$ of $(\mathfrak{g},K)$-modules, we denote by $\mathrm{Hom}_{\mathfrak{g},K}(M,N)$ the set of such morphisms $M \to N$.
    A $(\mathfrak{g},K)$\textbf{-submodule} of a given $(\mathfrak{g},K)$-module is an $A$-submodule stable by $\mathfrak{g}$ and $K$. A non-zero $(\mathfrak{g},K)$-module is called \textbf{simple} if its only proper $(\mathfrak{g},K)$-submodule is the zero submodule. Finally, a $(\mathfrak{g},K)$-module is said to be \textbf{admissible} if it is admissible as a representation of $K$.
\end{definition}

\subsection{Enveloping algebras of Harish-Chandra pairs}\label{UgK}
The modules over a given Harish-Chandra pair can be characterized as nondegenerate modules over a certain approximately unital algebra, referred to as the enveloping algebra of the considered pair. The aim of this subsection is to construct the latter. Let us fix a commutative unital $\mathbb{C}$-algebra $A$ and and a Harish-Chandra pair $(\mathfrak{g},K)$ over $A$.

Let $U(\mathfrak{g})$ be the universal enveloping $A$-algebra of $\mathfrak{g}$. The representation $\kappa$ of $K$ on $\mathfrak{g}$ extends to a representation of $K$ on $U(\mathfrak{g})$ by $A$-algebra automorphisms. We use the Sweedler notation for the associated comodule map $\kappa' : U(\mathfrak{g}) \to U(\mathfrak{g}) \otimes_\mathbb{C} \mathbb{C}[K] $, that is, we write $\kappa'(X) = \sum X_{(0)} \otimes X_{(1)}$ for all $X \in U(\mathfrak{g})$. Moreover, the Lie algebra homomorphism $\iota : \mathfrak{k} \to \mathfrak{g}$ extends to an algebra homomorphism $U(\mathfrak{k}) \to U(\mathfrak{g})$, which we also denote by $\iota$.

Let us endow $\mathbb{C}'[K]$ with the $\mathbb{C}[K]$-module structure defined by the formula
\begin{equation}\label{module_structure_C'}
    \langle \varphi x, f\rangle = \langle x, \varphi f \rangle \qquad (x \in \mathbb{C}'[K], \varphi, f \in \mathbb{C}[K]).
\end{equation}
In addition, since every element of $\mathbb{C}[K]$ is smooth, by identifying $U(\mathfrak{k})$ with the convolution algebra of distributions on $K$ supported at the identity element, we obtain an inclusion of algebras $U(\mathfrak{k}) \to \mathbb{C}'[K]$.

Now, consider the functor $\bm{\mathrm{F}}_{\mathfrak{g},K}$ associating to any unital $A$-algebra $B$ the set of pairs $(\beta_\mathfrak{g},\beta_K)$ consisting of a unital $A$-algebra morphism $\beta_\mathfrak{g} : U(\mathfrak{g}) \to B$ and a unital $\mathbb{C}$-algebra morphism $\beta_K : \mathbb{C}'[K] \to B$ satisfying $\beta_{K|U(\mathfrak{k})} = \beta_\mathfrak{g} \circ \iota$ and
\[\beta_K(x)\beta_\mathfrak{g}(X) =\sum \beta_\mathfrak{g}(X_{(0)})\beta_K( X_{(1)}x) \]
for every $X \in U(\mathfrak{g})$ and $x \in U(K)$ .

\begin{proposition} \label{representable}
    The functor $\bm{\mathrm{F}}_{\mathfrak{g},K}$ is representable. Moreover, if $(B, (\beta_\mathfrak{g}, \beta_K))$ is a universal element of $\bm{\mathrm{F}}_{\mathfrak{g},K}$, then the bilinear maps
    \[(X,x) \longmapsto \beta_\mathfrak{g}(X) \beta_K(x), \qquad (x,X) \longmapsto \beta_K(x) \beta_\mathfrak{g}(X)\]
    induce $A$-linear isomorphisms $U(\mathfrak{g}) \otimes_{U(\mathfrak{k})}\mathbb{C}'[K] \to B$ and $\mathbb{C}'[K]\otimes_{U(\mathfrak{k})}U(\mathfrak{g})\to B$, respectively.
    \[  \]
\end{proposition}

\begin{proof}
    Define $B = U(\mathfrak{g}) \otimes_{U(\mathfrak{k})}\mathbb{C}'[K]$ as an $A$-module and endow it with the following multiplication law:
    \[ (X \otimes x)(Y \otimes y) = \sum XY_{(0)} \otimes (Y_{(1)}x)y \qquad (X,Y \in U(\mathfrak{g}), x,y \in \mathbb{C}'[K]).\]
    The latter is well defined, associative and $1 \otimes 1$ is the unit of the resulting $A$-algebra. Moreover, the maps $\beta_\mathfrak{g} : X \to X \otimes 1$ and $\beta_K : x \to 1 \otimes x$ are unital morphisms of $A$-algebras and $\mathbb{C}$-algebras, respectively. Once these properties are checked, the fact that $(B, (\beta_\mathfrak{g}, \beta_K))$ is a universal element of $\bm{\mathrm{F}}_{\mathfrak{g},K}$ is clear.

    To conclude, note that the linear map \[\mathbb{C}'[K]\otimes_{U(\mathfrak{k})}U(\mathfrak{g})\to U(\mathfrak{g}) \otimes_{U(\mathfrak{k})}\mathbb{C}'[K], \quad x \otimes X \to \beta_K(x) \beta_\mathfrak{g}(X) = \sum X_{(0)} \otimes X_{(1)}x\]
    is a well defined isomorphism, with inverse $X \otimes x \mapsto \sum S(X_{(1)})x \otimes X_{(0)}$, where $S$ denotes the antipode of the Hopf algebra $\mathbb{C}[K]$.
\end{proof}

From now on, we denote by $(\Upsilon(\mathfrak{g},K), (\upsilon_\mathfrak{g}, \upsilon_K))$ a fixed universal element of $\bm{\mathrm{F}}_{\mathfrak{g},K}$.

\begin{definition}\label{definition_enveloping_algebra}
    The \textbf{enveloping algebra of the pair $(\mathfrak{g},K)$} is defined as the ideal of $\Upsilon(\mathfrak{g},K)$ generated by $\upsilon_K(U(K))$. It is denoted by $U(\mathfrak{g},K)$.
\end{definition}

\begin{proposition}\label{structureU}
    The bilinear maps \[(X,x) \longmapsto \upsilon_\mathfrak{g}(X) \upsilon_K(x), \qquad (x,X) \longmapsto \upsilon_K(x) \upsilon_\mathfrak{g}(X)\] induce $A$-linear isomorphisms \begin{equation} \label{isomorphisms} U(\mathfrak{g}) \otimes_{U(\mathfrak{k})}U(K) \to U(\mathfrak{g},K), \qquad U(K)\otimes_{U(\mathfrak{k})}U(\mathfrak{g})\to U(\mathfrak{g},K).\end{equation}
\end{proposition}

\begin{proof}
    Using the $K$ finiteness of the elements of $\mathbb{C}[K]$, one can easily check that $U(K)$ is a $\mathbb{C}[K]$-submodule of $\mathbb{C}'[K]$, with module structure (\ref{module_structure_C'}). Now, for every $x \in U(K)$ and $X \in U(\mathfrak{g})$, we have
    \begin{gather*}
        \upsilon_K(x)\upsilon_\mathfrak{g}(X) =\sum \upsilon_\mathfrak{g}(X_{(0)})\upsilon_K( X_{(1)}x), \\ \upsilon_\mathfrak{g}(X)\upsilon_K(x) = \sum \upsilon_K( S(X_{(1)})x)\upsilon_\mathfrak{g}(X_{(0)}),
    \end{gather*}
    where $S$ denote the antipode of the Hopf algebra $\mathbb{C}[K]$. Because of the surjectivity of the isomorphisms of Proposition \ref{representable}, this means that $U(\mathfrak{g},K)$ is both the left ideal of $\Upsilon(\mathfrak{g},K)$ generated by $\upsilon_K(U(K))$ and the right ideal of $\Upsilon(\mathfrak{g},K)$ generated by $\upsilon_K(U(K))$. The maps of (\ref{isomorphisms}) are thus surjective.

    Let us prove the injectivity of, say, $U(\mathfrak{g}) \otimes_{U(\mathfrak{k})}U(K) \to U(\mathfrak{g},K)$. If $u$ is an element of $U(\mathfrak{g}) \otimes_{\mathbb{C}}U(K)$ (which we consider as included in $U(\mathfrak{g})\otimes_{\mathbb{C}}\mathbb{C}'[K]$) whose image in $U(\mathfrak{g},K)$ vanishes, then it vanishes in $U(\mathfrak{g}) \otimes_{U(\mathfrak{k})}\mathbb{C}'[K]$ because of the injectivity of the first isomorphism of Proposition \ref{representable}. It follows that $u$ is of the form \[\sum_iX_i \,\iota(Y_i) \otimes x_i - X_i \otimes Y_ix_i,\] for elements $X_i \in U(\mathfrak{g})$, $Y_i \in U(\mathfrak{k})$ and $x_i \in \mathbb{C}'[K]$. Since $u \in U(\mathfrak{g})\otimes_{\mathbb{C}}U(K)$ and $U(K)$ is approximately unital, there is an idempotent $e \in U(K)$ such that $u = ue$ (we consider here the natural right action of $U(K)$ on $U(\mathfrak{g})\otimes_{\mathbb{C}}\mathbb{C}'[K]$). Hence we have \(u = \sum_i X_i \,\iota(Y_i) \otimes x_ie - X_i \otimes Y_ix_ie\),
    which implies that $u$ vanishes in $U(\mathfrak{g}) \otimes_{U(\mathfrak{k})}U(K)$. This proves the injectivity of the map $U(\mathfrak{g}) \otimes_{U(\mathfrak{k})}U(K) \to U(\mathfrak{g},K)$. The proof is the same for the other morphism.
\end{proof}

Let us make some remarks. First, $U(\mathfrak{g},K)$ is approximately unital (Appendix~\ref{AUA}). Indeed, the elements $\upsilon_K(\chi_F)$, for $F$ ranging over $\mathcal{F}(\widehat{K})$, form an approximate unit of $U(\mathfrak{g},K)$. Second, since $U(\mathfrak{g},K)$ is an ideal of $\Upsilon(\mathfrak{g},K)$, it has both a $U(\mathfrak{g})$-bimodule and a $\mathbb{C}'[K]$-bimodule structure via $\upsilon_\mathfrak{g}$ and $\upsilon_K$ respectively. In particular, there are natural left and right actions of both $\mathfrak{g}$ and $K$ on $U(\mathfrak{g},K)$. By construction, the combined actions of $\mathfrak{g}$ and $K$ on the left endow $U(\mathfrak{g},K)$ with a $(\mathfrak{g},K)$-module structure (compare the conditions in the definition of the $(\mathfrak{g},K)$-modules with the conditions satisfied by the pairs $(\beta_\mathfrak{g},\beta_K)$ in the definition of $\bm{\mathrm{F}}_{\mathfrak{g},K}$).

Let us now explain why $(\mathfrak{g},K)$-modules correspond to nondegenerate $U(\mathfrak{g},K)$-modules. Let $M$ be a $(\mathfrak{g},K)$-module. The representation of $K$ on $M$ defines a unital $\mathbb{C}$-algebra morphism $\beta_K : \mathbb{C}'[K] \to \mathrm{End}_A(M)$ (Subsection \ref{compactgroups}), while the action of $\mathfrak{g}$ extends to a unital $A$-algebra morphism $\beta_\mathfrak{g} : U(\mathfrak{g}) \to \mathrm{End}_A(M)$. Because of the compatibility conditions of Definition \ref{gKmodule}, we have $(\beta_\mathfrak{g},\beta_K) \in \bm{\mathrm{F}}_{\mathfrak{g},K}(\mathrm{End}_A(M))$. The latter thus defines a $\Upsilon(\mathfrak{g},K)$-module structure on $M$, which restricts to a nondegenerate $U(\mathfrak{g},K)$-module structure. Conversely, if $M$ is a nondegenerate $U(\mathfrak{g},K)$-module, then $M$ canonically identifies with $U(\mathfrak{g},K) \otimes_{U(\mathfrak{g},K)}M$ and the $(\mathfrak{g},K)$-module structure on $U(\mathfrak{g},K)$ induces a $(\mathfrak{g},K)$-module structure on $M$. One can easily check the following statement.

\begin{proposition}\label{correspondence_Umod_gKmod}
    The two operations defined in the above paragraph are functorial and inverse to each other.
\end{proposition}

Finally, let us discuss two basic examples. Consider the Harish-Chandra pair $(\mathfrak{k},K)$ over $\mathbb{C}$ characterized by $\iota_{\mathfrak{k},K} = \mathrm{id}_\mathfrak{k}$. Then a $(\mathfrak{k},K)$-module is just a representation of $K$. Moreover, the enveloping algebra $U(\mathfrak{k},K)$ identifies with $U(K)$. On the other hand, a $(\mathfrak{g},1)$-module (where we denote by $1$ the trivial group) is just a $\mathfrak{g}$-module, and $U(\mathfrak{g},1) = U(\mathfrak{g})$. For simplicity, we thus adopt the following conventions.

\begin{convention}\label{simpler_notation_HC_pairs}
    For the rest of the paper, the Harish-Chandra pair $(\mathfrak{k},K)$ is labeled by $K$, so that, for example, a $K$-module means a representation of $K$. Similarly, a Harish-Chandra pair of the form $(\mathfrak{g},1)$ is simply denoted by $\mathfrak{g}$.
\end{convention}

\subsection{Change of ring functors}\label{changeofring} To conclude this section, we introduce two different kinds of induction functors and show how they are related.

First, let us consider a commutative unital $\mathbb{C}$-algebra $A$ and a morphism of Harish-Chandra pairs $(\phi, \Phi) : (\mathfrak{h},L) \to (\mathfrak{g},K)$ over $A$. The precomposition by $\Phi$ defines a morphism of Hopf algebras $\Phi' : \mathbb{C}[K] \to \mathbb{C}[L]$. The adjoint map $\Phi'' : \mathbb{C}'[L] \to \mathbb{C}'[K]$ is then a unital $\mathbb{C}$-algebra homomorphism. Moreover, $\phi$ extends to a morphism of unital $A$-algebras $U(\mathfrak{h}) \to U(\mathfrak{g})$, which we again denote by $\phi$. Now the pair $(\upsilon_\mathfrak{g} \circ \phi, \upsilon_K \circ \Phi'')$ is an element of $\bm{\mathrm{F}}_{\mathfrak{h},L}(\Upsilon(\mathfrak{g},K))$, hence it induces a morphism of unital $A$-algebras $\upsilon(\phi, \Phi) : \Upsilon(\mathfrak{h},L) \to \Upsilon(\mathfrak{g},K)$. In particular, multiplication on the right via $\upsilon(\phi, \Phi)$ defines a right $U(\mathfrak{h},L)$-module structure on $U(\mathfrak{g},K)$. \textbf{The $(\phi,\Phi)$-induction} is defined as the functor $U(\mathfrak{g},K) \otimes_{U(\mathfrak{h},L)} -$, from the category of $(\mathfrak{h},L)$-modules to the category of $(\mathfrak{g},K)$-modules. It is denoted by $\mathrm{ind}_{\phi,\Phi}$, or $\mathrm{ind}_{\mathfrak{h},L}^{\mathfrak{g},K}$ if the morphism $(\phi,\Phi)$ is implicit.

Let us also denote by $\mathrm{res}_{\phi,\Phi}$, or $\mathrm{res}_{\mathfrak{h},L}^{\mathfrak{g},K}$, the \textbf{restriction functor} from the category of $(\mathfrak{g},K)$-modules to that of $(\mathfrak{h},L)$-modules, which consists in preserving the underlying $A$-module and taking the pullback of the action of $(\mathfrak{g},K)$ by $(\phi,\Phi)$.

Second, let us consider a morphism of unital $\mathbb{C}$-algebras $\psi : A \to B$. Let $(\mathfrak{g},K)$ denote any Harish-Chandra pair over $A$. The \textbf{specialization of $(\mathfrak{g},K)$ at $\psi$} is defined as the Harish-Chandra pair $(\mathfrak{g}_\psi,K) = (B \otimes_A\mathfrak{g}, K)$ over $B$ characterized by $\iota_{\mathfrak{g}_\psi,K} = \mathrm{id}_B \otimes \iota_{\mathfrak{g},K}$ and $\kappa_{\mathfrak{g}_\psi,K}(k) = \mathrm{id}_B \otimes \kappa_{\mathfrak{g},K}(k)$ for every $k \in K$. Since $U(\mathfrak{g}_\psi) = B \otimes_A U(\mathfrak{g})$, we have $\Upsilon(\mathfrak{g}_\psi,K) = B \otimes_A \Upsilon(\mathfrak{g},K)$ and thus \begin{equation}\label{BtensorA}
    U(\mathfrak{g}_\psi,K) = B \otimes_A U(\mathfrak{g},K)
\end{equation} by proposition \ref{structureU}. The specialization at $\psi$ defines a functor from the category of Harish-Chandra pairs over $A$ to the category of Harish-Chandra pairs over $B$.

Moreover, if $M$ is a $(\mathfrak{g},K)$-module then $M_\psi = B \otimes_A M$ has a natural $(\mathfrak{g}_\psi, K)$-module structure; we call it the \textbf{specialization of $M$ at $\psi$}. This defines a functor $\mathrm{sp}_\psi^{\mathfrak{g},K}$ from the category of $(\mathfrak{g},K)$-modules to the category of $(\mathfrak{g}_\psi,K)$-modules. In what follows, the \textbf{specialization map} $(M \to M_\psi, \, m \mapsto 1 \otimes m)$ will generically be denoted by $\mathrm{sp}_\psi$.

Conversely, if $M$ is a $(\mathfrak{g}_\psi, K)$-module, then it restricts to a $(\mathfrak{g}, K)$-module $M^\psi$. Indeed, a $(\mathfrak{g}_\psi, K)$-module is just a a $(\mathfrak{g}, K)$-module endowed with a $B$-module structure such that the action of $A$ factors through $\psi$ and the action of $B$ commutes with that of $\mathfrak{g}$ and $K$. The \textbf{forgetful functor} $M \mapsto M^\psi$ from the category of $(\mathfrak{g}_\psi, K)$-modules to that of $(\mathfrak{g}, K)$-modules is denoted by $\mathrm{for}_\psi^{\mathfrak{g},K}$.

We warn the reader that the words restriction and specialization may very well be used in contexts where the morphism $(\mathfrak{h},L) \to (\mathfrak{g},K)$ is not injective and when $\psi : A \to B$ is not surjective, respectively.

The above functors satisfy the following compatibility properties.

\begin{proposition}\label{properties_ind_res} Let $\psi : A \to B$ be as above and let $(\mathfrak{h},L) \to (\mathfrak{g},K)$ be a morphism of Harish-Chandra pairs over $A$.
    The specialization and forgetful functors commute with induction, that is, we have natural isomorphisms
        \begin{gather}
            \label{commute_sp_ind}\mathrm{sp}_\psi^{\mathfrak{g},K} \circ \mathrm{ind}_{\mathfrak{h},L}^{\mathfrak{g},K} = \mathrm{ind}_{\mathfrak{h}_\psi,L}^{\mathfrak{g}_\psi,K} \circ \mathrm{sp}_\psi^{\mathfrak{h},L},\\
            \label{commute_for_ind}\mathrm{for}_\psi^{\mathfrak{g},K} \circ \mathrm{ind}_{\mathfrak{h}_\psi,L}^{\mathfrak{g}_\psi,K} = \mathrm{ind}_{\mathfrak{h},L}^{\mathfrak{g},K} \circ \mathrm{for}_\psi^{\mathfrak{h},L}.
        \end{gather}
    
    Moreover, if $(\mathfrak{i},M) \to (\mathfrak{h},L)$ is another morphism of Harish-Chandra pairs over $A$, then the functor $\mathrm{ind}_{\mathfrak{h},L}^{\mathfrak{g},K} \circ \mathrm{ind}_{\mathfrak{i},M}^{\mathfrak{h},L}$ is naturally isomorphic to $\mathrm{ind}_{\mathfrak{i},M}^{\mathfrak{g},K}$.
\end{proposition}

\begin{proof}
    These statements follow from the associativity of the tensor product. For (\ref{commute_sp_ind}) and (\ref{commute_for_ind}), use (\ref{BtensorA}).
\end{proof}

Let us now introduce some notations concerning a particular type of specializations. Let $A$ and $B$ be two commutative unital $\mathbb{C}$-algebras.

\begin{convention}\label{notation_specialization_HC}
    Let $\mu_{A,B}$ be the canonical embedding $A \to B\otimes_\mathbb{C} A$. When denoting the specialization and forgetful functors related to the morphism to $\mu_{A,B}$, we will replace the subscript $\mu_{A,B}$ with a subscript $B$. If $(\mathfrak{g},K)$ is a Harish-Chandra pair over $A$, its specialization at $\mu_{A,B}$ is simply denoted by $(\mathfrak{g}_B, K)$, or $K_B$ if $A = \mathbb{C}$ and $(\mathfrak{g},K)$ is the Harish-Chandra pair defined by $K$ (Convention~\ref{simpler_notation_HC_pairs}). Similarly, the specialization of a $(\mathfrak{g},K)$-module $M$ at $\mu_{A,B}$ is denoted by $M_B$ and the associated specialization map $M \to M_B$ by $\mathrm{sp}_B$.
\end{convention}

\begin{remark}\label{compatible}
    We stress that a $(\mathfrak{g}_B, K)$-module is just a $(\mathfrak{g},K)$-module endowed with a \textbf{compatible} $B$-module structure, where by compatible we mean that the $\mathbb{C}$-vector space structure of $M$ inherited from the $A$-module and $B$-module structures are the same and that the action of $B$ commutes with the actions of $A$, $K$ and~$\mathfrak{g}$.
\end{remark}

Finally, the following lemma gives sufficient conditions in order to have a ``compact picture" for induced modules.

\begin{lemma}\label{lemma_compact_pic}
    Let $(\mathfrak{g},K)$ be a Harish-Chandra pair over $A$ and let $(\mathfrak{h},L)$ be a subpair of $(\mathfrak{g},K)$. Assume that the Lie algebra homomorphism $\iota :\mathfrak{k} \to \mathfrak{g}$ is injective, so that we may identify $\mathfrak{k}$ with a Lie $\mathbb{C}$-subalgebra of $\mathfrak{g}$. Assume also that there exists a free $A$-submodule $\mathfrak{s}$ of $\mathfrak{g}$ such that $\mathfrak{g} = \mathfrak{k}_A \oplus \mathfrak{s}$ and $\mathfrak{h} = \mathfrak{l}_A \oplus \mathfrak{s}$.
    
    Let $M$ be an $(\mathfrak{h},L)$-module. Assume that there is an $L$-module $V$ and an $A$-module $N$ such that the restriction of $M$ to an $L_A$-module identifies with $V \otimes_\mathbb{C}N$ via an isomorphism of $L_A$-modules $\rho : V \otimes_\mathbb{C}N \to M$.
    
    Then the restriction of $\mathrm{ind}_{\mathfrak{h},L}^{\mathfrak{g},K}(M)$ to a $K_A$-module identifies with $\mathrm{ind}_L^K(V) \otimes_\mathbb{C} N$ through
    \[\begin{array}{cccll}
         U(K) \otimes_{U(L)}V \otimes_\mathbb{C} N& \longrightarrow &U(\mathfrak{g},K) \otimes_{U(\mathfrak{h},L)}M  \\
         x \otimes v \otimes w& \longmapsto & \upsilon_K(x) \otimes \rho(v \otimes w), 
    \end{array}\]
    where $\upsilon_K$ denotes the canonical morphism $U(K) \to U(\mathfrak{g},K)$, see Definition \ref{definition_enveloping_algebra}.
\end{lemma}

\begin{proof}
    Let $(Y_\lambda)_{\lambda \in \Lambda}$ be a totally ordered basis of $\mathfrak{s}$. Let us denote by $U(\mathfrak{s})$ the $A$-submodule of $U(\mathfrak{g})$ spanned by the elements of the form
    \[Y^{\nu_1}_{\lambda_1}\dotsb Y^{\nu_n}_{\lambda_n}\]
    for $(\lambda_1, ...,\lambda_n)$ an increasing sequence of elements of $\Lambda$ and $(\nu_1, ..., \nu_n)$ a sequence of natural numbers. By the Poincaré-Birkhoff-Witt theorem \cite{GrAlLie1}*{§2 N°7}, the multiplication in $U(\mathfrak{g})$ yields an $A$-linear isomorphism
    \[U(\mathfrak{k}) \otimes_\mathbb{C}U(\mathfrak{s}) = U(\mathfrak{k}_A) \otimes_AU(\mathfrak{s}) \to U(\mathfrak{g}).\]
    Using (\ref{isomorphisms}), we obtain an isomorphism
    \[U(K) \otimes_\mathbb{C}U(\mathfrak{s}) \to U(\mathfrak{g},K),\]
    induced by multiplication in $\Upsilon(\mathfrak{g},K)$. Since $U(K)$ is a nondegenerate right $U(L)$-module for multiplication, we obtain from the above an isomorphism
    \begin{equation}\label{isomorphism_UgK}
        U(K)\otimes_{U(L)} U(L)\otimes_\mathbb{C} U(\mathfrak{s}) \to U(\mathfrak{g},K).
    \end{equation}

    Now, $\mathfrak{h}$ being a direct summand of $\mathfrak{g}$ (at the level of $A$-modules), the natural $A$-algebra homomorphism $U(\mathfrak{h}) \to U(\mathfrak{g})$ is injective (this is due to the Poincaré-Birkhoff-Witt theorem). Hence we can view $U(\mathfrak{h})$ as a subalgebra of $U(\mathfrak{g})$ and $U(\mathfrak{s})$ as a subalgebra of $U(\mathfrak{h})$. As above, the multiplication in $\Upsilon(\mathfrak{h},L)$ yields an isomorphism
    \begin{equation}\label{isomorphism_UhL}
        U(L) \otimes_\mathbb{C} U(\mathfrak{s}) \to U(\mathfrak{h},L).\
    \end{equation}
    Thus, by combining (\ref{isomorphism_UgK}) and (\ref{isomorphism_UhL}), we obtain that the multiplication in $\Upsilon(\mathfrak{g},K)$, together with the morphism of algebras $U(\mathfrak{h},L) \to \Upsilon(\mathfrak{g},K)$ induced by the morphism $(\mathfrak{h},L) \to (\mathfrak{g},K)$, induces an isomorphism
    \[U(K) \otimes_{U(L)} U(\mathfrak{h},L) \to U(\mathfrak{g},K).\]
    The latter implies a $K_A$-module isomorphism
    \begin{align*}
        \mathrm{ind}_{\mathfrak{h},L}^{\mathfrak{g},K}(M) = U(\mathfrak{g},K) \otimes_{U(\mathfrak{h},L)}M \approx U(K) \otimes_{U(L)}M.
    \end{align*}
    Now $\rho$ identifies the $L_A$-module $M$ to $V \otimes_\mathbb{C}N$, so we obtain the desired isomorphism of $K_A$-modules $\mathrm{ind}_{\mathfrak{h},L}^{\mathfrak{g},K}(M) \approx \mathrm{ind}_L^K(V) \otimes_\mathbb{C} N$.
\end{proof}

\section{The algebraic Mackey deformation and parabolic induction}\label{sec2}

Let $G$ be a real reductive group and let $K$ be a maximal compact subgroup of $G$. At the level of Lie groups, the central object underlying the Mackey deformation is the Lie groupoid $\mathcal{G} \rightrightarrows \mathbb{R}$ defined by the deformation to the normal cone relative to the inclusion $K\subset G$, which in this case is simply a smooth bundle of Lie groups. From the point of view of Harish-Chandra pairs though, the right analog of $\mathcal{G}$ is a certain Harish-Chandra pair $(\bm{\mathfrak{g}},K)$ over $\mathbb{C}[t]$, which corresponds to the kind of algebraic families of Harish-Chandra pairs introduced in \cite{BHS}*{§2.1.2, §2.3.1}. The goal of this section is, after some conventions (Subsection \ref{DeformationConvention}), to introduce this Harish-Chandra pair (Subsection \ref{DeformationHCpairs}), then to define an algebraic notion of parabolic induction along the Mackey deformation (Subsection \ref{DeformationPI}) and finally to study a certain central subalgebra of $U(\bm{\mathfrak{g}})$ (Subsection \ref{Center}), which we will use in the last section to derive our main results.

\subsection{Conventions on real reductive groups}\label{DeformationConvention}

Let us first introduce some notations and conventions regarding real reductive groups, lowest $K$-types and motions groups. In the following, $G$ denotes any real reductive group. Recall Convention~\ref{NotationComplexLie}. 

Let $\theta$ be a Cartan involution of $G$ and let $K$ denote the subgroup of fixed points of $\theta$. This Cartan involution induces, by differentiation and complexification, an involution of $\mathfrak{g}$ that we also denote by $\theta$; we write $\mathfrak{s} = \mathfrak{g}^{-\theta}$. We fix a bilinear form $(\cdot |\cdot)$ on $\mathfrak{g}$ which is invariant by $\theta$ and $\mathrm{Ad}_{\mathfrak{g}}(G)$ and such that $(\theta(\cdot) |\cdot)$ is positive definite (in particular real) on $\mathfrak{g}_\mathbb{R} \times \mathfrak{g}_\mathbb{R}$. The bilinear form on $\mathfrak{g}^*$ obtained by duality is also denoted by $(\cdot |\cdot)$. If $\mathfrak{h}$ is a conjugation-invariant linear subspace of $\mathfrak{s}$, we identify its linear dual $\mathfrak{h}^*$ with the space of linear forms on $\mathfrak{g}$ vanishing on $\mathfrak{h}^\perp$. The centralizer of $G$ in $\exp(\mathfrak{s}_\mathbb{R})$ is called the \textbf{split component} of $G$, while the kernel of all Lie group homomorphisms $G \to \mathbb{R}_+^*$ is called the \textbf{essential component} of $G$.

We denote by $\mathfrak{a}$ a fixed maximal conjugation-invariant abelian subspace of $\mathfrak{s}$. The corresponding Weyl group is denoted by $W(\mathfrak{g},\mathfrak{a})$ and the analytic subgroup of $G$ with Lie algebra $\mathfrak{a}_\mathbb{R}$ is denoted by $A$. Let $\Delta$ be the restricted root system associated to the pair $(\mathfrak{g}_\mathbb{R},\mathfrak{a}_\mathbb{R})$. For any $\alpha \in \Delta\cup\{0\}$, let us write $\mathfrak{g}_\mathbb{R}^\alpha$ for the associated root space. To each positive root system $\Delta^+ \subset \Delta$ we associate the parabolic subalgebra $\mathfrak{p}(\Delta^+)_\mathbb{R} = \sum_{\alpha \in \Delta^+\cup \{0\}} \mathfrak{g}_\mathbb{R}^\alpha$ of $\mathfrak{g}_\mathbb{R}$. In the following, we only consider parabolic subgroups of $G$ obtained as normalizers of subalgebras of $\mathfrak{g}_\mathbb{R}$ containing some $\mathfrak{p}(\Delta^+)_\mathbb{R}$, for a positive system $\Delta^+ \subset \Delta$. If $P$ is a parabolic subgroup of $G$, we refer to the subgroup $P \cap \theta(P)$ as \textit{the} Levi factor of $P$. In general, by \textbf{Levi subgroup} of $G$, we mean any subgroup of $G$ that can be obtained as the Levi factor of a parabolic subgroup. Alternatively, the Levi subgroups of $G$ are the centralizers of subsets of $\mathfrak{a}$ (or, equivalently, of subsets of $\mathfrak{a}^*$). This characterization is standard; a reader may find the material for the proof in \cite{Knapp}.

Let $L$ be a Levi subgroup of $G$ and let $A_L$ be the split component of $L$. The quotient of $L$ by its essential component defines a morphism of Harish-Chandra pairs $(\mathfrak{l}, K\cap L) \to \mathfrak{a}_L$. If $M$ is an $(\mathfrak{l}, K\cap L)$-module and $N$ is an $\mathfrak{a}_L$-module, we write
\begin{equation}\label{tensor_prod_module_1} M \otimes N = M \otimes_\mathbb{C} \mathrm{res}_{\mathfrak{l}, K\cap L}^{\mathfrak{a}_L}(N),\end{equation}
where the latter is considered as a tensor product $(\mathfrak{l}, K\cap L)$-module. The one-dimensional $\mathfrak{a}_L$-module defined by any $\lambda \in \mathfrak{a}_L^*$ is denoted by $\mathbb{C}_\lambda$. For any parabolic subgroup $P$ of $G$ having $L$ as Levi factor, we denote by $\rho_P$ the linear functional on $\mathfrak{a}_L$ defined by the half trace of the adjoint action on $\mathfrak{p}$.

Let us also recall Vogan's notion of lowest $K$-types \citelist{\cite{VoganAnnals}*{Definition 5.1} \cite{VoganGreen}*{Definition 5.4.18}}. Let $\pi$ be any $K$-type. Given a Cartan subalgebra $\mathfrak{t}_\mathbb{R} \subset \mathfrak{k}_\mathbb{R}$, a choice of positive root system $\Delta^+_c \subset i\mathfrak{t}_\mathbb{R}^*$ associated to the pair $(\mathfrak{k}_\mathbb{R}, \mathfrak{t}_\mathbb{R})$, and a highest weight $\mu$ of $\pi$, we define the \textbf{length} of $\pi$ as follows:
\[\mathrm{length}(\pi) = (\mu + 2\rho_c | \mu + 2\rho_c),\]
where $\rho_c$ denote the half sum of all the elements of $\Delta^+_c$. This definition does not actually depend on the choices of $\mathfrak{t}_\mathbb{R}$, $\Delta^+_c$ and $\mu$. By definition, a \textbf{lowest $K$-type} of a given representation $\sigma$ of $K$ is a $K$-type of $\sigma$ having minimal length among all $K$-types of $\sigma$. Since the number of $K$-types of length less than a fixed $r >0 $ is finite, the set of lowest $K$-types of a representation of $K$ is finite.

\begin{remark}
    In his green book \cite{VoganGreen}, from which we will extract several results, Vogan rather uses the notion of lambda-lowest $K$-type (Definition 5.4.1), which is more complicated. Nonetheless, the lambda-lowest $K$-types coincide with the lowest $K$-types for the standard $(\mathfrak{g},K)$-modules (the ones obtained from discrete series by parabolic induction), see \cite{VoganAnnals}*{Lemma 8.8}. We also refer to the discussion in \cite{VoganGreen}*{p.297}.
\end{remark}

We recall that the \textbf{Cartan motion group} of $G$, denoted by $G_0$, is the semi-direct product $K \ltimes \mathfrak{s}_\mathbb{R}$ for which the action of $K$ on $\mathfrak{s}_\mathbb{R}$ by conjugation coincides with the restriction of the adjoint action of $G$. As for $G$, the natural inclusion $K \subset G_0$ defines a Harish-Chandra pair $(\mathfrak{g}_0, K)$.

Let $P$ be a parabolic subgroup of $G$, let $L$ be its Levi factor and let $A_L$ be the split component of $L$. The Cartan motion group $L_0$ of $L$ is identified with $(K\cap L) \ltimes (\mathfrak{s}\cap \mathfrak{l})_\mathbb{R} \subset G_0$. We also write $P_0 = (K\cap P) \ltimes \mathfrak{s}_\mathbb{R}$. The quotients of $P_0$ by $\mathfrak{s}_\mathbb{R}$ and $(K\cap P) \ltimes \mathfrak{a}_{L,\mathbb{R}}^\perp$ define morphisms of Harish-Chandra pairs $(\mathfrak{p}_0,K\cap P) \to K\cap P$ and $(\mathfrak{p}_0,K\cap P) \to \mathfrak{a}_L$, respectively. If $V$ is a $K\cap P$-module and $M$ is an $\mathfrak{a}_L$-module, we define
\begin{equation}\label{tensor_prod_module_2}
    V \otimes M = \mathrm{res}^{K\cap P}_{\mathfrak{p}_0,K\cap P}(V) \otimes_\mathbb{C} \mathrm{res}_{\mathfrak{p}_0,K\cap P}^{\mathfrak{a}_L}(M),
\end{equation}
the latter being considered as a tensor product $(\mathfrak{p}_0, K\cap P)$-module.

\subsection{The deformation at the level of Harish-Chandra pairs}\label{DeformationHCpairs} Let us now introduce the Harish-Chandra pair $(\bm{\mathfrak{g}},K)$ as well as some of its subpairs attached to parabolic subgroups of $G$. For every $z \in \mathbb{C}$, the character of $\mathbb{C}[t]$ defined by evaluation at $z$ is itself labelled by $z$.

\begin{definition}
    Let $\bm{\mathfrak{g}}$ be the Lie $\mathbb{C}[t]$-subalgebra $\mathfrak{k}[t] \oplus t\mathfrak{s}[t]$ of $\mathfrak{g}[t]$. Consider the action of $K$ on $\mathfrak{g}[t]$ obtained from the adjoint action on $\mathfrak{g}$ by extension of scalars. This action restricts to an action on $\bm{\mathfrak{g}}$ by Lie $\mathbb{C}[t]$-algebra automorphisms. The Harish-Chandra pair $(\bm{\mathfrak{g}},K)$ defined by this action and the embedding $\mathfrak{k} \hookrightarrow \mathfrak{k}[t] \subset \bm{\mathfrak{g}}$ is called the \textbf{algebraic Mackey deformation pair}.
\end{definition}

Let us see how this pair specializes at points of $\mathbb{C}$. For every $z \in \mathbb{C}^\times$, we have $\bm{\mathfrak{g}}_z = \mathfrak{g}$ (notation defined in Subsection \ref{changeofring}), the specialization map $\mathrm{sp}_z : \bm{\mathfrak{g}}\to \mathfrak{g}$ being defined by $\mathrm{sp}_z (X) = X$ and $\mathrm{sp}_z( tY) = zY$, for every $X \in \mathfrak{k}$ and $Y \in \mathfrak{s}$. On the other hand, the specialization of $\bm{\mathfrak{g}}$ at $0$ is $\mathfrak{g}_0 = \mathfrak{k} \ltimes \mathfrak{s}$, with specialization map defined by $\mathrm{sp}_0(X) = X$ and $\mathrm{sp}_0(tY) = Y$, for every $X \in \mathfrak{k}$ and $Y \in \mathfrak{s}$. Then, from Formula (\ref{BtensorA}), it follows that $U(\bm{\mathfrak{g}},K)$ specializes at each $z \in \mathbb{C}^\times$ to $U(\mathfrak{g},K)$ and at $0$ to $U(\mathfrak{g}_0,K)$.

Similarly, we wish to define a Harish-Chandra pair $(\bm{\mathfrak{p}},K\cap P)$ over $\mathbb{C}[t]$ for every parabolic subgroup $P$ of $G$. Let us fix such a subgroup $P$, let $L$ be its Levi factor and let $\mathfrak{u}$ be the nilpotent radical of $\mathfrak{p}$. Let $(\bm{\mathfrak{l}},K\cap L)$ be the algebraic Mackey deformation pair associated to $L$, which we view as a subpair of $(\bm{\mathfrak{g}},K)$. Note that $t\mathfrak{g}[t]$, hence $t\mathfrak{u}[t]$, are Lie $\mathbb{C}[t]$-subalgebras of $\bm{\mathfrak{g}}$.

\begin{definition}
    The Lie $\mathbb{C}[t]$-subalgebra $\bm{\mathfrak{l}} \oplus t \mathfrak{u}[t]$ of $\bm{\mathfrak{g}}$ is denoted by $\bm{\mathfrak{p}}$. The Mackey deformation pair associated to $P$ is defined as the subpair $(\bm{\mathfrak{p}},K\cap P)$ of $(\bm{\mathfrak{g}},K)$.
\end{definition}

Let us show how these Harish-Chandra pairs specialize at points of $\mathbb{C}$. For the pair $(\bm{\mathfrak{l}},K\cap L)$, it is the same as for $(\bm{\mathfrak{g}},K)$: we have $(\bm{\mathfrak{l}}_z,K\cap L) = (\mathfrak{l},K\cap L)$ if $z\in \mathbb{C}^\times$ and $(\bm{\mathfrak{l}}_0,K\cap L) = (\mathfrak{l}_0,K\cap L)$, where $\mathfrak{l}_0 = (\mathfrak{k}\cap \mathfrak{l}) \ltimes (\mathfrak{s}\cap \mathfrak{l})$ is the complexification of the Lie algebra of the motion group $L_0$ of $L$.

The case of the pair $(\bm{\mathfrak{p}}, K\cap P)$ is slightly more complicated. Let us choose a positive root system associated to the pair $(\mathfrak{g},\mathfrak{a})$ so that the span of the positive root spaces $\mathfrak{n}$ contains $\mathfrak{u}$. We have $\bm{\mathfrak{g}} = \mathfrak{k}[t] \oplus t\mathfrak{a}[t] \oplus t\mathfrak{n}[t]$ and $\bm{\mathfrak{p}} = (\mathfrak{l}\cap \mathfrak{k})[t] \oplus t\mathfrak{a}[t] \oplus t\mathfrak{n}[t]$. It follows that $\bm{\mathfrak{p}}$ is a direct summand of $\bm{\mathfrak{g}}$. Hence, for every $z \in \mathbb{C}$, the specialization of $\bm{\mathfrak{p}}$ at $z$ identifies with the image of $\bm{\mathfrak{p}}$ in $\bm{\mathfrak{g}}_z$. We thus obtain the following.

\begin{proposition}\label{specialization_P_at_0}
    For every $z \in \mathbb{C}^\times$, we have $(\bm{\mathfrak{p}}_z, K\cap P) = (\mathfrak{p}, K\cap P)$, the specialization map being the restriction of $\mathrm{sp}_z : \bm{\mathfrak{g}} \to \mathfrak{g}$ to $\bm{\mathfrak{p}}$. Moreover, the specialization of $(\bm{\mathfrak{p}},K\cap P)$ at $0$ is $(\mathfrak{p}_0, K\cap P)$, where the specialization map $\mathrm{sp}_0 : \bm{\mathfrak{p}} \to \mathfrak{p}_0$ is defined by
\[\mathrm{sp}_0(Y_\mathfrak{k} + t Y_\mathfrak{s} +tZ) = Y_\mathfrak{k} \ltimes (Y_\mathfrak{s} + Z_\mathfrak{s})\qquad (Y \in \mathfrak{l}, Z \in \mathfrak{u}),\]
where for every $X \in \mathfrak{g}$, we write $X_\mathfrak{k} = \tfrac{1}{2} (Z+\theta(Z))$ and $X_\mathfrak{s} = \tfrac{1}{2} (Z-\theta(Z))$.
\end{proposition}

\begin{proof}
    Let us only check what concerns the specialization at $0$. Consider the specialization map $\mathrm{sp}_0 : \bm{\mathfrak{g}} \to \mathfrak{g}_0$. We saw above that the specialization of  $(\bm{\mathfrak{p}},K\cap P)$ at $0$ is $(\mathrm{sp}_0(\bm{\mathfrak{p}}),K\cap P)$, with specialization map $\mathrm{sp}_{0|\bm{\mathfrak{p}}} : \bm{\mathfrak{p}} \to \mathrm{sp}_0(\bm{\mathfrak{p}})$. Now, if $Y \in \mathfrak{l}$ and $Z \in \mathfrak{u}$, we have
    \[\mathrm{sp}_0(Y_\mathfrak{k} + tY_\mathfrak{s} +tZ) = \mathrm{sp}_0(Y_\mathfrak{k} + tZ_\mathfrak{k} +t (Y_\mathfrak{s} + Z_\mathfrak{s})) = Y_\mathfrak{k} \ltimes (Y_\mathfrak{s} + Z_\mathfrak{s}).\]
    This proves the formula in the proposition but also the equality $\mathrm{sp}_0(\bm{\mathfrak{p}}) = \mathfrak{p}_0$, since $\mathfrak{s} = (\mathfrak{s}\cap \mathfrak{l}) \oplus \{Z_\mathfrak{s} : Z \in \mathfrak{u}\}$.
\end{proof}

\subsection{Algebraic deformation of the parabolic induction}\label{DeformationPI}
Let us now define the notion of parabolic induction suited for our generalized Harish-Chandra pairs. Let $P$ be a parabolic subgroup of $G$ with Levi factor $L$ and let $M_L$ and $A_L$ be the essential component and the split component of $L$, respectively. Let us fix a unital commutative $\mathbb{C}$-algebra $S$.

Let us begin with some conventions. Recall the tensor product notations (\ref{tensor_prod_module_1}) and (\ref{tensor_prod_module_2}) for modules over $(\mathfrak{l},K\cap L)$ or $(\mathfrak{p}_0,K\cap P)$, as well as Convention~\ref{notation_specialization_HC} and Remark \ref{compatible}. If $M$ is an $(\mathfrak{l}_S,K\cap L)$-module and $N$ is an $\mathfrak{a}_L$-module, the $(\mathfrak{l},K\cap L)$-module $M \otimes N$ inherits from $M$ a compatible $S$-module structure and can therefore be considered as an $(\mathfrak{l}_S,K\cap L)$-module. In the same way, if $M$ is now only an $(\mathfrak{l},K\cap L)$-module but $N$ is an $\mathfrak{a}_{L,S}$-module, we view $M\otimes N$ as an $(\mathfrak{l}_S,K\cap L)$-module. Finally, if $V$ is a $K\cap P$-module and $N$ is an $\mathfrak{a}_{L,S}$-module, we view $V \otimes N$ as a $(\mathfrak{p}_{0,S}, K\cap P)$-module.

In the following, we will also often consider the Harish-Chandra pairs $(\bm{\mathfrak{g}}_S,K)$, $(\bm{\mathfrak{l}}_S,K\cap L)$ and $(\bm{\mathfrak{p}}_S, K\cap P)$. Regarding the specializations of these Harish-Chandra pairs, we will use the following notations. Every element $z \in \mathbb{C}$ defines an algebra homomorphism $S[t] \to S$, which we also denote by $z$. Similarly, every character $\lambda$ of $S$ defines an algebra homomorphism $S[t] \to \mathbb{C}[t]$ which we also denote by $\lambda$.

Let us now revisit the classical parabolic induction. We define $\mathrm{ext}_{L,S}^P$, the functor ``of extension to the parabolic subgroup", from the category of $(\mathfrak{l}_S,K\cap L)$-modules to the category of $(\mathfrak{p}_S, K)$-modules, by the formula
\begin{equation}\label{first_ext} \mathrm{ext}_{L,S}^P(-) = \mathrm{res}^{\mathfrak{l}_S, K\cap L}_{\mathfrak{p}_S, K\cap P}(- \otimes \mathbb{C}_{-\rho_P}),\end{equation}
where the morphism $(\mathfrak{p}_S, K\cap P) \to (\mathfrak{l}_S, K\cap L)$ derives from the Lie group homomorphism $P \to L$ consisting in the quotient of $P$ by its unipotent radical.
The \textbf{parabolic induction functor}, from the category of $(\mathfrak{l}_S,L\cap K)$-modules to the category of $(\mathfrak{g}_S,K)$-modules, is defined as follows: \begin{equation}\label{first_parabolic_induction}
    \mathrm{Ind}_{P,S}^G = \mathrm{ind}_{\mathfrak{p}_S,K\cap P}^{\mathfrak{g}_S,K} \circ \mathrm{ext}_{L,S}^P.
\end{equation}
If $S = \mathbb{C}$, we simply write $\mathrm{ext}_L^P$ and $\mathrm{Ind}_L^P$. The latter is precisely the algebraic counterpart of the usual parabolic induction, see Theorem~\ref{comparisonInd}.

In order to harmonize the induction notation for the Cartan motion group with that of parabolic induction, we will write
\begin{equation}\mathrm{Ind}_{P_0,S}^{G_0} \quad \text{and}\quad \mathbf{Ind}_{P_0,S}^{G_0}\end{equation}
for $\mathrm{ind}_{\mathfrak{p}_{0,S},K\cap P}^{\mathfrak{g}_{0,S},K}$ and $\mathrm{ind}_{\bm{\mathfrak{p}}_{0,S},K\cap P}^{\bm{\mathfrak{g}}_{0,S},K}$, respectively. We simply write $\mathrm{Ind}_{P_0}^{G_0}$ and $\mathbf{Ind}_{P_0}^{G_0}$ in the case where $S = \mathbb{C}$. We stress that $\mathrm{Ind}_{P_0}^{G_0}$ coincides, at the level of $K$-finite vectors, with the topological induction usually considered in the literature, from representations of $P_0$ to representations of $G_0$. This is due to the fact that $P_0$ acts trivially on the space of complex \(-\tfrac{1}{2}\) densities of $\mathfrak{g}/\mathfrak{p}$, see Appendix \ref{relInd}.

Our goal now is to generalize (\ref{first_ext}) and (\ref{first_parabolic_induction}), and define a parabolic induction from $(\mathfrak{l}_S,K\cap L)$-modules to $(\bm{\mathfrak{g}}_S, K)$-modules. Let $V$ be an $(\mathfrak{l}_S,K\cap L)$-module and let $\pi$ denote the associated representation of $\Upsilon(\mathfrak{l}_S, K\cap L)$ on $V$, see Subsection \ref{UgK}. Let us consider $\mathrm{End}_S(V)$ as the set  of $0$-degree elements of $\mathrm{End}_S(V)[t] \subset \mathrm{End}_{S[t]}(V[t])$. We consider the $(\bm{\mathfrak{l}}_S,K\cap L)$-module structure on $V[t]$ whose associated representation \(\bm{\pi}\) of \(\Upsilon(\bm{\mathfrak{l}}_S,K\cap L)\) is defined by $\bm{\pi}_{|K\cap L} = \pi_{|K\cap L}$ and
\begin{equation}\label{formula_second_ext}
    \bm{\pi}(X+tY + tH) = \pi(X) + t\pi(Y) +\pi(H) -t \rho_P(H),
\end{equation}
for every $X \in \mathfrak{l}\cap \mathfrak{k}$, $Y \in \mathfrak{m}_L \cap \mathfrak{s}$, $H \in \mathfrak{a}_L$.
We then define the $(\bm{\mathfrak{p}}_S, K\cap P)$-module
\begin{equation}\label{second_ext}
    \mathbf{ext}_{L,S}^P(V) = \mathrm{res}_{\bm{\mathfrak{p}}_S, K\cap P}^{\bm{\mathfrak{l}}_S,K\cap L}\left(V[t]\right).
\end{equation}
This defines a functor $\mathbf{ext}_{L,S}^P$ from $(\mathfrak{l}_S,K\cap L)$-modules to $(\bm{\mathfrak{p}}_S,K\cap P)$-modules.

\begin{definition}
    Consider the functor from the category of $(\mathfrak{l}_S,K\cap L)$-modules to the category of $(\bm{\mathfrak{g}}_S,K)$-modules defined by
    \[\mathbf{Ind}_{P,S}^G = \mathrm{ind}_{\bm{\mathfrak{p}}_S,K\cap P}^{\bm{\mathfrak{g}}_S,K} \circ \mathbf{ext}_{L,S}^P.\]
    If $V$ is any $(\mathfrak{l}_S,K\cap L)$-module, the $(\bm{\mathfrak{g}}_S,K)$-module $\mathbf{Ind}_{P,S}^G(V)$ is said to be \textbf{parabolically induced} from $V$ along the Mackey deformation.
\end{definition}

In what follows, when $S = \mathbb{C}$, we will simply write $\mathbf{ext}_L^P$ and $\mathbf{Ind}_P^G$ instead of $\mathbf{ext}_{L,\mathbb{C}}^P$ and $\mathbf{Ind}_{P,\mathbb{C}}^G$, respectively.

We now list several useful properties of these functors. The following proposition follows from the definitions (\ref{first_ext}) and (\ref{second_ext}) and Proposition \ref{specialization_P_at_0}.

\begin{proposition}\label{specializations_ext} Let $V$ be an $(\mathfrak{l}, K\cap L)$-module on which $\mathfrak{a}_{L}$ acts trivially and let $W$ be an $\mathfrak{a}_{L,S}$-module. For every $z \in \mathbb{C}^\times$, let $\delta_z$ denote the automorphism of $\mathfrak{a}_{L,S}$ given by multiplication by~$z^{-1}$
    \begin{enumerate}
        \item For every character $\lambda : S \to \mathbb{C}$, we have a natural identification
        \begin{gather*}
            \mathbf{ext}_{L,S}^P(V \otimes W)_\lambda = \mathbf{ext}_{L}^P(V \otimes W_\lambda),
        \end{gather*}
        given by the specialization map $\mathrm{id}_V \otimes_\mathbb{C} \mathrm{sp}_\lambda : V \otimes_\mathbb{C} W \to V \otimes_\mathbb{C} W_\lambda$, extended by $\mathbb{C}[t]$-linearity.
        \item The $(\bm{\mathfrak{l}}_S,K\cap L)$-module $\mathbf{ext}_{L,S}^P(V\otimes W)$ specializes at each $z \in \mathbb{C}$ as follows:
        \begin{gather*}
            \mathbf{ext}_{L,S}^P(V\otimes W)_z = \begin{cases}
                \mathrm{ext}_L^P\big(V \otimes \mathrm{res}_{\delta_z}(W)\big) &(z\neq 0), \\
                \mathrm{res}_{K\cap L}^{\mathfrak{l}, K\cap L}(V) \otimes W & (z = 0).
            \end{cases}
        \end{gather*}
        We recall that this equality is about $(\mathfrak{p}_S, K\cap P)$-modules in the first case and about $(\mathfrak{p}_{0,S}, K\cap P)$-modules in the second one. The specialization map $(V \otimes_\mathbb{C} W)[t] \to V \otimes_\mathbb{C} W$ is simply given by evaluation at $z$.
        \item We have a natural identification
        \begin{gather*}
            \mathrm{res}_{(K\cap P)_{S[t]}}^{\bm{\mathfrak{p}}_S, K\cap P}\left(\mathbf{ext}_{L,S}^P(V \otimes W) \right) = \mathrm{res}_{K\cap L}^{\mathfrak{l}, K\cap L}(V) \otimes_{\mathbb{C}} W[t],
        \end{gather*}
        where in the right part $W[t]$ is only considered as an $S[t]$-module. For the notation $(K\cap P)_{S[t]}$, see Convention \ref{notation_specialization_HC}.
    \end{enumerate}
\end{proposition}

From the above, we deduce the following.

\begin{proposition}\label{properties_induction}
    Let $V$ and $W$ be as above.
    \begin{enumerate}
        \item For every character $\lambda : S \to \mathbb{C}$, we have
        \begin{gather*}
            \mathbf{Ind}_{P,S}^G(V \otimes W)_\lambda = \mathbf{Ind}_{P}^G(V \otimes W_\lambda).
        \end{gather*}
        \item The $(\bm{\mathfrak{g}}_S, K)$-module $\mathbf{Ind}_{P,S}^G(V \otimes W)$ specializes at each $z \in \mathbb{C}^\times$ to the $(\mathfrak{g}_S,K)$-module
        \[\mathrm{Ind}_{P,S}^G(V \otimes \mathrm{res}_{\delta_z}W),\]
        and it specializes at $0$ to the $(\mathfrak{g}_{0,S},K)$ module
        \[\mathrm{Ind}^{G_0}_{P_0,S}(\mathrm{res}_{K\cap L}^{\mathfrak{l}, K\cap L}(V)\otimes W).\]
        \item We have
        \begin{gather*}
            \mathrm{res}_{K_{S[t]}}^{\bm{\mathfrak{g}}_S,K} \mathbf{Ind}_{P,S}^G(V \otimes W) = \mathrm{ind}_{K\cap L}^{K}(\mathrm{res}_{K\cap L}^{\mathfrak{l}, K\cap L} V) \otimes_\mathbb{C} W[t],
        \end{gather*}
        where $W[t]$ is only considered as an $S[t]$-module.
    \end{enumerate}
\end{proposition}

\begin{proof}
    The items \textit{(1)} and \textit{(2)} follow directly from the items \textit{(1)} and \textit{(2)} of Proposition \ref{specializations_ext} and the fact that induction commutes with specialization (Proposition~\ref{properties_ind_res}). The item \textit{(3)} is a consequence of Proposition~\ref{specializations_ext}\textit{(3)} and Lemma \ref{lemma_compact_pic}. The reason why we can apply the latter is the following: if $\mathfrak{n}$ is defined as just before Proposition \ref{specialization_P_at_0}, we have $\mathfrak{g}_S = \mathfrak{k}_S \oplus \mathfrak{a}_S \oplus \mathfrak{n}_S$ and $\mathfrak{p}_S = \mathfrak{l}_S\oplus\mathfrak{a}_S \oplus \mathfrak{n}_S$, as well as $\bm{\mathfrak{g}}_S = \mathfrak{k}_{S[t]} \oplus t\mathfrak{a}_S[t] \oplus t\mathfrak{n}_S[t]$ and $\bm{\mathfrak{p}}_S = \mathfrak{l}_{S[t]}\oplus t\mathfrak{a}_S[t] \oplus t\mathfrak{n}_S[t]$.
\end{proof}

To conclude this subsection, let us make a couple remarks.

\begin{remark}\label{commute_parInd_for}
    Using (\ref{commute_for_ind}), it is easy to see that the functor $\mathrm{for}_S^{\bm{\mathfrak{g}},K}\circ\mathbf{Ind}_{P,S}^G$ is naturally isomorphic to $\mathbf{Ind}_{P}^G \circ \mathrm{for}_S^{\mathfrak{l},K\cap L}$ and similarly for $\mathrm{for}_S^{\bm{\mathfrak{g}},K}\circ\mathrm{Ind}_{P,S}^G$ and $\mathrm{Ind}_{P}^G \circ \mathrm{for}_S^{\mathfrak{l},K\cap L}$.
\end{remark}

\begin{remark}\label{specialization_lambda}
    Applying the statement of Proposition \ref{properties_induction}\textit{(2)} concerning the specializations at $t = 1$ and at $t=0$ to Proposition \ref{properties_induction}\textit{(1)}, we obtain
        \begin{gather*}\mathrm{Ind}_{P,S}^G(V \otimes W)_\lambda = \mathrm{Ind}_{P}^G(V \otimes W_\lambda),\\
        \mathrm{Ind}_{P_0,S}^{G_0}(\mathrm{res}_{K\cap L}^{\mathfrak{l}, K\cap L}(V)\otimes W)_\lambda = \mathrm{Ind}_{P_0}^{G_0}(\mathrm{res}_{K\cap L}^{\mathfrak{l}, K\cap L}(V)\otimes W_\lambda)\end{gather*}
        for every $(\mathfrak{l},K\cap L)$-module $V$ annihilating $\mathfrak{a}_L$ and every $\mathfrak{a}_{L,S}$-module $W$.
\end{remark}

\subsection{A deformation of the center of the enveloping algebra}\label{Center}
The goal of this subsection is to study a central subalgebra $\mathcal{Z}(\bm{\mathfrak{g}})$ of $U(\bm{\mathfrak{g}})$, which somehow already appears in \cite{BHY_operator}*{Section 5.2} in the context of the C*-algebraic Mackey deformation. We describe the action of $\mathcal{Z}(\bm{\mathfrak{g}})$ on the parabolic induced modules and we also prove that the $K$-isotypic components of finitely generated $(\bm{\mathfrak{g}},K)$-modules are finitely generated over $\mathcal{Z}(\bm{\mathfrak{g}})$. The results thus obtained will be of great use in the last section.

What follows is essentially an adaptation to $(\bm{\mathfrak{g}},K)$-modules of the fundamentals of the theory of $(\mathfrak{g},K)$-modules, due to Harish-Chandra, for which we mainly refer to \cite{Wallach}*{Chapter 3}.

Let us begin with some notation on filtrations and gradings, which are always assumed to be defined over $\mathbb{N}$.
If $R$ is a filtered ring, we write $(\mathrm{f}_nR)_{n\in\mathbb{N}}$ for the corresponding filtration, and $\mathrm{f}_{-1}R = 0$. The graduation of any graded ring $R$ is denoted by $(\mathrm{gr}_nR)_{n\in\mathbb{N}}$. Every graded ring $R$ can (and will be) considered as a filtered ring, so that $\mathrm{f}_n R = \bigoplus_{k = 0}^n \mathrm{gr}_kR$. If $R$ is a filtered ring, its associated graded is denoted by $\mathrm{gr}R$ and we write $\mathrm{gr}_nR = \mathrm{f}_nR/\mathrm{f}_{n-1}R$ for each $n \in \mathbb{N}$. This notation brings no ambiguity in the case where $R$ is already graded, since it then canonicaly identifies with $\mathrm{gr}R$. If $\varphi : S \to R$ is a morphism of rings, with $R$ filtered, then $\mathrm{f}_nS = \varphi^{-1}(\mathrm{f}_nR)$ defines a filtration on $S$ and the associated graded morphism $\mathrm{gr}S \to \mathrm{gr}R$ is surjective (respectively injective) if and only if $\varphi$ is.

For a given commutative unital ring $A$, by a graded or filtered $A$-algebra we mean an associative $A$-algebra endowed with a grading or a filtration which is preserved by $A$. For any filtered $\mathbb{C}$-algebra $B$, we denote by $\mathcal{R}(B)$ the \textbf{Rees algebra} associated to $B$, defined as the filtered $\mathbb{C}[t]$-subalgebra of $B[t]$ spanned by $t^n\mathrm{f}_n B[t]$, where $n$ ranges over $\mathbb{N}$. Then $\mathcal{R}$ defines a functor from the category of filtered $\mathbb{C}$-algebras to the category of filtered $\mathbb{C}[t]$-algebras. Note that for any filtered $\mathbb{C}$-algebra $B$, we have \(\mathrm{gr}\mathcal{R}(B) = \textstyle\bigoplus_n t^n \mathrm{gr}_nB[t] = \mathcal{R}(\mathrm{gr}B).\) Here is another remark: if a $H$ is a group of filtered automorphisms of $B$, then we have $\mathcal{R}(B^H) = \mathcal{R}(B)^H$.

In the following, the universal enveloping algebras of Lie algebras are equipped with their canonical filtration \cite{GrAlLie1}*{§2, N°6} and the symmetric algebras of modules are endowed with their canonical graduation.

If $A$ is any commutative unital ring and $M$ is any $A$-module, we denote by $S(M)$ the symmetric $A$-algebra of $M$. We stress that if $A$ is a principal domain such as $\mathbb{C}[t]$, then every inclusion of $A$-modules $N \subset M$, with $M$ free, induces an embedding of $A$-algebras $S(N) \subset S(M)$.

\begin{proposition}
    The $\mathbb{C}[t]$-algebra homomorphism $U(\bm{\mathfrak{g}}) \to U(\mathfrak{g}[t])$ induced by the inclusion $\bm{\mathfrak{g}} \subset \mathfrak{g}[t]$ is injective. Moreover, the filtration of $U(\bm{\mathfrak{g}})$ induced by the canonical filtration of $U(\mathfrak{g}[t])$ via this inclusion is the canonical one.
\end{proposition}

\begin{proof}
    Since $\bm{\mathfrak{g}}$ and $\mathfrak{g}[t]$ are both free $\mathbb{C}[t]$-modules,  by the Poincaré-Birkhoff-Witt theorem \cite{GrAlLie1}*{Théorème 1 p.30}, the graded algebras associated to $U(\bm{\mathfrak{g}})$ and $U(\mathfrak{g}[t])$, with respect to their canonical filtrations, are $S(\bm{\mathfrak{g}})$ and $S(\mathfrak{g}[t]) = S(\mathfrak{g})[t]$, respectively. Moreover, the graded algebra homomorphism $S(\bm{\mathfrak{g}}) \to S(\mathfrak{g}[t])$ induced by the filtered homomorphism $U(\bm{\mathfrak{g}}) \to U(\mathfrak{g}[t])$ is exactly the one induced by the inclusion of $\mathbb{C}[t]$-modules $\bm{\mathfrak{g}} \subset \mathfrak{g}[t]$. The latter is injective and the two statements of the Lemma follow.
\end{proof}

Let $\mathcal{Z}(\mathfrak{g})$ be the center of the enveloping algebra of $\mathfrak{g}$, equipped with the filtration coming from $U(\mathfrak{g})$. Let us identify $\mathrm{gr} U(\mathfrak{g})$ with the symmetric algebra $S(\mathfrak{g})$ of $\mathfrak{g}$ in the usual way. The graded algebra of $\mathcal{Z}(\mathfrak{g})$ then identifies with the subalgebra $S(\mathfrak{g})^\mathfrak{g}$ of $\mathfrak{g}$-invariant elements of $S(\mathfrak{g})$, see \cite{Wallach}*{}. We define a subalgebra of $\mathcal{Z}(\mathfrak{g})[t] \subset U(\mathfrak{g}[t])$ as follows:
\[\mathcal{Z}(\bm{\mathfrak{g}}) = \mathcal{R} (\mathcal{Z}(\mathfrak{g})). \]
Since $\bm{\mathfrak{g}} \subset t\mathfrak{g}[t]$, we clearly have $\mathcal{Z}(\bm{\mathfrak{g}}) \subset U(\bm{\mathfrak{g}})$. Moreover, the graded algebra of $\mathcal{Z}(\bm{\mathfrak{g}})$ is $\mathcal{R}(S(\mathfrak{g})^\mathfrak{g}) = \mathcal{R}(S(\mathfrak{g}))^\mathfrak{g}$.

Our goal now is to find an analogue of the Harish-Chandra isomorphism for~$\mathcal{Z}(\bm{\mathfrak{g}})$. We begin with a brief account on the Harish-Chandra isomorphism and its relation to the Chevalley restriction theorem, see \cite{Wallach}*{§3.1-§3.2} for more information.
From now on, let us fix any Cartan subalgebra $\mathfrak{h}$ of $\mathfrak{g}$ and let us denote by $W(\mathfrak{g},\mathfrak{h})$ the associated Weyl group. We choose an associated positive root system and denote by $\mathfrak{n}^+$ and $\mathfrak{n}^-$ the sum of the positive and negative root spaces, respectively. Because of the direct sum decomposition $\mathfrak{g} = \mathfrak{n}^- \oplus \mathfrak{h} \oplus \mathfrak{n}^+$, the multiplication in $U(\mathfrak{g})$ induces a linear isomorphism $U(\mathfrak{n}^-) \otimes_\mathbb{C} U(\mathfrak{h}) \otimes_\mathbb{C} U(\mathfrak{n}^+) \to U(\mathfrak{g})$. Let $q_{\mathfrak{g}, \mathfrak{h}} : U(\mathfrak{g}) \to U(\mathfrak{h})$ be the linear map making the diagram
\[
\begin{tikzcd}
     U(\mathfrak{n}^-) \otimes_\mathbb{C} U(\mathfrak{h}) \otimes_\mathbb{C} U(\mathfrak{n}^+) \ar[r, "\sim"] \ar[dr, "\varepsilon \otimes \mathrm{id} \otimes \varepsilon"']& U(\mathfrak{g}) \ar[d, "q_{\mathfrak{g}, \mathfrak{h}}"] \\
    & U(\mathfrak{h}) 
\end{tikzcd}
\]
commute, where $\varepsilon$ denote both the counit of $U(\mathfrak{n}^-)$ and the counit of $U(\mathfrak{n}^+)$. For any $\lambda \in \mathfrak{h}^*$, we denote by $T_\lambda$ the automorphism of $U(\mathfrak{h})$ defined by the formula $T_\lambda(H) = H - \lambda(H)$, for every $H \in \mathfrak{h}$. We denote by $\rho_{\mathfrak{g},\mathfrak{h}}$ the half sum of the positive $(\mathfrak{g},\mathfrak{h})$-roots. The Harish-Chandra isomorphism $T_{\rho_{\mathfrak{g},\mathfrak{h}}} \circ q_{\mathfrak{g},\mathfrak{h}} : \mathcal{Z}(\mathfrak{g}) \to U(\mathfrak{h})^{W(\mathfrak{g},\mathfrak{h})}$ is denoted by $\gamma_{\mathfrak{g},\mathfrak{h}}$. We recall that $\gamma_{\mathfrak{g},\mathfrak{h}}$ is a morphism of filtered algebras, whose associated graded homomorphism $\mathrm{gr}(\gamma_{\mathfrak{g},\mathfrak{h}}) : S(\mathfrak{g})^\mathfrak{g} \to S(\mathfrak{h})^{W(\mathfrak{g},\mathfrak{h})}$ is the restriction of the algebra homomorphism $p_{\mathfrak{g}|\mathfrak{h}} : S(\mathfrak{g}) \to S(\mathfrak{h})$ induced by the orthogonal projection $\mathfrak{g} \twoheadrightarrow \mathfrak{h}$. If we identify $\mathfrak{g}$ and $\mathfrak{h}$ with their linear duals using the invariant bilinear form $(\cdot |\cdot )$ fixed in Subsection \ref{DeformationConvention}, the graded morphism $\mathrm{gr}(\gamma_{\mathfrak{g},\mathfrak{h}})$ coincides with the Chevalley restriction isomorphism. It follows that $\gamma_{\mathfrak{g},\mathfrak{h}}$ is an isomorphism of filtered algebras.

Let us extend all the maps introduced above by $\mathbb{C}[t]$-linearity. The Harish-Chandra isomorphism in particular extends to filtered isomorphisms
\begin{gather}\label{extended_HC_iso}\begin{gathered}
    \mathcal{Z}(\mathfrak{g})[t] \to U(\mathfrak{h}[t])^{W(\mathfrak{g},\mathfrak{h})},\\
    \mathcal{Z}(\bm{\mathfrak{g}}) \to \mathcal{R}\left(U(\mathfrak{h})^{W(\mathfrak{g},\mathfrak{h})}\right) = \mathcal{R}(U(\mathfrak{h}))^{W(\mathfrak{g},\mathfrak{h})}.\end{gathered}
\end{gather}
As a $\mathbb{C}$-algebra, $\mathcal{R}(U(\mathfrak{h})) = U(t\mathfrak{h}[t])$ is generated by $\{t\}\cup \,t\mathfrak{h}$, hence is finitely generated. From Hilbert's theorem on the finite generation of rings of invariants \cite{Mukai}*{Theorem 4.51}, it follows that $\mathcal{Z}(\bm{\mathfrak{g}})$ is a finitely generated $\mathbb{C}$-algbera.

Let us now express the action of $\mathcal{Z}(\bm{\mathfrak{g}})$ on parabolically induced modules. Let us assume that $\mathfrak{h}$ is the complexification of a Cartan subalgebra $\mathfrak{h}_\mathbb{R}$ containing $\mathfrak{a}_\mathbb{R}$. Let $P$ be a parabolic subgroup of $G$ with Levi factor denoted by $L$ and let $A_L$ be the split component of $L$. Let $\mathfrak{u}$ denote the nilpotent radical of $\mathfrak{p}$ and let $\bar{\mathfrak{u}} = \theta(\mathfrak{u})$. We assume that the chosen positive root system of $(\mathfrak{g},\mathfrak{h})$ is \textbf{compatible} with $P$ in the sense that $\mathfrak{u} \subset \mathfrak{n}^+$ and $\bar{\mathfrak{u}} \subset \mathfrak{n}^-$. Let $M_L$ be the essential component of $L$ and let $\mathfrak{h}_L = \mathfrak{h} \cap \mathfrak{m}_L$. Note that $\mathfrak{h} = \mathfrak{h}_L \oplus \mathfrak{a}_L$. Moreover, $\mathfrak{h}$ is a Cartan subalgebra of $\mathfrak{l}$ and $\mathfrak{h}_L$ is a Cartan subalgebra of $\mathfrak{m}_L$.

Let $q_{\mathfrak{g},\mathfrak{l}} : U(\mathfrak{g}) \to U(\mathfrak{l})$ be the linear map identifying with \[\varepsilon \otimes \mathrm{id} \otimes \varepsilon : U(\bar{\mathfrak{u}}) \otimes_\mathbb{C} U(\mathfrak{l})\otimes_\mathbb{C} U(\mathfrak{u})\to U(\mathfrak{l})\] under the isomorphism $U(\mathfrak{g}) \approx U(\bar{\mathfrak{u}}) \otimes_\mathbb{C} U(\mathfrak{l})\otimes_\mathbb{C} U(\mathfrak{u})$ given by multiplication. Again, $\varepsilon$ denotes the counits of $U(\bar{\mathfrak{u}})$ and $U(\mathfrak{u})$.

\begin{lemma}
    For every $n \in \mathbb{N}$ and $z \in \mathrm{f}_n\mathcal{Z}(\mathfrak{g})$, we have \[q_{\mathfrak{g},\mathfrak{l}}(z) \in \mathrm{f}_n\mathcal{Z}(\mathfrak{l}), \qquad z -q_{\mathfrak{g},\mathfrak{l}}(z) \in (\mathrm{f}_{n-1}U(\mathfrak{g}))\mathfrak{u}.\]
\end{lemma}

\begin{proof}
    Note first that $q_{\mathfrak{g},\mathfrak{l}}$ is a filtered linear homomorphism. Let $(X_i)_{i= 1...m}$ be a basis of $\mathfrak{l}$ and let $(Y_j)_{j= 1 ...n}$ be a basis of $\mathfrak{u}$. Let us write $\bar{Y}_j = \theta(Y_j)$, so that $(Y_j)_{j= 1 ...n}$ is a basis of $\bar{\mathfrak{u}}$. For any $\mu \in \mathbb{N}^m$, $\nu \in\mathbb{N}^n$ and $\bar{\nu} \in\mathbb{N}^n$, let us write
    \[\bar{Y}^{\bar{\nu}} X^\mu Y^\nu = \bar{Y}_1^{\bar{\nu}_1} \cdots\bar{Y}_n^{\bar{\nu}_n}X^{\mu_1}_1 \cdots X_m^{\mu_m} Y_1^{\nu_1}\cdots Y_n^{\nu_n} \in U(\mathfrak{g}). \]
    Let $d \in \mathbb{N}$ and $z \in \mathrm{f}_d\mathcal{Z}(\mathfrak{g})$. By writing $z$ as a linear combination of the above basis and expressing the fact that $z$ commutes with $\mathfrak{a}_L$, we easily see that $z -q_{\mathfrak{g},\mathfrak{l}}(z)$ is both an element of $(\mathrm{f}_{d-1}U(\mathfrak{g}))\mathfrak{u}$ and $\bar{\mathfrak{u}}(\mathrm{f}_{d-1}U(\mathfrak{g}))$. It follows that for every $x \in U(\mathfrak{l})$, we have $xq_{\mathfrak{g},\mathfrak{l}}(z) = q_{\mathfrak{g},\mathfrak{l}}(xz) = q_{\mathfrak{g},\mathfrak{l}}(zx) = q_{\mathfrak{g},\mathfrak{l}}(z)x$, which concludes the proof.
\end{proof}

Let us extend the map $q_{\mathfrak{g},\mathfrak{l}}$ by $\mathbb{C}[t]$-linearity. From the previous lemma, we deduce the following result.

\begin{corollary}\label{lemma_Z_n}
    For every $z \in \mathcal{Z}(\bm{\mathfrak{g}})$, we have \[q_{\mathfrak{g},\mathfrak{l}}(z) \in\mathcal{Z}(\bm{\mathfrak{l}}), \qquad z -q_{\mathfrak{g},\mathfrak{l}}(z) \in U(\bm{\mathfrak{g}})t\mathfrak{u}[t],\]
    where $\mathcal{Z}(\bm{\mathfrak{l}}) = \mathcal{R}(\mathcal{Z}(\mathfrak{l})) \subset U(\bm{\mathfrak{l}})$. 
\end{corollary}

Since $W(\mathfrak{l},\mathfrak{h}) = W(\mathfrak{m}_L,\mathfrak{h}_L)$ acts trivially on $U(\mathfrak{a}_L)$, the multiplication induces an isomorphism of algebras $U(\mathfrak{h}_L)^{W(\mathfrak{m}_L,\mathfrak{h}_L)} \otimes_\mathbb{C} U(\mathfrak{a}_L) \approx U(\mathfrak{h})^{W(\mathfrak{l},\mathfrak{h})}$. Let us also point out that $U(\mathfrak{h})^{W(\mathfrak{g},\mathfrak{h})} \subset U(\mathfrak{h})^{W(\mathfrak{l},\mathfrak{h})}$.

Since $\mathfrak{l}$ is the commuting direct sum $\mathfrak{m}_L \oplus \mathfrak{a}_L$, the multiplication in $U(\mathfrak{l})$ induces a linear isomorphism $U(\mathfrak{m}_L)\otimes_\mathbb{C} U(\mathfrak{a}_L) \approx U(\mathfrak{l})$. Under this isomorphism, the center $\mathcal{Z}(\mathfrak{l})$ of $U(\mathfrak{l})$ identifies with $\mathcal{Z}(\mathfrak{m}_L) \otimes_\mathbb{C} U(\mathfrak{a}_L)$, while the Harish-Chandra isomorphism \(\gamma_{\mathfrak{l},\mathfrak{h}}\) identifies with $\gamma_{\mathfrak{m}_L,\mathfrak{h}_L} \otimes \mathrm{id}_{U(\mathfrak{a}_L)}$.
Similarly, the multiplication in $U(\bm{\mathfrak{l}})$ induces an isomorphism of $\mathbb{C}[t]$-algebras $U(\bm{\mathfrak{m}}_L) \otimes_{\mathbb{C}[t]}U(t\mathfrak{a}_L[t]) \approx U(\bm{\mathfrak{l}})$ under which $\mathcal{Z}(\bm{\mathfrak{l}})$ identifies with $\mathcal{Z}(\bm{\mathfrak{m}}_L) \otimes_{\mathbb{C}[t]}U(t\mathfrak{a}_L[t])$. Then $\gamma_{\mathfrak{l},\mathfrak{h}}$ maps $\mathcal{Z}(\bm{\mathfrak{l}})$ onto 
\begin{equation} \label{isomorphism_UhW}
    \mathcal{R}\left(U(\mathfrak{h})^{W(\mathfrak{l},\mathfrak{h}_L)}\right) = U(t\mathfrak{h}[t])^{W(\mathfrak{l},\mathfrak{h})} \approx U(t\mathfrak{h}_L[t])^{W(\mathfrak{m}_L,\mathfrak{h}_L)} \otimes_{\mathbb{C}[t]} U(t\mathfrak{a}_L[t]),
\end{equation}
see Formula (\ref{extended_HC_iso}).

In the following, we denote by $\delta_t : U(t\mathfrak{a}_L[t]) \to U(\mathfrak{a}_L[t])$ the isomorphism of $\mathbb{C}[t]$-algebras defined by
\begin{equation*}
    \delta_t(tH) = H\qquad (H \in \mathfrak{a}_L).
\end{equation*}

Let $V$ be an $(\mathfrak{l},K\cap L)$-module annihilating $\mathfrak{a}_L$ and let $W$ be a $\mathfrak{a}_L$-module. Let us denote by $\pi_V$ and $\pi_W$ the associated representations of $\Upsilon(\mathfrak{m}_L,K\cap L)$ and $U(\mathfrak{a}_L)$ on $V$ and $W$, respectively. We denote by $\mathrm{Ch}_V : U(\mathfrak{h}_L)^{W(\mathfrak{m}_L,\mathfrak{h}_L)} \to \mathrm{End}_\mathbb{C}(V)$ the morphism of algebras defined by
\begin{equation*}
    \mathrm{Ch}_V \circ \gamma_{\mathfrak{m}_L,\mathfrak{h}_L} = \pi_{V|\mathcal{Z}(\mathfrak{m}_L)}.
\end{equation*}
Let us extend $\pi_V$, $\pi_W$ and $\mathrm{Ch}_V$ in a $\mathbb{C}[t]$-linear way, so that they form representations of $\Upsilon(\mathfrak{m}_L[t],K\cap L)$, $U(\mathfrak{a}_L[t])$ and $U(\mathfrak{h}_L[t])^{W(\mathfrak{m}_L,\mathfrak{h}_L)}$ on $V[t]$, $W[t]$ and $V[t]$, respectively.

\begin{proposition}\label{action_center_induced_modules}
    Let $\mathbf{Ch}$ be the representation of $\mathcal{R}[U(\mathfrak{h})^{W(\mathfrak{l},\mathfrak{h})}]$ on $\mathbf{ext}_L^P(V \otimes W)$, which we identify with $V[t] \otimes_{\mathbb{C}[t]}W[t]$, defined by
    \[\mathbf{Ch}(z \otimes h) = \mathrm{Ch}_V(z) \otimes\pi_W(\delta_th)\]
    for every $z \in U(t\mathfrak{h}[t])^{W(\mathfrak{l},\mathfrak{h})}$ and $ h\in U(t\mathfrak{a}_L[t])$,
    see Formula (\ref{isomorphism_UhW}).
    
    The action of $\mathcal{Z}(\bm{\mathfrak{g}})$ on $\mathbf{Ind}_P^G(V\otimes W)$ coincides with the one induced from the representation $\mathbf{Ch} \circ \gamma_{\mathfrak{g}, \mathfrak{h}}$ on $\mathbf{ext}_L^P(V \otimes W)$ by the functor $\mathrm{ind}_{\bm{\mathfrak{p}},K\cap P}^{\bm{\mathfrak{g}},K}$.
\end{proposition}

\begin{proof}
    Let $\bm{\pi}$ be the natural representation of $\Upsilon(\bm{\mathfrak{p}},K\cap L)$ on $\mathbf{ext}_L^P(V \otimes W)$.
    
    Let $x \in U(\bm{\mathfrak{g}},K)$ and let $m \in \mathbf{ext}_L^P(V \otimes W)$. We consider $x \otimes m$ as an element of $\mathbf{Ind}_P^G(V\otimes W)$. Using Corollary \ref{lemma_Z_n} and the fact that $\mathcal{Z}(\bm{\mathfrak{g}})$ is a central subalgebra of $\Upsilon(\bm{\mathfrak{g}},K)$, we obtain for every $z \in \mathcal{Z}(\bm{\mathfrak{g}})$:
    \[z(x \otimes m) = zx \otimes m = xz \otimes m = x(z - q_{\mathfrak{g},\mathfrak{l}}(z)) \otimes m + xq_{\mathfrak{g},\mathfrak{l}}(z) \otimes m = x \otimes \bm{\pi}(q_{\mathfrak{g},\mathfrak{l}}(z))m. \]
    Thus, the proposition will be proved if we check that $(\bm{\pi} \circ q_{\mathfrak{g},\mathfrak{l}})_{|\mathcal{Z}(\bm{\mathfrak{g}})} = \mathbf{Ch} \circ \gamma_{\mathfrak{g}, \mathfrak{h}}$.
    
    Let $\rho_P\in \mathfrak{h}^*$ denote the half trace of the adjoint action on $\mathfrak{p}$. From Formula (\ref{formula_second_ext}), it is clear that for every $z \in \mathcal{Z}(\bm{\mathfrak{m}}_L)$ and every $h \in U(t\mathfrak{a}_L[t])$, we have
    \[\bm{\pi}(zh) = \pi_V(z) \otimes (\pi_W \circ\delta_t)(T_{\rho_P}(h)) = (\mathrm{Ch}_V \circ \gamma_{\mathfrak{m}_L,\mathfrak{h}_L})(z) \otimes (\pi_W \circ\delta_t)(T_{\rho_P}(h)).\]
    Now, since $\rho_P$ vanishes on $\mathfrak{h}_L$ and \(\gamma_{\mathfrak{l},\mathfrak{h}}\) coincides with $\gamma_{\mathfrak{m}_L,\mathfrak{h}_L} \otimes \mathrm{id}_{U(\mathfrak{a}_L)}$:
    \begin{gather*}
        \left(\gamma_{\mathfrak{m}_L,\mathfrak{h}_L} \otimes \delta_t \circ T_{\rho_P|U(t\mathfrak{a}_L[t])}\right) \circ q_{\mathfrak{g},\mathfrak{l}} =(\mathrm{id}_{U(\mathfrak{h}_L[t])}\otimes \delta_t)\circ T_{\rho_P} \circ \gamma_{\mathfrak{l},\mathfrak{h}}\circ q_{\mathfrak{g},\mathfrak{l}}.
    \end{gather*}
    Using $\rho_{\mathfrak{g},\mathfrak{h}} = \rho_P + \rho_{\mathfrak{l},\mathfrak{h}}$ and $q_{\mathfrak{g},\mathfrak{h}} = q_{\mathfrak{l},\mathfrak{h}} \circ q_{\mathfrak{g},\mathfrak{l}}$, we can make the following simplification:
    \begin{gather*}
        T_{\rho_P} \circ \gamma_{\mathfrak{l},\mathfrak{h}}\circ q_{\mathfrak{g},\mathfrak{l}} = T_{\rho_P} \circ T_{\rho_{\mathfrak{l},\mathfrak{h}}} \circ q_{\mathfrak{l},\mathfrak{h}} \circ q_{\mathfrak{g},\mathfrak{l}} = T_{\rho_{\mathfrak{g},\mathfrak{h}}} \circ q_{\mathfrak{g},\mathfrak{h}} = \gamma_{\mathfrak{g},\mathfrak{h}},
    \end{gather*}
    Thus, we obtain $(\bm{\pi} \circ q_{\mathfrak{g},\mathfrak{l}})_{|\mathcal{Z}(\bm{\mathfrak{g}})} = \mathbf{Ch} \circ \gamma_{\mathfrak{g}, \mathfrak{h}}$,
    which concludes the proof.
\end{proof}

We wish now to prove the following theorem.

\begin{theorem} \label{finitely_gen_center}
    Let $V$ be a finitely generated $(\bm{\mathfrak{g}},K)$-module. Then every $K$-isotypical component of $V$ is finitely generated over $\mathcal{Z}(\bm{\mathfrak{g}})$.
\end{theorem}

This is the analogue of Harish-Chandra's theorem \cite{Wallach}*{§3.4.1} for $(\bm{\mathfrak{g}},K)$-modules. The proof will essentially be the same.

\begin{proof}
    Let $W$ be a finite $\mathbb{C}$-dimensional $K$-submodule of $V$ such that $U(\bm{\mathfrak{g}})W = V$.
    
    We define an $\mathbb{N}$-filtration on $V$ as follows:
    \[\mathrm{f}_{n}V = [\mathrm{f}_nU(\bm{\mathfrak{g}})]W \qquad (n \in \mathbb{N}).\]
    Since $\mathrm{f}_m U(\bm{\mathfrak{g}}) \mathrm{f}_n V \subset \mathrm{f}_{n+m}V$ for every $n,m \in \mathbb{N}$, the graded $\mathbb{C}[t]$-module $\mathrm{gr}V$ inherits an $S(\bm{\mathfrak{g}})$-module structure. Note that $S(\mathfrak{k}[t])$ acts trivially on $\mathrm{gr} V$ for this structure, because $\mathrm{f}_n V$ is a $K_{\mathbb{C}[t]}$-submodule of $V$ for every $n \in \mathbb{N}$. This implies that \[\mathrm{gr}_n V = [\mathrm{gr}_nU(\mathfrak{g})] W =  [\mathrm{gr}_nS(t\mathfrak{s}[t]) ]W\qquad (n \in \mathbb{N}),\] and hence $\mathrm{gr} V = S(t\mathfrak{s}[t]) W$. In particular, $\mathrm{gr} V$ is finitely generated over $S(t\mathfrak{s}[t])$.

    Since $\mathrm{gr} V$ is a direct sum of $K_{\mathbb{C}[t]}$-subquotients of $V$, it has a natural $K_{\mathbb{C}[t]}$-module structure. Let $\mu$ be any $K$-type. The isotypical component of $\mathrm{gr}(V)$ of type $\mu$ is $\mathrm{gr}(V^\mu)$, the associated graded of the isotypical component of $V$ of type $\mu$. Arguing exactly as in \cite{Wallach}*{§3.4.3}, one can check that $\mathrm{gr}(V^\mu)$ is finitely generated over $S(t\mathfrak{s}[t])^K$.

    Let $p_{\mathfrak{g}|\mathfrak{s}}$ be the algebra homomorphism $S(\mathfrak{g}) \to S(\mathfrak{s})$ induced by the orthogonal projection $\mathfrak{g} \twoheadrightarrow \mathfrak{s}$, which we extend by $\mathbb{C}[t]$-linearity to a morphism $S(\bm{\mathfrak{g}}) \to S(t\mathfrak{s}[t])$. For every $n,m \in \mathbb{N}$ and every $\mathbb{C}$-subspace $N \subset \mathrm{f}_nV$, we have
    \begin{equation}\label{center_versus_grading_V}
        \mathrm{gr}_{n+m}\left([\mathrm{f}_m \mathcal{Z}(\bm{\mathfrak{g}})] N\right) = [\mathrm{gr}_m \mathcal{Z}(\bm{\mathfrak{g}})] \mathrm{gr}_n N = p_{\mathfrak{g}|\mathfrak{s}}(\mathrm{gr}_m S(t\mathfrak{g}[t])^\mathfrak{g}) \mathrm{gr}_n N.
    \end{equation}
    Let $J$ denote the subalgebra $p_{\mathfrak{g}|\mathfrak{s}}(S(\mathfrak{g})^\mathfrak{g})$ of $S(\mathfrak{s})$. We have $p_{\mathfrak{g}|\mathfrak{s}}(S(t\mathfrak{g}[t])^\mathfrak{g})  = \mathcal{R}(J)$, the Rees algebra of $J$. From (\ref{center_versus_grading_V}), it follows that for every $\mathbb{C}$-linear subspace $N$ of $V$, the associated graded of $\mathcal{Z}(\bm{\mathfrak{g}}) N$ contains $\mathcal{R}(J) \mathrm{gr}N$.

    Let $\delta_t : S(t\mathfrak{s}[t]) \to S(\mathfrak{s}[t])$ be the isomorphism of $\mathbb{C}[t]$-algebras induced by the isomorphism of $\mathbb{C}[t]$-modules $t\mathfrak{s}[t] \to \mathfrak{s}[t]$ consisting in dividing by $t$. Since $S(\mathfrak{s})^K$ is a finitely generated $J$-module, see \cite{Wallach}*{§3.4.4}, $S(\mathfrak{s})^K[t] = \delta_t(S(t\mathfrak{s}[t])^K)$ is a finitely generated $J[t] = \delta_t(\mathcal{R}(J))$-module, therefore $S(t\mathfrak{s}[t])^K$ is finitely generated as an $\mathcal{R}(J)$-module.

    Let us conclude the proof. Let $\mu$ be any $K$-type. Then $\mathrm{gr}(V^\mu)$ is finitely generated over $S(t\mathfrak{s}[t])^K$ and since the latter is finitely generated over $\mathcal{R}(J)$, we obtain that $\mathrm{gr}(V^\mu)$ is finitely generated over $\mathcal{R}(J)$. Hence, there exists a finite dimensional $\mathbb{C}$-linear subspace $N$ of $V^\mu$ such that $\mathrm{gr}(V^\mu) = \mathcal{R}(J) \mathrm{gr}N$. The latter is contained in $\mathrm{gr}(\mathcal{Z}(\bm{\mathfrak{g}})N)$, and given that $\mathcal{Z}(\bm{\mathfrak{g}})N \subset \mathcal{Z}(\bm{\mathfrak{g}})V^\mu \subset V^\mu$, we have $\mathrm{gr}(\mathcal{Z}(\bm{\mathfrak{g}})N) = \mathrm{gr}(V^\mu)$. We conclude that $V^\mu = \mathcal{Z}(\bm{\mathfrak{g}})N$, thus proving that $V^\mu$ is a finitely generated $\mathcal{Z}(\bm{\mathfrak{g}})$-module.
\end{proof}

\begin{corollary}\label{finite_generacy_chiUchi}
    For every finite subset $F$ of $\widehat{K}$, the algebra $\chi_F U(\bm{\mathfrak{g}},K) \chi_F$ is finitely generated as a $\mathcal{Z}(\bm{\mathfrak{g}})$-module.
\end{corollary}

\begin{proof}
    If $V = U(\bm{\mathfrak{g}},K)\chi_F$, then $V$ is generated as a $(\bm{\mathfrak{g}},K)$-module by $\chi_F$, hence $\chi_F V = \oplus_{\mu \in F}V^\mu$ is finitely generated over $\mathcal{Z}(\bm{\mathfrak{g}})$ by Theorem \ref{finitely_gen_center}.
\end{proof}

\begin{corollary}\label{admissibility_simple}
    Every simple $(\bm{\mathfrak{g}},K)$-module is admissible.
\end{corollary}

\begin{proof}
    By the Poincaré-Birkhoff-Witt theorem \cite{GrAlLie1}*{Théorème 1 p.30}, $U(\bm{\mathfrak{g}})$ is countably generated as a $\mathbb{C}[t]$-module and has thus countable $\mathbb{C}$-dimension. Proposition \ref{structureU} and the fact that $\widehat{K}$ is countable therefore shows that $U(\bm{\mathfrak{g}},K)$ has countable $\mathbb{C}$-dimension. Then, since $\mathcal{Z}(\bm{\mathfrak{g}})$ is central in $\Upsilon(\bm{\mathfrak{g}},K)$, Schur's lemma \cite{Bourbaki}*{§3, N°2} ensures that $\mathcal{Z}(\bm{\mathfrak{g}})$ acts on every simple $(\bm{\mathfrak{g}},K)$-module by scalars. In conclusion, if $M$ is any simple $(\bm{\mathfrak{g}},K)$-module, every $K$-isotypic component of $M$ is finitely generated over $\mathcal{Z}(\bm{\mathfrak{g}})$ by Theorem \ref{finitely_gen_center}, hence finite $\mathbb{C}$-dimensional.
\end{proof}

\section{Spectral structure of the algebraic Mackey deformation}\label{sec3}

As shown by \cites{Afgoustidis3,AfgoustidisClare}, the unitary Mackey analogy essentially relies on the piecewise triviality of the spectrum of the continuous field of reduced C*-algebra over $\mathbb{R}$ which underlies the contraction of a real reductive group onto its motion group. In this section, we prove that the same phenomena appears for the algebraic Mackey deformation. After an introductory subsection on the spectra of rings (Subsection \ref{spectral_properties_rings}) and a characterization of the Mackey bijection in terms of lowest $K$-types (Subsection \ref{the_mackey_bijection_revisited}), we show that that Delorme's generalized Harish-Chandra homomorphisms \cite{Delorme84} are in a sense constant along the Mackey deformation (Subsection \ref{HC_homomorphisms}). Eventually, we derive from this all the results announced in the introduction (Subsection \ref{main_theorems}).

\subsection{Spectral properties of ring homomorphisms}\label{spectral_properties_rings}

We begin this section by giving some vocabulary on the spectra of rings and algebras. The definitions introduced here will be used in the statement of our main theorems.

Let $R$ be any ring. The primitive spectrum of $R$, denoted by $\mathrm{Prim}(R)$, is always endowed with the Jacobson topology \cite{Jacobson}. By the \textbf{simple spectrum} of $R$, we mean the set $\mathrm{Sim}(R)$ of isomorphism classes of simple $R$-modules. There is a canonical (surjective) map $\mathrm{Sim}(R) \to \mathrm{Prim}(R)$, which associates to any simple module $M$ its annihilator in $R$. We recall that $\mathrm{Prim}(R)$ is a $\mathrm{T_1}$ space if and only if each of its primitive ideals are maximal among all primitive ideals, see \cite{Dixmier}*{§3.1}.

\begin{definition}\label{isospectral}
    A ring homomorphism $f : R \to B$ is said to be \textbf{isospectral} (\textbf{epispectral}, respectively) if
    \begin{enumerate}
        \item every simple module $M$ of $B$ restricts via $f$ to a simple $R$-module $f^*M$,
        \item the map $\mathrm{Sim}(B) \to \mathrm{Sim}(R)$ induced by $M \mapsto f^*M$ is bijective (injective, respectively).
    \end{enumerate}
\end{definition}

\begin{remark}\label{two_out_of_three}
    Let $f: R \to B$ and $\varphi : S \to R$ be two ring homomorphisms. In order for one of the three morphisms $f$, $\varphi$, $f\circ \varphi$ to be isospectral, it is enough that the two other are isospectral.
\end{remark}

In typical cases, these properties can be checked pointwise, as follows.

\begin{lemma}\label{isospectral_specializations}
    Let $A$ be a unital commutative $\mathbb{C}$-algebra and let $f : R \to B$ be a morphism of unital $A$-algebras. For every character $\varphi : A \to \mathbb{C}$, we consider the specializations $R_\varphi$ and $B_\varphi$ of $R$ and $B$ at $\varphi$, as well as the morphism of $\mathbb{C}$-algebras $f_\varphi : R_\varphi \to B_\varphi$ induced by $\varphi$.

    Assume that $R$ and $B$ have countable $\mathbb{C}$-dimension. Then $f$ is isospectral, or epispectral, if and only $f_\varphi$ is, for every character $\varphi : A \to \mathbb{C}$.
\end{lemma}

\begin{proof}
    For each $S \in \{R,B\}$ and every character $\varphi : A \to \mathbb{C}$, let  $\mathrm{Sim}_\varphi(S)$ denote the set of isomorphism classes of simple $S$-module on which $A$ acts through $\varphi$.
    Since $R$ and $B$ have countable $\mathbb{C}$-dimension, all of their simple modules have countable dimension and by Schur's lemma \cite{Bourbaki}*{§3, N°2}, $A$ must act on each of them through a character. In other words, for each $S \in \{R,B\}$, $\mathrm{Sim}(S)$ is the disjoint union of all $\mathrm{Sim}_\varphi(S)$. One can then easily check that $f$ is isospectral (respectively, epispectral) if and only if, for every character $\varphi$ of $A$, the morphism $f_\varphi$ is isospectral (respectively, epispectral).
\end{proof}

Following \cite{BaumNistor}, we also introduce the terminology below.

\begin{definition}[Baum--Nistor]\label{definition_spectrum_preserving}
    A ring homomorphism $f : R \to B$ is said to be \textbf{spectrum-preserving} if
    \begin{enumerate}
        \item for every $J_B \in \mathrm{Prim}(B)$, there exists a unique $J_R \in \mathrm{Prim}(R)$ such that $f^{-1}(J_B) \subset J_R $,
        \item the map $f^* : J_B \mapsto J_R$ thus defined is a bijection of $\mathrm{Prim}(B)$ onto $\mathrm{Prim}(R)$.
    \end{enumerate}
\end{definition}

In practice, we will use the isospectral property to detect the spectrum-preserving property, as follows. Let us call \textbf{point-primitive} any ring $R$ such that $\mathrm{Prim}(R/J)$ is a one-point set for all $J \in \mathrm{Prim}(R)$. Equivalently, a ring $R$ is point-primitive if $\mathrm{Prim}(R)$ is a $\mathrm{T_1}$ space and the canonical map $\mathrm{Sim}(R) \to \mathrm{Prim}(R)$ is injective. Note that being point-primitive is stable under taking quotients. The rings involved in the next subsections will always be point-primitive thanks to the admissibility criterion given by Proposition \ref{criterion_point_primitive}.

\begin{lemma}\label{isospectral_implies_spectrum_preserving}
    Let $f : R \to B$ be an isospectral ring homomorphism. If $R$ is point-primitive, then $f$ is spectrum preserving and $f^*(J) = f^{-1}(J)$ for every primitive ideal $J$ of $B$.
\end{lemma}

\begin{proof}
    Assume that $R$ is point primitive. Since $f$ is isospectral, every simple $B$-module restricts to a simple $R$-module, hence $f^{-1}(J)$ is primitive for every $J \in \mathrm{Prim}(B)$. Now every primitive ideal of $R$ is maximal among all primitive ideals, because $R$ is point-primitive. This proves the first condition of Definition~\ref{definition_spectrum_preserving}. Moreover, we have a commutative diagram
    \[\begin{tikzcd}
        \mathrm{Sim}(B) \ar[r, "f^*"] \ar[d, twoheadrightarrow] & \mathrm{Sim}(R) \ar[d, twoheadrightarrow]\\
        \mathrm{Prim}(B) \ar[r, "f^*"] & \mathrm{Prim}(R)
    \end{tikzcd},\]
    where $f^*(J) = f^{-1}(J)$ for every $J \in \mathrm{Prim}(B)$. Since $f$ is isospectral and $R$ is point-primitive, the upper horizontal arrow and the right vertical arrow are bijections. This implies that the left vertical arrow is injective, hence bijective. The map $f^* : \mathrm{Prim}(B) \to \mathrm{Prim}(R)$ is thus bijective, which proves the second condition of Definition~\ref{definition_spectrum_preserving}.
\end{proof}

Under finite type conditions, a spectrum-preserving algebra homomorphism induces a homeomorphism of primitive spectra and an isomorphism in periodic cyclic homology, see \cite{BaumNistor}. We state this precisely for future use.

\begin{theorem}[Baum--Nistor]\label{baum_nistor_consequences}
    Let $A$ be a unital finitely generated commutative $\mathbb{C}$-algebra. Let $f : R \to B$ be a morphism of $A$-algebras which are finitely generated as $A$-modules. If $f$ is spectrum preserving, then $f^* : \mathrm{Prim}(B) \to \mathrm{Prim}(R)$ is a homeomorphism and the map $\mathrm{HP}_*(R) \to \mathrm{HP}_*(B)$ induced by $f$ in periodic cyclic homology is an isomorphism.
\end{theorem}

For an introduction to periodic cyclic homology, see for example \cite{Loday}*{§5.1}.

\subsection{The Mackey bijection revisited}\label{the_mackey_bijection_revisited}
Let us now review the classification of the simple $(\mathfrak{g},K)$-modules and that of the simple $(\mathfrak{g}_0,K)$-modules. We will then give a new description of the admissible Mackey bijection \cite{Afgoustidis1}*{§5}.

We start by Vogan's classification of the simple $(\mathfrak{g},K)$-modules, for which we refer to \citelist{\cite{VoganAnnals} \cite{VoganGreen}*{Chapter 6, §5-6}}. Let $\mathrm{Sim}(\mathfrak{g},K)$ denote the set of isomorphism classes of simple $(\mathfrak{g},K)$-modules.

\begin{definition}
    A \textbf{cuspidal datum of $G$} is a pair $\sigma = (L_\sigma, D_\sigma)$ consisting of a Levi subgroup $L_\sigma$ of $G$ and a discrete series $(\mathfrak{m}_\sigma, K\cap M_\sigma)$-module $D_\sigma$, where $M_\sigma$ denotes the essential component of $L_\sigma$.

    Two cuspidal data $\sigma$ and $\sigma'$ are said to be \textbf{equivalent} if there exists $k \in K$ such that $L_{\sigma'} = kL_\sigma k^{-1}$ and $D_{\sigma}$ is isomorphic to $D_{\sigma'}$ via to conjugation by~$k$.
\end{definition}

In the context of this definition, the split component of $L_\sigma$ is denoted by $A_\sigma$ and the intersection of $L_\sigma$ with $K$ is denoted by $K_\sigma$. Moreover, the $(\mathfrak{m}_\sigma,K_\sigma)$-module $D_\sigma$ will be considered as a $(\mathfrak{l}_\sigma,K_\sigma)$-module annihilating $\mathfrak{a}_\sigma$.

Let us fix a cuspidal datum $\sigma$ of $G$. Given a parabolic subgroup $P$ of $G$ with Levi factor $L_\sigma$ and $\lambda \in \mathfrak{a}_\sigma^*$, we define
\begin{equation}\label{definition_I^P} I_{\sigma,\lambda}^P = \mathrm{Ind}_P^G(D_{\sigma} \otimes \mathbb{C}_\lambda).\end{equation}
The common set of lowest $K$-types of all $I_{\sigma,\lambda}^P$ for $\lambda \in \mathfrak{a}_\sigma^*$ is denoted by
\begin{equation*}
    \widehat{K}(\sigma).
\end{equation*}
We have the following.

\begin{proposition}\label{multiplicityone}
    If $\lambda\in \mathfrak{a}_\sigma^*$ and $P$ is a parabolic subgroup of $G$ with Levi factor $L_\sigma$, then every lowest $K$-type of $I_{\sigma,\lambda}^P$ occurs with multiplicity one.
\end{proposition}

For every $\lambda$ and $P$ as above and $\mu \in \widehat{K}(\sigma)$, let
\begin{equation*}
    J_{\sigma,\lambda}(\mu) \in \mathrm{Sim}(\mathfrak{g},K)
\end{equation*} denote the unique simple subquotient of $I_{\sigma,\lambda}^P$ containing $\mu$. This subquotient does not depend on the choice of $P$ by \cite{VoganGreen}*{Proposition 4.20}.

Let us introduce some notation. We put \begin{equation}\label{notation_SgKsigma}\mathrm{Sim}(\mathfrak{g},K)_\sigma = \{J_{\sigma,\lambda}(\mu) : \lambda \in \mathfrak{a}_\sigma^*, \mu \in \widehat{K}(\sigma) \}.\end{equation}
Let us also write
\[W_\sigma = \{k \in N_K(A_\sigma) : (k\cdot k^{-1})_\ast D_\sigma \approx D_\sigma\}/Z_K(A_\sigma).\]
From \cite{VoganGreen}*{Chapter 4, §3}, we have a distinguished split exact sequence of groups
\[1 \to W_\sigma' \to W_\sigma \to R_\sigma \to 1, \]
where $R_\sigma$, called the $R$-group of $\sigma$, is a product of copies of $\mathbb{Z}/2\mathbb{Z}$; we also have a simply transitive action of the dual group $\hat{R}_\sigma$ on $\widehat{K}(\sigma)$. For every $\lambda \in \mathfrak{a}_\sigma^*$, let $W_{\sigma, \lambda}$ be the stabilizer of $\lambda$ in $W_\sigma$ and let $\hat{R}_\sigma(\lambda)$ be the orthogonal of $W_{\sigma, \lambda}$ in $\hat{R}_\sigma$. Vogan's classification is given by the two following theorems.

\begin{theorem}\label{vogan_classification} Given two pairs $(\lambda, \mu)$ and $(\lambda', \mu')$ of elements of $\mathfrak{a}_\sigma^* \times \widehat{K}(\sigma)$, we have $J_{\sigma,\lambda}(\mu) = J_{\sigma,\lambda'}(\mu')$ if and only if $\lambda' \in W_\sigma\lambda$ and $\mu' \in \hat{R}_\sigma(\lambda) \mu$.
\end{theorem}

\begin{theorem}\label{partition}
    Let $\mathcal{C}$ be a complete set of representatives for the equivalence classes of cuspidal data of $G$. The collection of all $\mathrm{Sim}(\mathfrak{g},K)_\sigma$ for $\sigma \in \mathcal{C}$ partitions $\mathrm{Sim}(\mathfrak{g},K)$. Also, the collection of all $\widehat{K}(\sigma)$ for $\sigma \in \mathcal{C}$ partitions $\widehat{K}$.
\end{theorem}

Now, since our goal is to rephrase the admissible Mackey bijection \cite{Afgoustidis1}*{§5}, we also recall the classification of the $(\mathfrak{g}_0,K)$-modules, due to Champetier--Delorme \cite{DelormeCh}. See also \cite{GaudillotEstrada} and Appendix \ref{relInd} for a new proof. Let $\mathrm{Sim}(\mathfrak{g}_0, K)$ denote the set of isomorphism classes of simple $(\mathfrak{g}_0, K)$-modules.

Let $\lambda\in \mathfrak{s}^*$. We denote by $K^\lambda$ and $G_0^\lambda$ the centralizers of $\lambda$ in $K$ and $G$, respectively. Note that $G_0^\lambda = K^\lambda \ltimes \mathfrak{s}_\mathbb{R}$. In fact, if $L^\lambda$ denotes the centralizer of $\lambda$ (viewed as a linear form on $\mathfrak{g}$ vanishing on $\mathfrak{k}$) in $G$, with split component $A_\lambda$, and if $P$ is any parabolic subgroup of $G$ with Levi factor $L^\lambda$, then $\lambda \in \mathfrak{a}_\lambda^*$ and $P_0 = G^\lambda_0$. We can therefore use the notations introduced at the end of Subsection \ref{DeformationConvention}.

\begin{theorem}
    For every $\lambda \in \mathfrak{a}^* \subset \mathfrak{s}^*$ and every $K^\lambda$-type $\pi$,
    \[J^0(\lambda,\pi) = \mathrm{Ind}_{G_0^\lambda}^{G_0} (\pi \otimes \mathbb{C}_\lambda)\]
    defines an element of $\mathrm{Sim}(\mathfrak{g}_0,K)$. Two such elements $J^0(\lambda, \pi)$ and $J^0(\lambda',\pi')$ are equal if and only if $(\lambda, \pi)$ is related to $(\lambda',\pi')$ by an element of $W(\mathfrak{g},\mathfrak{a})$. Moreover, every element of $\mathrm{Sim}(\mathfrak{g}_0,K)$ is of this form, for some such $(\lambda, \pi)$.
\end{theorem}

Finally, let us introduce the Mackey bijection. It is a one-to-one map \[\mathcal{M} : \mathrm{Sim}(\mathfrak{g},K) \to \mathrm{Sim}(\mathfrak{g}_0,K),\] characterized as follows (we refer to \cite{Afgoustidis1} for more information). Let $\lambda \in \mathfrak{a}^*$ and let $\pi$ be any $K^\lambda$-type. Let $L^\lambda$ be the centralizer of $\lambda$ in $G$. There exists a unique isomorphism class $V(\pi)$ of simple, tempered $(\mathfrak{l}^\lambda, K^\lambda)$-modules with real infinitesimal character (in short, \textbf{tempiric}) and whose unique lowest $K^\lambda$-type is $\pi$, see \cite{Afgoustidis1}*{§3.1} and the references therein. Let $P$ be any parabolic subgroup of $G$ with Levi factor $L^\lambda$. Then $\mathrm{Ind}_P^G(\mathrm{V}(\pi) \otimes \mathbb{C}_\lambda)$ has a unique subquotient
\begin{equation*}
    J(\lambda,\pi) \in \mathrm{Sim}(\mathfrak{g},K)
\end{equation*}
containing all of the lowest $K$-types, see \cite{VoganGreen}*{Proposition 4.20, Theorem 6.6.15, Theorem 6.5.10(d)} and the proof of \cite{Afgoustidis1}*{Theorem 5.1}. We define
\[ \mathcal{M} \left( J(\lambda,\pi)\right) = J^0(\lambda,\pi). \]
In the following, we will need the following lemma.

\begin{lemma}\label{preservesKtypes}
    The Mackey bijection $\mathcal{M}$ preserves lowest $K$-types.
\end{lemma}

\begin{proof}
    Let $\lambda \in \mathfrak{a}^*$ and let $\pi$ be a $K^\lambda$-type. Assume that there exists a purely imaginary parameter $\nu \in i\mathfrak{a}^*_\mathbb{R}$ such that $L^\nu = L^\lambda$, or equivalently, $\mathfrak{l}^\nu = \mathfrak{l}^\lambda$. Then $J(\lambda,\pi)$ and $J(\nu, \pi)$ define the same isomorphism class of $K$-modules. The same is true for $J^0(\lambda,\pi)$ and $J^0(\nu, \pi)$.
    Since the Mackey bijection preserves lowest $K$-types for unitary representations \cite{Afgoustidis1}*{Proposition 4.1}, the set of lowest $K$-types of $J^0(\nu, \pi)$ and $J(\nu, \pi)$ are the same. Hence $J(\lambda,\pi)$ and $\mathcal{M}(J(\lambda,\pi)) = J^0(\lambda,\pi)$ have the same set of lowest $K$-types. Let us therefore check that such $\nu$ exists.
    
    Let us write $\lambda = \lambda_0 + i \lambda_1$ with $\lambda_0,\lambda_1 \in \mathfrak{a}_\mathbb{R}^*$ and let us assume that $\lambda_0 \neq 0$ (otherwise we take $\nu = \lambda$). Let $\Delta \subset \mathfrak{a}_\mathbb{R}^*$ be the root system associated to the pair $(\mathfrak{g}_\mathbb{R},\mathfrak{a}_\mathbb{R})$. Since $\Delta$ is finite, there exists $\lambda' \in \mathrm{span}_\mathbb{R}(\lambda_0,\lambda_1)$ such that
    \[(\alpha |\lambda) = 0 \Longleftrightarrow \{\lambda_0,\lambda_1\} \subset \alpha^\perp \Longleftrightarrow (\alpha | \lambda') = 0\]
    for every $\alpha \in \Delta$. Putting $\nu = i \lambda'$, we obtain
    \(\mathfrak{l}^\lambda = \oplus_{\alpha \in \lambda^\perp}\mathfrak{g}^\alpha = \oplus_{\alpha \in \nu^\perp}\mathfrak{g}^\alpha = \mathfrak{l}^\nu,\)
    which concludes the proof.
\end{proof}

Let us now use the Vogan classification in order to characterize the Mackey bijection in a slightly different way. Let $\sigma$ be a cuspidal datum of $G$ and $\lambda \in \mathfrak{a}^*_\sigma$. Let us choose any parabolic subgroup $P$ of $G$ with Levi factor $L_\sigma$. Recall that $P_0$ only depends on $L_\sigma$. Let us write
\begin{equation} \label{definition_I^0} I^0_{\sigma,\lambda} = \mathrm{Ind}_{P_0}^{G_0}\left(\mathrm{res}_{K_\sigma}^{\mathfrak{l}_\sigma,K_\sigma}(D_\sigma) \otimes \mathbb{C}_\lambda\right).\end{equation}
Then $I^0_{\sigma,\lambda}$ and $I^P_{\sigma,\lambda}$ are the specializations of $\mathbf{Ind}_P^G(D_\sigma \otimes \mathbb{C}_\lambda)$ at $t=0$ and at $t= 1$, respectively, see \ref{properties_induction}\textit{(2)}. By specializing Proposition \ref{properties_induction}\textit{(3)} at $t=0$ and at $t= 1$, it follows that $I^0_{\sigma,\lambda}$ is isomorphic to $I^{P}_{\sigma,\lambda}$ as a $K$-module. The set of lowest $K$-types of $I^0_{\sigma,\lambda}$ is therefore $\widehat{K}(\sigma)$ and all of these $K$-types occur with multiplicity one (because of Proposition \ref{multiplicityone}). For every $\mu \in \widehat{K}(\sigma)$, let
\begin{equation*}
    J^0_{\sigma,\lambda}(\mu) \in \mathrm{Sim}(\mathfrak{g}_0,K)
\end{equation*}
denote the unique simple subquotient of $I^0_{\sigma,\lambda}$ containing the $K$-type $\mu$.

\begin{proposition}\label{MackeyBijDiff}
    Let $\sigma$ be a cuspidal datum of $G$ and let $\lambda \in \mathfrak{a}^*_\sigma$. For every $\mu \in \widehat{K}(\sigma)$, the Mackey bijection $\mathcal{M}$ maps $J_{\sigma,\lambda}(\mu)$ to $J^0_{\sigma,\lambda}(\mu)$.
\end{proposition}

\begin{proof}
    Let us pick parabolic subgroups $P_\sigma$ and $P_\lambda$ of $G$ with Levi factors $L_\sigma$ and $L^\lambda$ respectively and such that $P_\sigma \subset P_\lambda$. Then $P_{\sigma,\lambda} = P_\sigma \cap L^\lambda$ is a parabolic subgroup of $L_\lambda$ and we have
    \[I^P_{\sigma,\lambda} = \mathrm{Ind}_{P_\lambda}^{G}\left( \mathrm{Ind}^{L^\lambda}_{P_{\sigma,\lambda}}(D_{\sigma}) \otimes \mathbb{C}_\lambda\right).\]
    Since $D_{\sigma}$ is a tempiric $(\mathfrak{l}_\sigma,K_\sigma)$-module, $\mathrm{Ind}^{L^\lambda}_{P_{\sigma,\lambda}}(D_{\sigma})$ decomposes into a finite direct sum $\bigoplus_\gamma X_\gamma$ of tempiric $(\mathfrak{l}^\lambda,K^\lambda)$-modules. For every index $\gamma$, let $\pi_\gamma$ be the unique lowest $K^\lambda$-type of $X_\gamma$ and let $X_\gamma^\mathrm{low}$ be the corresponding isotypical component.

    Let us now fix a $\mu \in \widehat{K}(\sigma)$. Since $\mu$ occurs in $I^P_{\sigma,\lambda}$ with multiplicity one and
    \[I^P_{\sigma,\lambda} = \bigoplus_\gamma \mathrm{Ind}_{P_\lambda}^{G}(X_{\gamma} \otimes \mathbb{C}_\lambda),\]
    there exists a unique index $\gamma_\mu$ such that $\mu$ is a $K$-type of $\mathrm{Ind}_{P_\lambda}^{G}(X_{\gamma_\mu}\otimes \mathbb{C}_\lambda)$, and $\mu$ is then a lowest $K$-type of the latter. Thus, we have $J_{\sigma,\lambda}(\mu) = J(\pi_{\gamma_\mu},\lambda)$.

    On the motion group's side, we claim that
    \begin{equation}\label{equation_mackey_bij_K_types}
        I^0_{\sigma,\lambda} = \mathrm{Ind}_{G_0^\lambda}^{G_0}\!\left(\mathrm{ind}^{K^\lambda}_{K_\sigma}\!\left(\mathrm{res}_{K_\sigma}^{\mathfrak{l}_\sigma,K_\sigma}(D_\sigma)\right) \otimes \mathbb{C}_\lambda \right).
    \end{equation}
    Indeed, multiplication induces an isomorphism $U(K^\lambda) \otimes_\mathbb{C}U(\mathfrak{s}) \to U(\mathfrak{g}_0^\lambda,K^\lambda)$ of algebras, which in turn induces a natural isomorphism of $(\mathfrak{g}_0^\lambda,K^\lambda)$-modules
    \begin{equation*}
    \mathrm{ind}^{K^\lambda}_{K_\sigma}\!\left(\mathrm{res}_{K_\sigma}^{\mathfrak{l}_\sigma,K_\sigma}(D_\sigma)\right) \otimes \mathbb{C}_\lambda = \mathrm{ind}^{\mathfrak{g}_0^\lambda,K^\lambda}_{\mathfrak{k}_\sigma \ltimes \mathfrak{s}_\mathbb{R},K_\sigma}\!\left( \mathrm{res}_{K_\sigma}^{\mathfrak{l}_\sigma,K_\sigma}(D_\sigma)\otimes \mathbb{C}_\lambda\right).\end{equation*}
    The above identification (\ref{equation_mackey_bij_K_types}) is then obtained by induction in stages, see Proposition \ref{properties_ind_res}\textit{(2)}.
    
    Since the restriction of $\mathrm{Ind}^{L^\lambda}_{P_{\sigma,\lambda}}(D_{\sigma})$ to $K^\lambda$ is isomorphic to $\mathrm{ind}^{K^\lambda}_{K_\sigma}(\mathrm{res}_{K_\sigma}^{\mathfrak{l}_\sigma,K_\sigma}(D_\sigma))$, we obtain from (\ref{equation_mackey_bij_K_types}) the following:
    \[I^0_{\sigma,\lambda} =\bigoplus_\gamma\mathrm{Ind}_{G_0^\lambda}^{G_0}\!\left(\mathrm{res}_{K^\lambda}^{\mathfrak{l}_\lambda,K_\lambda}(X_\gamma) \otimes \mathbb{C}_\lambda \right).\]
     In particular, we have an embedding of $(\mathfrak{g}_0,K)$-modules
     \[\mathrm{Ind}_{G_0^\lambda}^{G_0}(X^\mathrm{low}_{\gamma_\mu} \otimes\mathbb{C}_\lambda) \hookrightarrow \mathrm{Ind}_{G_0^\lambda}^{G_0}\!\left(\mathrm{res}_{K^\lambda}^{\mathfrak{l}_\lambda,K_\lambda}(X_\gamma) \otimes \mathbb{C}_\lambda\right) \subset I^0_{\sigma,\lambda}.\]
     Now, the isomorphism class of $\mathrm{Ind}_{G_0^\lambda}^{G_0}(X^\mathrm{low}_{\gamma_\mu} \otimes\mathbb{C}_\lambda)$ is $J^0(\pi_{\gamma_\mu},\lambda)$ and the set of lowest $K$-types of $J^0(\pi_{\gamma_\mu},\lambda)$, being the same as that of $J(\pi_{\gamma_\mu},\lambda) = J_{\sigma,\lambda}(\mu)$ by Lemma \ref{preservesKtypes}, contains $\mu$. In conclusion, we have $J^0_{\sigma,\lambda}(\mu) = J^0(\pi_{\gamma_\mu},\lambda) = \mathcal{M}(J_{\sigma,\lambda}(\mu))$.
\end{proof}

\subsection{Harish-Chandra homomorphisms along the Mackey deformation}\label{HC_homomorphisms}

We now come to the very goal of this paper. In what follows, we algebraically deform Delorme's generalized Harish-Chandra homomorphisms \cite{Delorme84} along the Mackey parameter $t$ and we prove that the $\mathbb{C}[t]$-algebra homomorphism thus defined is epispectral. We consider the setting defined by Subsection \ref{DeformationConvention}. In what follows, we will often use certain basic results on approximately unital algebras from Appendix~\ref{AUA}.

We begin with some conventions that we will use for the rest of the section, as well as some general remarks. For every $z \in \mathbb{C}$, we write
    \[U^z = U(\bm{\mathfrak{g}}_z, K) = \begin{cases}
        U(\mathfrak{g},K) & (z \neq 0),\\
        U(\mathfrak{g}_0,K) & (z = 0).
    \end{cases}\]
Moreover, we simply denote $U(\bm{\mathfrak{g}},K)$ and $U(\mathfrak{g},K)$ by $\bm{U}$ and $U$, respectively. Let us write $\mathrm{Sim}(\bm{\mathfrak{g}},K)$ for the set of isomorphism classes of simple $(\bm{\mathfrak{g}},K)$-modules. In what follows, we always identify $\mathrm{Sim}(\bm{U})$ and $\mathrm{Sim}(U^z)$ with $\mathrm{Sim}(\bm{\mathfrak{g}},K)$ and $\mathrm{Sim}(\bm{\mathfrak{g}}_z,K)$, respectively, according to Lemma \ref{correspondence_Umod_gKmod}. For every $z \in \mathbb{C}$, we denote by $\mathrm{sp}_z$ the map $\mathrm{Sim}(\bm{\mathfrak{g}}_z,K) \to \mathrm{Sim}(\bm{\mathfrak{g}},K)$ induced by pullback from the quotient morphism $\bm{U} \to U^z$.

Since $\bm{U}$ has countable $\mathbb{C}$-dimension (see the proof of Corollary \ref{admissibility_simple}), it follows from Schur's lemma \cite{Bourbaki}*{§3, N°2} that $\mathbb{C}[t]$ acts on each simple $(\bm{\mathfrak{g}},K)$-module through evaluation at some point $z \in \mathbb{C}$. We thus have
\begin{equation}\label{fibers_Sim}\mathrm{Sim}(\bm{\mathfrak{g}},K) = \bigsqcup_{z \in \mathbb{C}} \mathrm{sp}_z^*\big( \mathrm{Sim}(\bm{\mathfrak{g}}_z,K)\big).\end{equation}
For every cuspidal datum $\sigma$, let us write $\mathrm{Sim}(\mathfrak{g}_0,K)_\sigma  = \mathcal{M}(\mathrm{Sim}(\mathfrak{g},K)_\sigma)$ and
\begin{equation}\label{notation_sim_sigma}\mathrm{Sim}(\bm{\mathfrak{g}},K)_\sigma = \bigcup_{z \in \mathbb{C}} \mathrm{sp}_z^*(\mathrm{Sim}(\bm{\mathfrak{g}}_z,K)_\sigma).\end{equation}
From Theorem \ref{partition} and Lemma \ref{preservesKtypes}, it is clear that $\mathrm{Sim}(\bm{\mathfrak{g}},K)_\sigma$ consists of the isomorphism classes of simple $(\bm{\mathfrak{g}},K)$-modules whose set of lowest $K$-types is contained in $\widehat{K}(\sigma)$, or equivalently, intersects $\widehat{K}(\sigma)$.

Let us now introduce the setting of the generalized Harish-Chandra homomorphisms along the Mackey deformation. Let $\sigma$ be a cuspidal datum and let $P$ be a parabolic subgroup of $G$ with Levi factor $L_\sigma$. We simply denote $S(\mathfrak{a}_\sigma)$ by $S$. The action of $\mathfrak{a}_\sigma$ on $S$ by multiplication endows the latter with an $\mathfrak{a}_\sigma$-module structure. The natural $S$-module structure of $S$ is compatible with with this $\mathfrak{a}_\sigma$-module structure, so that we view $S$ as an $\mathfrak{a}_{L,S}$-module (recall Notation \ref{notation_specialization_HC}). We introduce the $(\bm{\mathfrak{g}}_S,K)$-module
\[\bm{I}_\sigma = \mathbf{Ind}_{P,S}^G(D_\sigma \otimes S).\]
We really consider the latter as a $(\bm{\mathfrak{g}},K)$-module with a compatible $S$-module structure, see Remark \ref{compatible}. Let us write $\chi_\sigma = \chi_{\widehat{K}(\sigma)}$. The main aim of this subsection is to describe the image of $\chi_{\sigma}\bm{U}\chi_\sigma$ through its representation on $\bm{I}_\sigma$, up to isospectral extension.

In order to state our theorem, we need some preparation. Let $\bm{\pi}_\sigma$ be the natural representation of $\Upsilon(\bm{\mathfrak{g}},K)$ on $\bm{I}_\sigma$. For each $z \in \mathbb{C}$, let $I_\sigma^z$ denote the $(\bm{\mathfrak{g}}_{z,S},K)$-module obtained from $\bm{I}_\sigma$ by specialization at $t = z$ and let $\pi_\sigma^z$ be the associated representation of $\Upsilon(\bm{\mathfrak{g}}_z,K)$ on the underlying $(\bm{\mathfrak{g}}_{z},K)$-module. For simplicity, we write $I_\sigma$ and $\pi_\sigma$ for $I_\sigma^1$ and $\pi_\sigma^1$, respectively.

Let us write
\[I_{\sigma|K} = \mathrm{ind}_{K\cap L}^{K}(\mathrm{res}_{K\cap L}^{\mathfrak{l}, K\cap L} D_\sigma).\]
The restriction of $\bm{I}_\sigma$ to a $K_{S[t]}$-module is given by Proposition \ref{properties_induction}\textit{(3)}:
\begin{equation}\label{compact_piucture_Isigma}\mathrm{res}_{K_{S[t]}}^{\bm{\mathfrak{g}}_S, K}(\bm{I}_\sigma) = I_{\sigma|K}  \otimes_\mathbb{C} S[t].\end{equation} Since $I_{\sigma|K}$ is the restriction to $K$ of any $I^P_{\sigma,\lambda}$ for $\lambda \in \mathfrak{a}_\sigma^*$ (see the previous subsection), its set of lowest $K$-types is $\widehat{K}(\sigma)$ and each of these $K$-types occurs with multiplicity one (Proposition \ref{multiplicityone}). For every $\mu \in \widehat{K}(\sigma)$, let $E_\mu$ denote the isotypic component of $I_{\sigma|K}$ of type $\mu$, which is thus an irreducible representation of $K$. We also write
\[E_\sigma  = \oplus_{\mu \in \widehat{K}(\sigma)} E_\mu = \chi_\sigma (I_{\sigma|K}).\]
Then, as a $K_{S[t]}$-module, $\chi_\sigma \bm{I}_\sigma$ identifies with $E_\sigma \otimes_\mathbb{C} S[t]$, while $\chi_\sigma \mathrm{End}_{S[t]}(\bm{I}_\sigma) \chi_\sigma$ identifies with $\mathrm{End}_\mathbb{C}(E_\sigma) \otimes_\mathbb{C} S[t]$. Accordingly, for every $z \in \mathbb{C}$, we can identify $\chi_\sigma \mathrm{End}_{S}(I^z_\sigma) \chi_\sigma$ with $\mathrm{End}_\mathbb{C}(E_\sigma) \otimes_\mathbb{C} S$.

For every pair $(\mu,\nu)$ of elements of $\widehat{K}(\sigma)$, let $r(\mu,\nu)$ be the unique element of $\hat{R}_\sigma$ such that $\mu = r(\mu,\nu) \nu$. We consider the representation of $W_\sigma$ on $\mathrm{End}_\mathbb{C}(E_\sigma)$ defined as follows: for every $(\mu,\nu)\in \widehat{K}(\sigma)^2$, the group $W_\sigma$ acts on $\mathrm{Hom}_\mathbb{C}(E_\nu, E_\mu)$ through $r(\mu,\nu)$, which we consider as a character of $W_\sigma$ via pullback by the group homomorphism $W_\sigma \to R_\sigma$. On the other hand, the contragredient of the natural representation of $W_\sigma$ on $\mathfrak{a}_\sigma^*$ extends to a representation of $W_\sigma$ on $S$.

\begin{definition}\label{def_Bsigma}
    Considering the tensor product of the two above representations of $W_\sigma$, we define
    \[B_\sigma = \left(\mathrm{End}_\mathbb{C}(E_\sigma) \otimes_\mathbb{C} S \right)^{W\sigma}.\]

\end{definition}
For every $T \in \{\bm{U}, U, U^0\}$, let us write $T_\sigma = \chi_\sigma T\chi_\sigma$. In \cite{Delorme84}, Delorme proves the following theorem, although it is stated a bit differently.

\begin{theorem}[Delorme \cite{Delorme84}*{Théorème 3}]\label{Delorme_theorem} The image of $U_\sigma$ through $\pi_\sigma$ is $B_\sigma$.
\end{theorem}

Delorme's theorem was originally stated in the connected semisimple case, but it extends to general real reductive groups, see \cite{DelormeAnnals}*{Remarque 1}. In this subsection, our goal is to prove a weak analogue of this theorem but for the whole Mackey deformation.

\begin{theorem}\label{epispectral}
    We have $\bm{\pi}_\sigma(\bm{U}_\sigma) \subset B_\sigma[t]$ and $\bm{\pi}_\sigma : \bm{U}_\sigma \to B_\sigma[t]$ is epispectral.
\end{theorem}

\begin{proof}
    We divide the proof of this theorem into four steps.

\smallskip

\textit{Step 1: the inclusion.}
    First, we prove that $\bm{\pi}_\sigma(\bm{U}_\sigma) \subset B_\sigma[t]$. For every $z \in \mathbb{C}$, we have a commutative diagram
    \[\begin{tikzcd}
        \bm{U}_\sigma \ar[r, "\bm{\pi}_\sigma"] \ar[d] & \mathrm{End}_\mathbb{C}(E_\sigma) \otimes_\mathbb{C} S[t] \ar[d, "\mathrm{sp}_z"] \ar[r, hookrightarrow] & \mathrm{End}_{S[t]}(\bm{I}_\sigma) \ar[d] \\
        U^z_\sigma \ar[r, "\pi_\sigma^z"] &\mathrm{End}_\mathbb{C}(E_\sigma) \otimes_\mathbb{C} S \ar[r, hookrightarrow] & \mathrm{End}_S(I_\sigma^z)
    \end{tikzcd},\]
    where the left vertical arrow is the restriction of the specialization map $\bm{U} \to U^z$ and the right vertical arrow is induced by the functor of specialization at $z$ defined on the category of $S[t]$-modules.
    This diagram shows that it is enough to prove that $\bm{\pi}^z_\sigma(U^z_\sigma) \subset B_\sigma$ for every $z$ in a Zariski dense subset of $\mathbb{C}$ (such as $\mathbb{C}^\times$) in order to obtain the inclusion $\bm{\pi}_\sigma(\bm{U}_\sigma) \subset B_\sigma[t]$.

    So let us fix $z \in \mathbb{C}^\times$. Recall that $U^z = U$. By Proposition \ref{properties_induction}\textit{(2)}, we have
    \[I_\sigma^z = \mathrm{Ind}_{P,S}^G(D_\sigma \otimes \mathrm{res}_{\delta_z}S),\]
    where $\delta_z$ is the automorphism of $\mathfrak{a}_{L,S}$ defined by dilation by $z^{-1}$. The linear automorphism $\delta'_z$ of $S$ sending any homogeneous element $x \in S$ of degree $n$ to $z^{-n}x$ defines an isomorphism of $\mathfrak{a}_L$-modules $S \to \mathrm{res}_{\delta_z}S$. It thus induces an isomorphism of $(\mathfrak{g},K)$-modules
    \[\tilde{\delta}'_z : \mathrm{Ind}_P^G(D_\sigma \otimes S) \to \mathrm{Ind}_P^G(D_\sigma  \otimes \mathrm{res}_{\delta_{z}}S),\]
    that is, given Remark \ref{commute_parInd_for}, an isomorphism of $(\mathfrak{g},K)$-modules $I_\sigma \to I_\sigma^z$. The restriction of $\tilde\delta'_z$ to a linear map $\chi_\sigma I_\sigma \to \chi_\sigma I_\sigma^z$ identifies with \(\mathrm{id} \otimes \delta'_z : E_\sigma \otimes_\mathbb{C} S \to E_\sigma \otimes_\mathbb{C} S\).
    It follows that $\pi_\sigma^z(U_\sigma) = (\mathrm{id} \otimes \delta'_z) \pi_\sigma(U_\sigma)(\mathrm{id} \otimes \delta'_z)^{-1}$. Now the map from $\mathrm{End}_S(\chi_\sigma I_\sigma)$ to $\mathrm{End}_S(\chi_\sigma I_\sigma^z)$ defined by conjugation by $(\mathrm{id} \otimes \delta'_z)$ identifies with \begin{equation}\label{rescaling}\mathrm{id} \otimes \delta'_z : \mathrm{End}_\mathbb{C}(E_\sigma) \otimes_\mathbb{C} S \to \mathrm{End}_\mathbb{C}(E_\sigma) \otimes_\mathbb{C} S.\end{equation}
    Therefore, using Theorem \ref{Delorme_theorem}, we have $\pi^z_\sigma(U^z_\sigma) = (\mathrm{id} \otimes \delta'_z)(B_\sigma) = B_\sigma$, which concludes this first step.
    
\smallskip

\textit{Step 2 : reduction to the fiber at $t = 0$ for the epispectral property.} Let $z \in \mathbb{C}$. Since $\bm{U}_\sigma$ is a corner subalgebra of $\bm{U}$, it is in particular a direct summand of the $\mathbb{C}[t]$-module $\bm{U}$. The specialization of $\bm{U}_\sigma$ at $z$ therefore coincides with the image of $\bm{U}_\sigma$ through the specialization map $\bm{U} \to U^z$, that is, $U^z_\sigma$. It follows that $\pi_\sigma^z : U_\sigma^z \to B_\sigma$ is precisely the morphism of $\mathbb{C}$-algebras induced from $\bm{\pi}_\sigma$ by specialization at $z$. Since $\bm{U}$ has countable $\mathbb{C}$-dimension, so does $\bm{U}_\sigma$ and by Lemma~\ref{isospectral_specializations}, $\bm{\pi}_\sigma : \bm{U}_\sigma \to B_\sigma[t]$ is epispectral if $\pi_\sigma^z : U_\sigma^z \to B_\sigma$ is epispectral for every $z \in \mathbb{C}$. For $z \in \mathbb{C}^\times$, this is clear: we proved in the first step of the proof that $\pi_\sigma^z : U^z_\sigma \to B_\sigma$ is surjective. We are thus left to prove that $\pi_\sigma^0 : U_\sigma^0 \to B_\sigma$ is epispectral.

\smallskip

\textit{Step 3: the simple spectrum of $B_\sigma$}. Let us now describe the simple spectrum of~$B_\sigma$. We could obtain such a description by hand (see Remark \ref{morita_equivalence} below) but, since we have Theorem~\ref{Delorme_theorem} and $\pi_\sigma(U_\sigma) = B_\sigma$, we rather use the representation theory of~$G$. The simple spectrum of $B_\sigma$ identifies via pullback by $\pi_\sigma$ to the isomorphism classes of simple $U_\sigma$-modules for which the action of $U_\sigma$ factors through~$\pi_\sigma$. On the other hand, Lemma \ref{TransferSimple} states that the simple $U_\sigma$-modules correspond, up to isomorphism, to the simple $U$-modules on which $\chi_\sigma$ acts non-trivially, that is, the simple $(\mathfrak{g},K)$-modules containing at least one of the $K$-types of $\widehat{K}(\sigma)$. We claim that, under this correspondence, we have $\pi_\sigma^*\big(\mathrm{Sim}(B_\sigma)\big) = \mathrm{Sim}(\mathfrak{g},K)_\sigma$, the latter being defined by (\ref{notation_SgKsigma}).

If $\mu$ is a $K$-type of strictly smaller length than the elements of $\widehat{K}(\sigma)$, then $\pi_\sigma(\chi_\mu) = 0$. Consequently, each of the elements of $\mathrm{Sim}(\mathfrak{g},K)$ corresponding to simple $B_\sigma$-modules via $\pi_\sigma^*$ must contain an element of $\widehat{K}(\sigma)$ as a lowest $K$-type. The inclusion $\pi_\sigma^*\big(\mathrm{Sim}(B_\sigma)\big) \subset \mathrm{Sim}(\mathfrak{g},K)_\sigma$ then follows from Theorem \ref{partition}.

Let us now prove the reverse inclusion. Let $\lambda \in \mathfrak{a}_\sigma^*$ and $\mu \in \widehat{K}(\sigma)$. Recall (\ref{definition_I^P}). Let $V_1 \subset V_2$ be $(\mathfrak{g},K)$-submodules of $I^P_{\sigma,\lambda}$ such that the isomorphism class of $V_2/V_1$ is $J_{\sigma,\lambda}(\mu)$. According to Theorem \ref{vogan_classification}, the set of lowest $K$-types of $V_2/V_1$ is $\hat{R}_\sigma(\lambda)\mu$. Given Proposition \ref{multiplicityone}, writing $\chi_{\mu,\lambda} = \chi_{\hat{R}_\sigma(\lambda)\mu}$, we thus have
\[\chi_{\mu,\lambda} V_2 = \chi_{\mu,\lambda} I^P_{\sigma,\lambda}, \quad \chi_{\mu,\lambda} V_1 = 0, \quad \chi_{\mu,\lambda} (V_2/V_1) = \chi_\sigma(V_2/V_1),\]
which proves that the quotient map $V_2 \to V_2/V_1$ induces a linear isomorphism $\chi_{\mu,\lambda} I^P_{\sigma,\lambda} \to \chi_\sigma(V_2/V_1)$. This implies that $\chi_{\mu,\lambda} I^P_{\sigma,\lambda}$ is a simple $U_\sigma$-submodule of $\chi_\sigma I^P_{\sigma,\lambda}$ which is of course isomorphic to $\chi_\sigma(V_2/V_1)$. Now, according to Remark \ref{specialization_lambda}, $I^P_{\sigma,\lambda}$ is the specialization of $I_\sigma$ at $\lambda$, the latter being viewed as a character of $S$. It follows that the action of $U_\sigma$ on $\chi_\sigma(V_2/V_1)$ factors through $\pi_\sigma$. This proves the converse inclusion $\mathrm{Sim}(\mathfrak{g},K)_\sigma \subset \pi_\sigma^*\big(\mathrm{Sim}(B_\sigma)\big)$.

\smallskip

\textit{Step 4: the specialization at $t = 0$ is isospectral.} By the very definition of $B_\sigma$, the above description of $\mathrm{Sim}(B_\sigma)$ only involves the $K_S$-module structure of $I_\sigma$. The latter is the specialization at $t=1$ of $\bm{I}_\sigma$, which identifies as a $K_{S[t]}$-module with $I_{|K} \otimes_\mathbb{C} S[t]$, see (\ref{compact_piucture_Isigma}). As a $K_S$-module, $I_\sigma$ can thus be identified with $I_\sigma^0$, the specialization of $\bm{I}_\sigma$ at $t=0$. This way, viewing $B_\sigma$ as a subalgebra of $\chi_\sigma \mathrm{End}_{S}(I^0_\sigma) \chi_\sigma$ rather than as a subalgebra of $\chi_\sigma \mathrm{End}_{S}(I_\sigma) \chi_\sigma$, we can re-express the classification of simple $B_\sigma$-modules using $I_\sigma^0$ in place $I_\sigma$, as follows.

Recall that, according to Remark \ref{specialization_lambda}, the specialization of $I_\sigma^0$ at any $\lambda \in \mathfrak{a}_\sigma^*$ is $I^0_{\sigma,\lambda}$ (see (\ref{definition_I^0}) for the notation). For every $\lambda \in \mathfrak{a}_\sigma^*$ and every $\mu \in \widehat{K}(\sigma)$, the action of $B_\sigma$ on $\chi_{\sigma} I^0_{\sigma,\lambda}$, obtained from the natural action on $\chi_\sigma I_\sigma^0 = E_\sigma$ by specialization at $\lambda$, stabilize $\chi_{\mu,\lambda} I^0_{\sigma,\lambda}$ and the latter is a simple $B_\sigma$-module. Moreover, two such $B_\sigma$-modules $\chi_{\mu,\lambda} I^0_{\sigma,\lambda}$ and $\chi_{\sigma,\mu'} I^0_{\sigma,\lambda'}$ are isomorphic if and only if $\chi_{\mu,\lambda} I^P_{\sigma,\lambda}$ and $\chi_{\sigma,\mu'} I^P_{\sigma,\lambda'}$ are isomorphic as $U_\sigma$-modules, that is, if and only if $J_{\sigma,\lambda}(\mu) = J_{\sigma,\lambda'}(\mu')$.

On the other hand, the simple spectrum of $U_\sigma^0$ corresponds via Lemma \ref{TransferSimple} to the subset of $\mathrm{Sim}(\mathfrak{g}_0,K)$ consisting of simple modules containing a $K$-type in $\widehat{K}(\sigma)$. Therefore, using Proposition \ref{MackeyBijDiff} and Lemma \ref{preservesKtypes} combined with Theorems \ref{vogan_classification} and \ref{partition}, imitating the two last paragraphs of \textit{step 3} with $U_\sigma^0$, $I^0_{\sigma,\lambda}$, $J^0_{\sigma,\lambda}(\mu)$ instead of $U_\sigma$, $I^P_{\sigma,\lambda}$, $J_{\sigma,\lambda}(\mu)$, respectively, we obtain that $\pi_\sigma^0 : U_\sigma^0 \to B_\sigma$ is epispectral with $(\pi_\sigma^0)^*(\mathrm{Sim}(B_\sigma)) = \mathrm{Sim}(\mathfrak{g}_0,K)_\sigma$. This concludes the proof of the theorem.
\end{proof}

\begin{remark}\label{facts_after_proof} In the course of the above proof, we established the following facts.
\begin{enumerate}
    \item For every $z \in \mathbb{C}$, the specialization at $t= z$ of $\bm{\pi}_\sigma : \bm{U}_\sigma \to B_n[t]$ is the well-defined epispectral morphism $\pi_\sigma^z : U_\sigma^z \to B_\sigma$.
    \item For every $z \in \mathbb{C}^\times$, the automorphism $\mathrm{id} \otimes \delta'_z$ of $B_\sigma$ defined by Formula (\ref{rescaling}) satisfies $\pi_\sigma^z = (\mathrm{id} \otimes \delta'_z)\circ \pi_\sigma$.
    \item Let us identify the simple spectra of $U_\sigma$ and $U_\sigma^0$ with subsets of $\mathrm{Sim}(\mathfrak{g},K)$ and $\mathrm{Sim}(\mathfrak{g}_0,K)$, respectively, via Lemma \ref{TransferSimple}. The diagram
    \[\begin{tikzcd}
        \mathrm{Sim}(\mathfrak{g},K) \ar[r, "\mathcal{M}"] & \mathrm{Sim}(\mathfrak{g}_0,K) \\
        \mathrm{Sim}(B_\sigma) \ar[u, "(\pi_\sigma^1)^*"', hookrightarrow] \ar[r, equals] & \mathrm{Sim}(B_\sigma) \ar[u, "(\pi_\sigma^0)^*"', hookrightarrow]
    \end{tikzcd}\]
    commutes and we have $(\pi_\sigma^1)^*\big(\mathrm{Sim}(B_\sigma) \big) = \mathrm{Sim}(\mathfrak{g},K)_\sigma$.
\end{enumerate}
Moreover, combining the above two last facts, we obtain the following result.
    \begin{enumerate}[resume]
        \item Identifying the simple spectrum of $\bm{U}_\sigma$ with a subset of $\mathrm{Sim}(\bm{\mathfrak{g}},K)$ via Lemma~\ref{TransferSimple}, we have
    \[\bm{\pi}_\sigma^*\big(\mathrm{Sim}(B_\sigma[t])\big) = \mathrm{Sim}(\bm{\mathfrak{g}},K)_\sigma,\]
    the latter being defined by (\ref{notation_sim_sigma}).
    \end{enumerate}
Finally, from the second fact and the description in \textit{Step 3} of the simple spectrum of $B_\sigma$, we deduce the following fact.
    \begin{enumerate}[resume]
        \item Identifying the simple spectrum of $U_\sigma$ with a subset of $\mathrm{Sim}(\mathfrak{g},K)$ via Lemma \ref{TransferSimple}, for every $z \in \mathbb{C}$, we have a commutative diagram
        \[\begin{tikzcd}
        \mathrm{Sim}(\mathfrak{g},K)_\sigma \ar[r, "\Lambda_z"] & \mathrm{Sim}(\mathfrak{g},K)_\sigma \\
        \mathrm{Sim}(B_\sigma) \ar[u, "(\pi_\sigma^1)^*"] \ar[ur, "(\pi_\sigma^z)^*", bend right]&
    \end{tikzcd},\]
    where $\Lambda_z$ is defined by
    \(J_{\sigma,\lambda}(\mu) \mapsto J_{\sigma,z^{-1}\lambda}(\mu)\).
    \end{enumerate}
\end{remark}

\begin{remark}\label{morita_equivalence}
    For a given finite group $\Gamma$ acting by automorphisms on a unital ring $R$, the skew group ring $\Gamma \ltimes R$ is defined as a universal element of the functor associating to any unital ring $B$ the set of pairs $(\varphi : \Gamma \to B^\times, \phi : R \to B)$ of group and unital ring homomorphisms, respectively, such that $\varphi(\gamma)\phi(r) = \phi(\gamma \cdot r)\varphi(\gamma)$ for every $(\gamma, r) \in \Gamma \times R$.

    The algebra $B_\sigma$ is in fact Morita equivalent to the skew group ring $R_\sigma \ltimes \mathbb{C}[\mathfrak{a}^*_\sigma/W_\sigma']$, where $\mathbb{C}[\mathfrak{a}^*_\sigma/W_\sigma'] = S^{W_\sigma'}$ is the algebra of regular functions on the affine variety $\mathfrak{a}^*_\sigma/W_\sigma'$. Thus, in the proof of Theorem \ref{epispectral}, we could also have used \cite{GaudillotEstrada}*{Theorem 5} to describe the simple spectrum of $B_\sigma$.

    Let us briefly see how to construct such a Morita equivalence. Let $\mu \in \widehat{K}(\sigma)$. Let us choose rank one elements $e_{\rho,r} \in \mathrm{End}_\mathbb{C}(E_{r\mu}, E_{\rho \mu})$ for all $r,\rho \in \hat{R}_\sigma$, satisfying
    \[e_{\tau,\rho} e_{\rho,r} = e_{\tau,r} \qquad (\tau, \rho,r \in \hat{R}_\sigma).\]
    If $e = \sum_r e_{r,r}$, then $e B_\sigma e$ is isomorphic to $R_\sigma \ltimes \mathbb{C}[\mathfrak{a}^*_\sigma/W_\sigma']$ via the ring and group homomorphisms
    \begin{gather*}
    \begin{cases}
        S^{W_\sigma'} \to e B_\sigma e\\
        f \mapsto \textstyle\sum_{\rho,r} e_{\rho r, r} \otimes f^\rho
    \end{cases},\qquad
    \begin{cases}
        R_\sigma \to (e B_\sigma e)^\times\\ w\mapsto \textstyle\sum_r r(w)e_{r,r}
    \end{cases},
    \end{gather*}
    where for every $f \in S^{W_\sigma'}$ and $\rho \in \hat{R}_\sigma$, we are denoting by $f^\rho$ the $R_\sigma$-isotypic component of $f$ of type $\rho$.
\end{remark}

Later we will need to apply Theorem \ref{baum_nistor_consequences} to $\bm{\pi}_\sigma$, so we now check that the associated finite generation condition is satisfied. Recall from Subsection \ref{Center} the definition of the central subalgebra $\mathcal{Z}(\bm{\mathfrak{g}})$ of $U(\bm{\mathfrak{g}})$. Since the elements of $\mathcal{Z}(\bm{\mathfrak{g}})$ commute with that of $U(K)$ and since $E_\mu$ is irreducible for each $\mu \in \widehat{K}(\sigma)$, the action of $\mathcal{Z}(\bm{\mathfrak{g}})$ on $\mathrm{End}_{S[t]}(\bm{I}_\sigma)$ by multiplication on the left through $\bm{\pi}_\sigma$ endows $\chi_\sigma\mathrm{End}_{S[t]}(\bm{I}_\sigma) \chi_\sigma= \mathrm{End}_\mathbb{C}(E_\sigma) \otimes_\mathbb{C} S[t]$ with a $\mathcal{Z}(\bm{\mathfrak{g}})$-algebra structure. Now we have
\[\bm{\pi}_\sigma(\mathcal{Z}(\bm{\mathfrak{g}})) B_\sigma = \bm{\pi}_\sigma(\mathcal{Z}(\bm{\mathfrak{g}})\chi_\sigma)B_\sigma \subset \pi_\sigma(\bm{U}_\sigma)B_\sigma \subset B_\sigma,\]
so $B_\sigma$ is a $\mathcal{Z}(\bm{\mathfrak{g}})$-subalgebra of $\mathrm{End}_\mathbb{C}(E_\sigma) \otimes_\mathbb{C} S[t]$.

\begin{proposition}\label{Z(g)_linear}
    The morphism $\bm{\pi}_\sigma : \bm{U}_\sigma \to B_\sigma[t]$ is $\mathcal{Z}(\bm{\mathfrak{g}})$-linear. Moreover, $B_\sigma[t]$ is finitely generated as a $\mathcal{Z}(\bm{\mathfrak{g}})$-module.
\end{proposition}

\begin{proof}
    The fact that $\bm{\pi}_\sigma$ is $\mathcal{Z}(\bm{\mathfrak{g}})$-linear is clear. Let us prove that $B_\sigma[t]$ is a finitely generated $\mathcal{Z}(\bm{\mathfrak{g}})$-module. The idea is to use Proposition \ref{action_center_induced_modules} to describe the action of $\mathcal{Z}(\mathfrak{g})$ on $\bm{I}_\sigma$, which, as a $(\bm{\mathfrak{g}},K)$-module, is $\mathbf{Ind}_{P}^G(D_\sigma \otimes S)$, see Remark \ref{commute_parInd_for}
    
    For that, we consider a Cartan subalgebra $\mathfrak{h}$ of $\mathfrak{g}$ obtained as the complexification of a Cartan subalgebra $\mathfrak{h}_\mathbb{R}$ of $\mathfrak{g}_\mathbb{R}$ containing $\mathfrak{a}_\mathbb{R}$, as well as a positive root system associated to $(\mathfrak{g},\mathfrak{h})$ which is compatible with the choice of $P$, as in Subsection~\ref{Center}. Let $\mathfrak{h}_\sigma = \mathfrak{m}_\sigma \cap \mathfrak{h}$ and let us identify $U(\mathfrak{h}[t])$ with $U(\mathfrak{h}_\sigma[t]) \otimes_{\mathbb{C}}U(\mathfrak{a}_\sigma[t])$ via multiplication. Let $\mathrm{Ch}_\sigma : U(\mathfrak{h}_\sigma) \to \mathbb{C}$ be a character through which the center of $U(\mathfrak{m}_\sigma)$ acts on $D_\sigma$ via the Harish-Chandra isomorphism. We extend it by $\mathbb{C}[t]$-linearity to a morphism $U(\mathfrak{h}_\sigma[t]) \to \mathbb{C}[t]$. Identifying $U(\mathfrak{a}_\sigma[t])$ with $S[t]$, we consider the $\mathbb{C}[t]$-algebra homomorphism $\mathbf{Ch}_\sigma : U(\mathfrak{h}[t]) \to S[t]$ defined as $\mathrm{Ch}_\sigma \otimes \mathrm{id}_{S[t]}$.

    According to \ref{action_center_induced_modules}, expressing $\bm{I}_\sigma$ in the compact picture (\ref{compact_piucture_Isigma}), the action of $ \mathcal{Z}(\mathbf{g})$ on $\bm{I}_\sigma$ is given by
    \[\bm{\pi}_\sigma(z) = \mathrm{id} \otimes (\delta_t\circ\mathbf{Ch}_\sigma\circ\gamma_{\mathfrak{g},\mathfrak{h}})(z) \qquad (z \in \mathcal{Z}(\bm{\mathfrak{g}})).\]
    The morphism $\delta_t : U(t\mathfrak{h}[t]) \to S[t]$ is surjective and since $U(t\mathfrak{h}[t])$ is integral over $\gamma_{\mathfrak{g},\mathfrak{h}}(\mathcal{Z}(\bm{\mathfrak{g}})) = U(t\mathfrak{h}[t])^{W(\mathfrak{g},\mathfrak{h})}$, see (\ref{extended_HC_iso}) and \cite{AtiyahMacDonald}*{Chapter 5 Exercice 12}, the ring $S[t]$ is integral over $(\delta_t\circ\mathbf{Ch}_\sigma\circ\gamma_{\mathfrak{g},\mathfrak{h}})(\mathcal{Z}(\bm{\mathfrak{g}}))$. This implies that $\mathrm{End}_\mathbb{C}(E_\sigma) \otimes_\mathbb{C} S[t]$ is a finitely generated $\mathcal{Z}(\bm{\mathfrak{g}})$-module. Since $\mathcal{Z}(\bm{\mathfrak{g}})$ is finitely generated as a $\mathbb{C}$-algebra, see the paragraph after (\ref{extended_HC_iso}), it is in particular noetherian and it follows that $B_\sigma[t]$ is finitely generated as a $\mathcal{Z}(\bm{\mathfrak{g}})$-module.
    
\end{proof}

\subsection{The main theorems}\label{main_theorems}

To conclude, we combine all the results from the previous subsections and prove the theorems announced in the introduction.

Let $(l_n)_{n \in \mathbb{N}}$ be the collection of all possible lengths of $K$-types, listed in ascending order. For every $n \in \mathbb{N}$, we denote by $\chi_n \in U(K)$ the projection onto the isotypical components of the $K$-types having length $\leq l_n$ and we choose a complete set of representatives $\mathcal{C}_n$ for the equivalence classes of cuspidal data $\sigma$ whose associated $K$-types, the elements of $\widehat{K}(\sigma)$, have length $l_n$. For any integers $n,m\geq -1$ , we write $\bm{U}_n = \chi_n \bm{U} \chi_n$, $\bm{J}_n = \bm{U}\chi_n \bm{U}$, $\bm{J}_{n,m} = \chi_m \bm{J}_n \chi_m$ and $U^z_n = \chi_n U^z \chi_n$, $J^z_n = U^z \chi_n U^z$, $J^z_{n,m} = \chi_m J^z_n \chi_m$ for every $z \in \mathbb{C}$, with the convention that $\chi_{-1} = 0$. Moreover, we denote by
\[\mathrm{Sim}(\bm{\mathfrak{g}},K)_n, \qquad \mathrm{Sim}(\bm{\mathfrak{g}}_z,K)_n \quad (z \in \mathbb{C}),\]
the subsets of $\mathrm{Sim}(\bm{\mathfrak{g}},K)$ and $\mathrm{Sim}(\bm{\mathfrak{g}}_z,K)$, respectively, consisting of elements whose lowest $K$-types have length $l_n$.

Let us fix $n \in \mathbb{N}$. Let $\bm{\pi}_n$ denote the representation of $\Upsilon(\bm{\mathfrak{g}},K)$ on the $(\bm{\mathfrak{g}},K)$-module $\bm{I}_n = \oplus_{\sigma \in \mathcal{C}_n} \bm{I}_\sigma$ and let $B_n = \prod_{\sigma \in \mathcal{C}_n} B_\sigma$. Since $\bm{\pi}_\sigma(\chi_n) = \bm{\pi}_\sigma(\chi_\sigma)$ and hence $\bm{\pi}_\sigma(\bm{U}_n) = \bm{\pi}_\sigma(\bm{U}_\sigma)$ for every $\sigma \in \mathcal{C}_n$, we have $\bm{\pi}_n(\bm{U}_n) \subset  B_n[t]$. In the same way, for every $z \in \mathbb{C}$, we denote by $\pi_n^z$ the representation of $\Upsilon(\bm{\mathfrak{g}}_z,K)$ on the $(\bm{\mathfrak{g}}_z,K)$-module $\oplus_{\sigma \in \mathcal{C}_n} I^z_\sigma$; we have $\pi_n^z(U_n^z) \subset B_n$.

\begin{proposition}\label{isospectral_augmented}
    The morphism $\bm{\pi}_n : \bm{U}_n \to B_n[t]$ is epispectral. Moreover, if we identify $\mathrm{Sim}(\bm{U}_n)$ with a subset of $\mathrm{Sim}(\bm{\mathfrak{g}},K)$ via Lemma \ref{TransferSimple}, we have $$\bm{\pi}_n^*\big( \mathrm{Sim}(B_n[t])\big) = \mathrm{Sim}(\bm{\mathfrak{g}},K)_n.$$
\end{proposition}

\begin{proof}
    From Theorem \ref{epispectral} and Remark \ref{facts_after_proof}(4), we can describe the simple spectrum of $B_n[t] = \prod_{\sigma \in \mathcal{C}_n}B_\sigma[t]$ as follows:
    \begin{equation*}
        \mathrm{Sim}(B_n[t]) = \bigsqcup_{\sigma\in \mathcal{C}_n} \{\chi_\sigma M : M \in \mathrm{Sim}(\bm{\mathfrak{g}},K)_\sigma\}.
    \end{equation*}
    Since $\mathrm{Sim}(\bm{\mathfrak{g}},K)_n = \bigsqcup_{\sigma \in \mathcal{C}_n}\mathrm{Sim}(\bm{\mathfrak{g}},K)_\sigma$ and $\chi_n M = \chi_\sigma M$ for every $\sigma\in \mathcal{C}_n$ and $M\in \mathrm{Sim}(\bm{\mathfrak{g}},K)_\sigma$, we actually have
    \begin{equation*}
        \mathrm{Sim}(B_n[t]) = \{\chi_n M : M \in \mathrm{Sim}(\bm{\mathfrak{g}},K)_n\},
    \end{equation*}
    which concludes the proof.
\end{proof}

Since $\bm{I}_n$ does not contain any $K$-type of length $<l_n$, we have $\bm{\pi}_n(\chi_{n-1}) =  0$. The map $\bm{\pi}_n$ therefore induces by quotient a morphism
\[\tilde{\bm{\pi}}_n : \bm{U}_n/\bm{J}_{n-1,n} \to B_n[t].\]
For the same reason, $\pi^z_n$ induces a morphism
\[\tilde{\pi}^z_n : U^z_n/J^z_{n-1,n} \to B_n\] for any $z \in \mathbb{C}$.

For every $U(K)$-bidmodule $M$ and every finite set of $K$-types $F$, let us write $M_F = \chi_F M\chi_F$. For every $\sigma \in \mathcal{C}_n$, the subalgebra $B_\sigma$ of $\mathrm{End}_\mathbb{C}(E_\sigma)\otimes_\mathbb{C} S$ contains the image of the representation of $U(K)$ on the $K$-module $E_\sigma$, see Definition \ref{def_Bsigma}. We can thus consider $B_n$ as a $U(K)$-bimodule. Since $\tilde{\bm{\pi}}_n$ and $\tilde{\pi}_n^z$ are morphisms of $U(K)$-bimodules, they restrict to algebra homomorphisms
\[\tilde{\bm{\pi}}_{n,F} : (\bm{U}_n/ \bm{J}_{n-1,n})_F \to B_{n,F}[t]\quad \text{and} \quad \tilde{\pi}_{n,F}^z : (U^z_n/J^z_{n-1,n})_F \to B_{n,F},\] for every finite set of $K$-types $F$. Note that if every element of $F$ has length $l_n$, we have $(\bm{U}_n/ \bm{J}_{n-1,n})_F = \bm{U}_F / \bm{J}_{n-1,F}$ and $(U^z_n/J^z_{n-1,n})_F = U^z_F / J^z_{n-1,F}$.

\begin{theorem}\label{big_theorem}
    Let $F$ be a set of $K$-types having length $l_n$ and let \mbox{$z \in \mathbb{C}$}. Let us write $\bm{U}_F(n) = \bm{U}_F/\bm{J}_{n-1,F}$ and $U^z_F(n) = U^z_F/J^z_{n-1,F}$. Consider the commutative diagram
    \[\begin{tikzcd}
        \bm{U}_F(n) \ar[r, "\tilde{\bm{\pi}}_{n,F}"] \ar[d] & B_{n,F}[t] \ar[d, "\mathrm{sp}_z"]\\
        U^z_F(n) \ar[r, "\tilde{\pi}_{n,F}^z"] & B_{n,F}
    \end{tikzcd},\]
    where the left vertical arrow is induced by quotient from the specialization morphism $\mathrm{sp}_z : \bm{U}_F \to U^z_F$. The following holds:
    \begin{enumerate}
        \item the two rings on the left hand side are point-primitive,
        \item the horizontal maps are isospectral and spectrum preserving,
        \item each arrow is a $\mathcal{Z}(\bm{\mathfrak{g}})$-linear map between finitely generated $\mathcal{Z}(\bm{\mathfrak{g}})$-modules,
        \item the map $\tilde{\bm{\pi}}^*_{n,F} : J \mapsto \tilde{\bm{\pi}}_{n,F}^{-1}(J)$ defines a homeomorphism of the primitive spectrum of $B_{n,F}[t]$ onto the primitive spectrum of $\bm{U}_F(n)$,
        \item the horizontal maps induce isomorphisms in periodic cyclic homology.
    \end{enumerate}
\end{theorem}

\begin{proof}
    The statement \textit{(1)} follows from Proposition \ref{criterion_point_primitive}, knowing that every simple $(\bm{\mathfrak{g}},K)$-module is admissible (Corollary \ref{admissibility_simple}).
    
    Since $\bm{U}_F = \chi_F \bm{U}_n \chi_F$ and $B_{n,F}[t] = \bm{\pi}_n(\chi_F) B_n[t] \bm{\pi}_n(\chi_F)$, it follows from Lemma \ref{TransferSimple} and Proposition \ref{isospectral_augmented} that $\bm{\pi}_n : \bm{U}_F \to B_{n,F}[t]$ is epispectral with
    \[\bm{\pi}_n^*\big( \mathrm{Sim}(B_{n,F}[t])\big) = \{ \chi_{F} M : M \in \mathrm{Sim}(\bm{\mathfrak{g}},K)_n,\, \chi_{F} M \neq 0 \},\]
    or, put differently :
    \[\bm{\pi}_n^*\big( \mathrm{Sim}(B_{n,F}[t])\big) = \{ \chi_{F} M : M \in \mathrm{Sim}(\bm{\mathfrak{g}},K),\, \chi_{F} M \neq 0, \, \chi_{n-1}M = 0 \}.\]
    The right hand side of this expression coincides with $\mathrm{Sim}(\bm{U}_F(n))$, which proves that $\tilde{\bm{\pi}}_{n,F}$ is isospectral. Similarly, applying Lemma \ref{isospectral_specializations} to $\bm{\pi}_n$, we obtain that $\pi_n^z : U_F^z \to B_{n,F}$ is epispectral with
    \begin{align*}
        (\pi_n^{z})^*\big(\mathrm{Sim} (B_{n,F})\big) &= \{ \chi_{F} M : M \in \mathrm{Sim}(\bm{\mathfrak{g}}_z,K)_n,\, \chi_{F} M \neq 0 \}\\
        & = \{ \chi_{F} M : M \in \mathrm{Sim}(\bm{\mathfrak{g}}_z,K),\, \chi_{F} M \neq 0, \, \chi_{n-1}M = 0 \}\\
        & = \mathrm{Sim}(U_F^z(n)),
    \end{align*}
    proving that $\pi_{n,F}^z : U_F^z(n) \to B_{n,F}$ is isospectral. By Lemma \ref{isospectral_implies_spectrum_preserving}, the morphisms $\tilde{\bm{\pi}}_{n,F}$ and $\pi_{n,F}^z$ are then also spectrum preserving. So the statement \textit{(2)} is proved.

    Note that Lemma \ref{isospectral_implies_spectrum_preserving} also shows that $\tilde{\bm{\pi}}^*_{n,F}(J)  = \tilde{\bm{\pi}}^{-1}_{n,F}(J)$ for every primitive ideal $J$ of $\mathrm{Prim}(B_{n,F}[t])$.
    
    Let us now check the statement \textit{(3)}. By Proposition \ref{Z(g)_linear}, the horizontal arrows of the diagram are $\mathcal{Z}(\bm{\mathfrak{g}})$-linear (the vertical arrows are $\mathcal{Z}(\bm{\mathfrak{g}})$-linear by definition) and $B_{n,F}[t]$ and $B_{n,F}$ are finitely generated as $\mathcal{Z}(\bm{\mathfrak{g}})$-modules (we use here that $\mathcal{Z}(\bm{\mathfrak{g}})$ is noetherian). The finite generation of $\bm{U}_{F}(n)$ and $U_F^z(n)$ as $\mathcal{Z}(\bm{\mathfrak{g}})$-modules follows from Corollary \ref{finite_generacy_chiUchi}. Applying Theorem \ref{baum_nistor_consequences}, we obtain the remaining statements, \textit{(4)} and~\textit{(5)}.
\end{proof}

We shall now point out some consequences of this theorem. To do so, we introduce some extra conventions. First, we shall endow $\mathbb{C}$ with the Zariski topology, that is, the cofinite topology. Moreover, for every countably dimensional $\mathbb{C}$-algebra~$B$, we consider $\mathrm{Prim}(B[t])$ as a continuous bundle over $\mathbb{C}$ via the map associating to any primitive ideal $J$ of $B$ the element $z \in \mathbb{C}$ such that $J \subset (t-z)B[t]$. Finally, if $I$ is an ideal of a ring $R$, we shall identify $\mathrm{Prim}(R/I)$ with the closed subset of $\mathrm{Prim}(R)$ consisting of primitive ideals $J$ containing $I$.

Due to Proposition \ref{criterion_point_primitive}, since every $(\bm{\mathfrak{g}},K)$-module is admissible (Corollary \ref{admissibility_simple}), the natural map $\mathrm{Sim}(\bm{\mathfrak{g}},K) \to \mathrm{Prim}(\bm{U})$ is a bijection. We endow $\mathrm{Sim}(\bm{\mathfrak{g}},K)$ with the topology making this bijection a homeomorphism. Similarly, for every $z \in \mathbb{C}$, we equip $\mathrm{Sim}(\bm{\mathfrak{g}}_z,K)$ with the topology making $\mathrm{Sim}(\bm{\mathfrak{g}}_z,K) \to \mathrm{Prim}(U^z)$ a homeomorphism. Since $\mathrm{sp}_z : \bm{U} \to U^z$ is surjective, the map $\mathrm{sp}_z^* : \mathrm{Sim}(\bm{\mathfrak{g}}_z,K) \hookrightarrow \mathrm{Sim}(\bm{\mathfrak{g}},K)$ is a closed embedding for every $z \in \mathbb{C}$. In what follows, we consider $\mathrm{Sim}(\bm{\mathfrak{g}},K)$ as a continuous bundle over $\mathbb{C}$, so that the fiber over each $z \in \mathbb{C}$ is $\mathrm{sp}_z^*(\mathrm{Sim}(\bm{\mathfrak{g}}_z,K))$.

For every finite set $F$ of $K$-types having the same length, let us denote by \[\mathrm{Sim}(\bm{\mathfrak{g}},K)_{F}\] the subset of $\mathrm{Sim}(\bm{\mathfrak{g}},K)$ consisting of elements whose set of lowest $K$-type intersects $F$. From Theorem \ref{big_theorem}\textit{(5)}, we deduce the following.

\begin{corollary}\label{remark_transfer_primitive}
     Let $F$ be a finite set of $K$-types having length $l_n$. Let us identify $\mathrm{Prim}(\bm{U}_F(n))$ with $\mathrm{Sim}(\bm{\mathfrak{g}},K)_{F}$ via Proposition \ref{TransferPrimitive}. The pullback by $\tilde{\bm{\pi}}_{n,F}$ defines an isomorphism \[\mathrm{Prim}(B_{n,F}[t]) \to \mathrm{Sim}(\bm{\mathfrak{g}},K)_{F}\] of continuous bundles over $\mathbb{C}$.
\end{corollary}

When $F$ is a one-point set, we can be even more explicit.

\begin{corollary}\label{corollary1}
    Let $\mu$ be any $K$-type. If $\sigma$ denote a cuspidal datum such that $\mu \in \widehat{K}(\sigma)$, the space $\mathrm{Sim}(\bm{\mathfrak{g}},K)_{\{\mu\}}$ is isomorphic as a continuous bundle over $\mathbb{C}$ to the product of affine varieties $\mathbb{C} \times \mathfrak{a}_\sigma^*/W_\sigma$.
\end{corollary}

\begin{proof}
    Let $n \in \mathbb{N}$ be such that the elements of $\widehat{K}(\sigma)$ have length $l_n$. From Definition \ref{def_Bsigma}, one can check that
    \[B_{n,\{\mu\}}[t] = (\chi_\mu B_\sigma \chi_\mu)[t] = \mathrm{End}_\mathbb{C}(E_\mu) \otimes_\mathbb{C} \left(S(\mathfrak{a}_\sigma)^{W_\sigma}\right)[t].\]
    The latter is Morita equivalent to $S(\mathfrak{a}_\sigma)^{W_\sigma}[t]$, which is precisely the algebra of regular functions on $\mathbb{C} \times \mathfrak{a}_\sigma^*/W_\sigma$. Applying Corollary \ref{remark_transfer_primitive} for $F = \{\mu\}$ concludes the proof.
\end{proof}

Let us now explore the consequences for the Mackey bijection $\mathcal{M}$. Let us denote by $\mathrm{Sim}(\mathfrak{g}[t],K)$ the set of isomorphism classes of simple $(\mathfrak{g}[t],K)$-modules. The enveloping algebra of $(\mathfrak{g}[t],K)$ is $U[t]$, see (\ref{BtensorA}). Moreover, the specialization of $(\mathfrak{g}[t],K)$ at any $z \in \mathbb{C}$ is $(\mathfrak{g},K)$. Since every every simple $(\mathfrak{g},K)$-module is admissible and $\mathbb{C}[t]$ acts on every simple $(\mathfrak{g}[t],K)$-module through a character (because of Schur's lemma and the fact that $U$ has countable dimension), every simple $(\mathfrak{g}[t],K)$-module is admissible. It follows that $U[t]$ is point-primitive (Proposition \ref{criterion_point_primitive}). Thus, as above, we can consider $\mathrm{Sim}(\mathfrak{g}[t],K)$ as a topological space naturally homeomorphic to $\mathrm{Prim}(U[t])$. The specialization maps $\mathrm{sp}_z : U[t] \to U$ induce closed embeddings $\mathrm{sp}^*_z : \mathrm{Sim}(\mathfrak{g},K) \to \mathrm{Sim}(\mathfrak{g}[t],K)$ for every $z \in \mathbb{C}$. In that way, if we forget about the topology, $\mathrm{Sim}(\mathfrak{g}[t],K)$ identifies with $\mathbb{C} \times \mathrm{Sim}(\mathfrak{g},K)$ as a bundle over $\mathbb{C}$. For any $n \in \mathbb{N}$, let us denote by $\mathrm{Sim}(\mathfrak{g}[t],K)_n$ the set of elements of $\mathrm{Sim}(\mathfrak{g}[t],K)$ whose lowest $K$-types have length $l_n$.

Let us consider the map $\bm{\mathcal{M}} : \mathrm{Sim}(\mathfrak{g}[t],K) \to \mathrm{Sim}(\bm{\mathfrak{g}},K)$ defined by
\[\bm{\mathcal{M}}(z, J_{\sigma,\lambda}(\mu)) = \begin{cases}
         \mathrm{sp}_z^*\big(J_{\sigma,z^{-1}\lambda}(\mu)\big) & (z\neq 0)\\
         \mathrm{sp}_0^*\big(\mathcal{M}(J_{\sigma,\lambda}(\mu))\big) & (z=0)
     \end{cases},
 \]
for any $z \in \mathbb{C}$, cuspidal datum $\sigma$, $\lambda \in \mathfrak{a}^*_\sigma$ and $\mu \in \widehat{K}(\sigma)$. Then the following holds.

\begin{theorem}
    For every $n \in \mathbb{N}$, the map $\bm{\mathcal{M}}$ restricts to an isomorphism $$\mathrm{Sim}(\mathfrak{g}[t],K)_n \to \mathrm{Sim}(\bm{\mathfrak{g}},K)_n$$ of continuous bundles over $\mathbb{C}$ and in particular, the Mackey bijection restricts to a homeomorphism $\mathrm{Sim}(\mathfrak{g},K)_n \to \mathrm{Sim}(\mathfrak{g}_0,K)_n$
\end{theorem}

\begin{proof}
    By Theorem \ref{big_theorem}, the morphism $\tilde{\pi}^1_n : U_{n}/J^1_{n-1,n} \to B_n$ is isospectral and spectrum preserving and it induces an homeomorphism of primitive spectra. The same is true for the morphism $U_{n}[t]/J^1_{n-1,n}[t] \to B_n[t]$ obtained by extension of scalars, which we also denote by $\tilde{\pi}^1_n$. Now, since $J^1_{n-1,n}[t] = \chi_n U[t]\chi_{n-1} U[t]\chi_n$,
    the primitive spectrum of $U_n[t]/J^1_{n-1,n}[t]$ is homeomorphic to $\mathrm{Sim}(\mathfrak{g}[t],K)_n$ via Proposition \ref{TransferPrimitive}. To sum up, we have an isomorphism $\tilde{\pi}^{1*}_n : \mathrm{Prim}(B_{n}[t]) \to \mathrm{Sim}(\mathfrak{g}[t],K)_n$ of continuous bundles over~$\mathbb{C}$.
    
    On the other hand, by Remark \ref{remark_transfer_primitive}, the morphism $\tilde{\bm{\pi}}_n$ induces an isomorphism $\tilde{\bm{\pi}}_n^* : \mathrm{Prim}(B_{n}[t]) \to \mathrm{Sim}(\bm{\mathfrak{g}},K)_n$ of continuous bundles over $\mathbb{C}$. From the items (3) and (5) of Remark \ref{facts_after_proof}, one can check that $\bm{\mathcal{M}} : \mathrm{Sim}(\mathfrak{g}[t],K)_n \to \mathrm{Sim}(\bm{\mathfrak{g}},K)_n$ coincides with $\tilde{\bm{\pi}}_n^* \circ (\tilde{\pi}^{1*}_n)^{-1}$, which proves the main statement of the theorem. The statement about the Mackey bijection is obtained by considering the fiber at $t = 0$ on both sides.
\end{proof}

Finally, we derive an isomorphism in periodic cyclic homology which is reminiscent of the Connes-Kasparov isomorphism \cite{BaumConnesHigson}*{p.263}, see also \cites{MAhigson, Afgoustidis3} and \cite{Bradd}.

\begin{theorem}\label{iso_HP}
    For every $n\in \mathbb{N}$ and every $z \in \mathbb{C}$, the specialization morphism $\mathrm{sp}_z : \bm{U}_n \to U^z_n$ induces an isomorphism in periodic cyclic homology and in particular, we have an isomorphism $(\mathrm{sp}_1)_* \circ (\mathrm{sp}_0)_*^{-1} : \mathrm{HP}_*(U^0_n) \to \mathrm{HP}_*(U_n)$.
\end{theorem}

For the proof, we shall use the following.

\begin{lemma}\label{invariance_HP_poly}
    If $B$ is a any $\mathbb{C}$-algebra, then for every $z \in \mathbb{C}$, the specialization map $\mathrm{sp}_z : B[t] \to B$ induces an isomorphism in periodic cyclic homology.
\end{lemma}

\begin{proof}
    Note that the natural injection $\iota : B \hookrightarrow B[t]$ is a section of $\mathrm{sp}_z : B[t] \to B$ for every $z \in \mathbb{C}$. It is thus enough to prove that $\iota \circ \mathrm{sp}_0 : B[t] \to B[t]$ induces an isomorphism in periodic cyclic homology.

    For every $z \in \mathbb{C}$, consider the ring homomorphism morphism $\lambda_z : B[t] \to B[t]$ mapping every homogeneous element $b$ of degree $d \in \mathbb{N}$ to $z^db$. Then $(\lambda_z)_{z \in \mathbb{C}}$ is a polynomial family of homomorphisms. By \cite{CuntzQuillen}*{Proposition 8.3}, the maps induced by $\lambda_1 = \mathrm{id}_{B[t]}$ and $\lambda_0 = \iota \circ \mathrm{sp}_0$ in periodic cyclic homology are the same, which concludes the proof.
\end{proof}

\begin{proof}[Proof of Theorem \ref{iso_HP}]
    Let $N \in \mathbb{N}$. We have finite filtrations of $\bm{U}_N$ and $U^z_N$ by ideals:
    \begin{gather*}
        0 = \bm{J}_{-1,N} \subset \bm{J}_{0,N} \subset ... \subset \bm{J}_{N-1,N} \subset \bm{J}_{N,N} = \bm{U}_N,\\
        0 = J^z_{-1,N} \subset J^z_{0,N} \subset ... \subset J^z_{N-1,N} \subset J^z_{N,N} = U^z_N.
    \end{gather*}
    For every $n \in \mathbb{N}$, let us denote by $\tilde{\mathrm{sp}}_z$ the map $\bm{J}_{n,N}/\bm{J}_{n-1,N} \to J^z_{n,N}/J^z_{n-1,N}$ induced from $\mathrm{sp}_z : \bm{U}_N \to U_N^z$ by subquotient.
    By the excision property of periodic cyclic homology \cite{CuntzQuillenExcision}*{§5} (short exact sequences of algebras induce six-term exact sequences in periodic cyclic homology, in a functorial way), to prove that $\mathrm{sp}_z : \bm{U}_N \to U^z_N$ induces an isomorphism in periodic cyclic homology, it is enough the prove that
    \begin{equation}\label{map_HP}
        (\tilde{\mathrm{sp}}_z)_* : \mathrm{HP}_*(\bm{J}_{n,N}/\bm{J}_{n-1,N}) \to \mathrm{HP}_*(J^z_{n,N}/J^z_{n-1,N})
    \end{equation}
    is an isomorphism for every $n \in \mathbb{N} \cap [0,N]$.

    Let $n \in \mathbb{N} \cap [0,N]$. Let us write $\bm{J}'_{n,N} = \bm{J}_{n,N}/\bm{J}_{n-1,N}$ and $J'^z_{n,N} = J^z_{n,N}/J^z_{n-1,N}$. We have
    $$\bm{J}'_{n,N} \chi_{n} \bm{J}'_{n,N} = \bm{J}'_{n,N}, \qquad \chi_{n}\bm{J}'_{n,N}\chi_{n} = \bm{U}_n/\bm{J}_{n-1,n}$$
    and
    $$J'^z_{n,N}\chi_{n} J'^z_{n,N} = J'^z_{n,N},\qquad \chi_{n}J'^z_{n,N}\chi_{n}= U^z_{n}/J^z_{n-1,n}.$$
    From the Morita invariance of periodic cyclic homology \cite{CuntzMorita}, it follows that the map induced by $\tilde{\mathrm{sp}}_z : \bm{J}'_{n,N} \to J'^z_{n,N}$ in periodic cyclic homology will be an isomorphism if this is the case for the map induced by $\tilde{\mathrm{sp}}_z : \bm{U}_n/\bm{J}_{n-1,n} \to U^z_{n}/J^z_{n-1,n}$.

    By Theorem \ref{big_theorem}\textit{(5)} and Lemma \ref{invariance_HP_poly}, in the commutative diagram
    \[
    \begin{tikzcd}
        \mathrm{HP}_*(\bm{U}_n/\bm{J}_{n-1,n}) \ar[rr, "(\tilde{\bm{\pi}}_n)_*"] \ar[d, "(\tilde{\mathrm{sp}}_z)_*"'] && \mathrm{HP}_*(B_{n}[t]) \ar[d, "(\mathrm{sp}_z)_*"]\\
        \mathrm{HP}_*(U^z_{n}/J^z_{n-1,n}.) \ar[rr, "(\tilde{\pi}_n^z)_*"] & &\mathrm{HP}_*(B_{n})
    \end{tikzcd},
    \]
    the horizontal and right vertical arrows are isomorphisms. It follows that the left vertical arrow is also an isomorphism, which concludes the proof.
\end{proof}

\appendix

\section{Background on approximately unital rings}\label{AUA}
This first appendix gathers some definitions and elementary facts about approximately unital rings. We refer to \cite{Renard}*{Chapitre 1} and \cite{KnappVogan}*{Chapters I, II} for more information. In the following, rings are only assumed to be associative.

Let $A$ be a ring. Let $\mathrm{Idem}(A)$ denote the set of idempotents of $A$. We equip $\mathrm{Idem}(A)$ with the order defined by $e \prec p \Leftrightarrow e = pep$. The ring $A$ is said to be \textbf{approximately unital} if every finite subset of $A$ is contained in a subring of the form $eAe$ for some $e \in \mathrm{Idem}(A)$.

Assume that $A$ is approximately unital. Then $\mathrm{Idem}(A)$ is directed. An \textbf{approximate unit} of $A$ is defined as a increasing net $(e_i)_{i \in I}$ of idempotents of $A$ such that $\{e_i : i \in I\}$ is cofinal in $\mathrm{Idem}(A)$. An $A$-module $M$ is said to be \textbf{nondegenerate} if $M = AM$, that is, if $M = A \otimes_AM$. Note that any submodule of a nondegenerate module is itself nondegenerate. If $M$ is any $A$-module, its \textbf{nondegenerate part} is defined as $AM$; it is the largest nondegenerate submodule of $M$.

For the rest of the section, all modules over approximately unital algebras are assumed to be nondegenerate.

The set of isomorphism classes of simple $A$-modules, called the simple spectrum of $A$, is denoted by $\mathrm{Sim}(A)$. A \textbf{primitive ideal} of $A$ is by definition the annihilator of a simple $A$-module. The set of all primitive ideals of $A$ is referred to as the \textbf{primitive spectrum} of $A$ and is denoted by $\mathrm{Prim}(A)$. We equip $\mathrm{Prim}(A)$ with the Jacobson topology, see \cite{Jacobson}.

Let $e$ be an idempotent of $A$. Let us denote by $\mathrm{Sim}_e(A)$ the set of isomorphism classes of simple $A$-modules not annihilating $e$. The following lemma relates the simple $eAe$- modules to simple $A$-modules (see \cite{Renard}*{p.23} for a proof).

\begin{lemma}\label{TransferSimple}
    If $M$ is a simple $A$-module which does not annihilate $e$, then $eM$ is a simple $eAe$-module. This defines a bijection $\mathrm{Sim}_e(A) \to \mathrm{Sim}(eAe)$.
\end{lemma}

This lemma can be restated in terms of primitive ideals as follows. For every ideal $J$ of $eAe$, let us write $J_A = \{ a \in A : eAaAe \subset J\}$. Let $\mathrm{Prim}_e(A)$ denote the set of primitive ideals $I$ of $A$ such that $e \notin I$.

\begin{proposition}\label{TransferPrimitive}
    $\mathrm{Prim}_e(A)$ is an open subset of $\mathrm{Prim}(A)$. The map $I \mapsto (eAe)\cap I$ is a homeomorphism of $\mathrm{Prim}_e(A)$ onto $\mathrm{Prim}(eAe)$ with inverse $J \mapsto J_A$.
\end{proposition}

\begin{proof}
    The complement of $\mathrm{Prim}_e(A)$ in $\mathrm{Prim}(A)$ is the set of primitive ideals of $A$ containing $e$, hence it is a closed subset of $\mathrm{Prim}(A)$. This proves the first statement.
    
    Let $M$ be a simple $A$-module such that $eM \neq 0$ and let $I$ be its annihilator in~$A$. By Lemma \ref{TransferSimple}, $eM$ is a simple $eAe$-module. Its annihilator, $(eAe)\cap I$, is therefore a primitive ideal of $eAe$. Moreover, for every $a \in A$, we have
    \[a \in ((eAe) \cap I)_A \iff eAaAe M = 0.\]
    Since $eM \neq 0$ and $M$ is simple, we have
    \[AaAeM = AaM = \begin{cases}
        M & \text{if $aM \neq 0$},\\
        0 & \text{if $aM = 0$}.
    \end{cases}\]
    It follows that
    \[a \in ((eAe) \cap I)_A \iff eAaM = 0 \iff AaM = 0 \iff aM = 0,\]
    proving that $((eAe)\cap I)_A = I$. Moreover, Lemma \ref{TransferSimple} shows that every simple $eAe$-module is isomorphic to $eM$ for some simple $A$-module $M$. Therefore, the map $I \mapsto (eAe)\cap I$ is a bijection of $\mathrm{Prim}_e(A)$ onto $\mathrm{Prim}(eAe)$ with inverse $J \mapsto J_A$.
    
    Defined on all ideals of $A$ and $eAe$ respectively, the maps $I \mapsto (eAe)\cap I$ and $J \mapsto J_A$ preserve inclusions and arbitrary intersections. It follows that the map $\mathrm{Prim}_e(A) \to \mathrm{Prim}(eAe)$ preserve closures and hence that it is a homeomorphism.
\end{proof}

From now on, let us also assume that $A$ is an algebra over an algebraically closed field $k$. Every $A$-module $M$ is naturally isomorphic to $A \otimes_A M$, hence every $A$-module is naturally equipped with a $k$-vector space structure. An $A$-module $M$ is said to be \textbf{admissible} if $eM$ is finite-dimensional over $k$, for every idempotent $e$ of $A$. Recall from Subsection \ref{spectral_properties_rings} the meaning of point-primitive.

\begin{proposition} \label{criterion_point_primitive}
    If each simple $A$-module is admissible, then $A$ is point-primitive.
\end{proposition}

\begin{proof}
    We need to prove that for any primitive ideal $I$ of $A$, the quotient ring $A/I$ has up to isomorphism a unique simple module. Let $M$ be a simple $A$-module and let $I$ be the annihilator of $M$. We identify $A/I$ with the image of $A$ in $\mathrm{End}_k(M)$. For every idempotent $e$ of $A$, let us identify $e\mathrm{End}_k(M)e$ with $\mathrm{End}_k(eM)$, via restriction to $eM$. If $eM \neq 0$, the $eAe$-module $eM$ is simple by Lemma \ref{TransferSimple} and finite-$k$-dimensional by assumption; as a consequence of Jacobson density theorem and Schur's lemma, we thus have $eAe/(I\cap(eAe)) = \mathrm{End}_k(eM)$. It follows that $A/I$ is the increasing union of all $\mathrm{End}_k(eM)$, where $e$ ranges over $\mathrm{Idem}(A)$. Consequently, if $p \in A/I$ is any rank-one projection, we have $A/I = (A/I) p(A/I)$ and $A/I$ is Morita equivalent to $p(A/I)p = kp$ by \cite{AnhMarki}*{Proposition 3.5}. Since $A/I$ and $k$ have equivalent categories of modules, there is indeed a unique simple $A/I$-module up to isomorphism.
\end{proof}

\section{Relations between various inductions}\label{relInd}

In this section, we consider a Lie group $G$ and a closed subgroup $H \subset G$, as well as a compact subgroup $K \subset G$ such that $G = KH$. We write $L = K\cap H$. Our goal is to compare the induction functor $\mathrm{ind}_{\mathfrak{h}, L}^{\mathfrak{g}, K}$ defined in Subsection \ref{changeofring} with the normalized induction for topological group representations \cite{vandenBan}*{§1.6}.

Recall that for a given $K$-admissible \cite{Warner}*{Definition p. 304} continuous representation $\sigma$ of $G$ on a Banach space $E$, all $K$-finite vectors of $E$ are smooth, see \cite{Warner}*{Lemma 4.5.5.1}. The space of $K$-finite vectors $E_K$ thus inherits a $(\mathfrak{g},K)$-module structure. We then refer to it as the $(\mathfrak{g},K)$-module associated to $\sigma$ and denote it by $\sigma_{\mathfrak{g},K}$. The map $\sigma \mapsto \sigma_{\mathfrak{g},K}$ has several key properties. If $(\sigma, E_\sigma)$ and $(\tau, E_\tau)$ are two continuous Banach representations of $G$ which are Naimark equivalent \cite{Warner}*{Definition p. 232} via a densely defined closed linear map $\phi : E_\sigma \to E_\tau$, then $\phi$ restricts to an isomorphism of $(\mathfrak{g},K)$-modules $\sigma_{\mathfrak{g},K} \to \tau_{\mathfrak{g},K}$. For this, we refer to \cite{Warner}*{Theorem 4.5.5.2}. Moreover, for every continuous Banach representation $(\sigma, E_\sigma)$ of $G$, the $(\mathfrak{g},K)$-module associated to its contragredient representation $(\sigma^*, E_\sigma^*)$ identifies with the contragredient $\sigma^*_{\mathfrak{g},K}$ of $\sigma_{\mathfrak{g},K}$.

Regarding topological induction, we use \cite{vandenBan}*{Section 1} as a reference. Let $\xi$ be a continuous representation of $H$ on a Hilbert space $V$. The representation of $G$ \textbf{(co)induced} from $V$ (see \cite{vandenBan}*{§1.1-§1.2}), denoted by $\text{C}_H^G(\xi)$, is defined as the space $C(G,\xi)$ of $V$-valued continuous right $\xi$-equivariant functions on $G$, endowed with the compact open topology and with the continuous action of $G$ by left translations. Because of the right equivariance, the hypothesis $KH = G$ implies that $C(G,\xi)$ is normable and hence a Banach space (consider for instance the norm of uniform convergence on $K$). For every $\alpha \in \mathbb{R}$, let $\delta^{\alpha}$ denote the representation of $H$ on the space of complex-valued $\alpha$-densities on $\mathfrak{g}_\mathbb{R}/\mathfrak{h}_\mathbb{R}$ via the character $h  \mapsto |\mathrm{det}\, \mathrm{Ad}_{\mathfrak{g}_\mathbb{R}/\mathfrak{h}_\mathbb{R}}(h)|^{-\alpha}$. As explained in \cite{vandenBan}*{§1.6}, the space $C(G, \xi\otimes\delta^{1/2})$ inherits a natural inner product for which the action of $G$ is continuous; the representation of $G$ on the associated Hilbert completion is said to be obtained from $V$ by \textbf{normalized induction} and is denoted by $\text{I}_H^G(\xi)$.

Let us point out two standard facts. First, if $\xi$ is $L$-admissible then both $\mathrm{C}_H^G(\xi)$ and $\mathrm{I}_H^G(\xi)$ are admissible. Indeed, the hypothesis $G= KH$ implies that the restriction of $\mathrm{C}_H^G(\xi)$ and $\mathrm{I}_H^G(\xi)$ to $K$ are $\mathrm{C}_L^K(\xi_{|L})$ and $\mathrm{I}_L^K(\xi_{|L})$, respectively; we can then apply Frobenius reciprocity to conclude that the $K$-type multiplicities are finite. Second, the operation $\mathrm{I}_H^G$ commutes with the operation consisting in taking the contragredient representation: $\mathrm{I}_H^G(\xi^*)$ naturally identifies with $\mathrm{I}_H^G(\xi)^*$ (this explicitly use the Hilbert space structure of $V$). The goal of this appendix is to prove the following theorem.

\begin{theorem}\label{comparisonInd}
    Let $\xi$ be an $L$-admissible continuous representation of $H$ on some Hilbert space. There is an isomorphism of $(\mathfrak{g},K)$-modules
    \[\mathrm{I}_G^H(\xi)_{\mathfrak{g},K} \approx \mathrm{ind}_{\mathfrak{h},L}^{\mathfrak{g},K}(\xi_{\mathfrak{h},L} \otimes \delta^{-1/2}_{\mathfrak{h},L}).\]
\end{theorem}

Before diving into the proof, let us introduce another functor. For any morphism $(\mathfrak{i},M) \to (\mathfrak{j},N)$ of Harish-Chandra pairs over $\mathbb{C}$, let $I^{\mathfrak{j},N}_{\mathfrak{i},M}$ be the functor from the category of $(\mathfrak{i},M)$-modules to the category of $(\mathfrak{j},N)$-modules defined by
\[\text{coind}^{\mathfrak{j},N}_{\mathfrak{i},M}: W \longmapsto \mathrm{Hom}_{U(\mathfrak{i},M)}(U(\mathfrak{j},N), W)_\mathrm{ndg},\]
where the subscript $\mathrm{ndg}$ refers to the non-degenerate part of the considered $U(\mathfrak{j},N)$-module; see \cite{KnappVogan}*{Formula (2.10)}.
The above functor is the right adjoint of the forgetful functor $\mathrm{res}^{\mathfrak{j},N}_{\mathfrak{i},M}$ restricting a $(\mathfrak{j},N)$-module to a $(\mathfrak{i},M)$-module, see \cite{KnappVogan}*{Proposition 2.21}. More precisely, we have a natural transformation \[\varepsilon : \mathrm{res}^{\mathfrak{j},N}_{\mathfrak{i},M}\circ \text{coind}^{\mathfrak{j},N}_{\mathfrak{i},M} \to \mathrm{id}_{\mathfrak{i},M},\]
the latter being the identity functor of the category of $(\mathfrak{i},M)$-module, which implements a natural isomorphism
\begin{equation}\label{Frobenius}\begin{aligned}\mathrm{Hom}_{\mathfrak{j},N}\big(X, \text{coind}^{\mathfrak{j},N}_{\mathfrak{i},M}(W)\big) &\approx \mathrm{Hom}_{\mathfrak{i},M}\big(\mathrm{res}^{\mathfrak{j},N}_{\mathfrak{i},M}(X), W\big)\\
\phi& \mapsto \varepsilon(W) \circ \phi.\end{aligned} 
\end{equation}

In order to prove the above theorem, we will need the following proposition.

\begin{proposition}\label{intermediateComparison}
    Let $\xi$ be an $L$-admissible continuous representation of $H$ on some Hilbert space $V$. There is an isomorphism of $(\mathfrak{g},K)$-modules
    \[\mathrm{C}_G^H(\xi)_{\mathfrak{g},K} \approx \mathrm{coind}_{\mathfrak{h},L}^{\mathfrak{g},K}(\xi_{\mathfrak{h},L}).\]
\end{proposition}

\begin{proof}
    Consider the evaluation $\mathrm{ev}_e : \mathrm{ind}_H^G(\xi) \to V$ at the identity element of $G$. Since it is $H$-equivariant, this defines a morphism of $(\mathfrak{h},L)$-modules \begin{equation*}
        \mathrm{ev}_e : \mathrm{C}_H^G(\xi)_{\mathfrak{g},K} \to \xi_{\mathfrak{h},L}.
    \end{equation*}
    Applying the infinitesimal Frobenius reciprocity (\ref{Frobenius}) to the latter, we obtain a morphism of $(\mathfrak{g},K)$-modules
    \[\widetilde{\mathrm{ev}_e} : \mathrm{C}_H^G(\xi)_{\mathfrak{g},K} \longrightarrow \text{coind}_{\mathfrak{h},L}^{\mathfrak{g},K}(\xi_{\mathfrak{h},L})\]
    such that $\varepsilon(\xi_{\mathfrak{h},L}) \circ \widetilde{\mathrm{ev}_e} =  \mathrm{ev}_e$. Our goal is to prove that $\widetilde{\mathrm{ev}_e}$ is an isomorphism.

    Let us first prove that $\widetilde{\mathrm{ev}_e}$ is injective. Let $\psi$ be an element of $\mathrm{C}_H^G(\xi)_{\mathfrak{g},K}$ such that $\widetilde{\mathrm{ev}_e}(\psi) = 0$. Then for every $k \in K$, we have \[\psi(k) = \mathrm{ev}_e(k^{-1}\psi) = \varepsilon(\xi_{\mathfrak{h},L})(\widetilde{\mathrm{ev}_e}(k^{-1}\psi)) = \varepsilon(\xi_{\mathfrak{h},L})(k^{-1}\widetilde{\mathrm{ev}_e}(\psi)) = 0.\] Since $G = KH$ and $\psi$ is right $H$-equivariant, we have $\psi = 0$. This proves the injectivity.

    In order to conclude that $\widetilde{\mathrm{ev}_e}$ is an isomorphism, it is enough to check that the two involved $(\mathfrak{g},K)$-modules have the same $K$-type multiplicities. Since $\mathfrak{g} = \mathfrak{h} + \mathfrak{k}$ and $\mathfrak{l} = \mathfrak{h} \cap \mathfrak{k}$, the multiplication $U(\mathfrak{h}) \times U(\mathfrak{k}) \to U(\mathfrak{g})$ induces an isomorphism $U(\mathfrak{h}) \otimes_{U(\mathfrak{l})} U(\mathfrak{k}) \to U(\mathfrak{g})$. Moreover, by Proposition \ref{structureU}, we have a linear isomorphism $U(\mathfrak{g}) \otimes_{U(\mathfrak{k})} U(K) \to U(\mathfrak{g},K)$. It follows that the multiplication in $\Upsilon(\mathfrak{g},K)$ induces a linear isomorphism \begin{equation*} U(\mathfrak{h}) \otimes_{U(\mathfrak{l})} U(K) \longrightarrow U(\mathfrak{g},K).\end{equation*} Let $V_L$ denote the space of $L$-finite vectors of $V$, viewed as a representation of $L$. Because of the above isomorphism, the restriction of
    \(\text{coind}_{\mathfrak{h},L}^{\mathfrak{g},K}(\xi_{\mathfrak{h},L})\) to a representation of $K$ identifies with $\text{coind}^{K}_{L}(V_L)$, see Convention \ref{simpler_notation_HC_pairs}. By (\ref{Frobenius}), for every $K$-type $\mu$, we have a linear isomorphism
    \[\mathrm{Hom}_{K}(\mu, \text{coind}^{K}_{L}(V_L)) \approx \mathrm{Hom}_{L}(\mu_{|L}, V_L).\]
    Hence, the multiplicity of each $K$-type $\mu$ in $\text{coind}_{\mathfrak{h},L}^{\mathfrak{g},K}(\xi_{\mathfrak{h},L} )$ is the dimension of the vector space $\mathrm{Hom}_{L}(\mu_{|L}, V_L)$.

    A similar analysis can be done for $\mathrm{C}_H^G(\xi)_{\mathfrak{g},K}$. Since $G = KH$, the restriction of $\mathrm{C}_H^G(\xi)$ to $K$ identifies with $\mathrm{C}_L^K(\xi_{|L})$. Applying the Frobenius reciprocity \cite{vandenBan}*{Lemma 1.2} to the latter, we get that the multiplicity of any $K$-type $\mu$ in $\mathrm{C}_H^G(\xi)$ is equal to the dimension of $\mathrm{Hom}_L(\mu, \xi_{|L})$. The $K$-type multiplicities of $\text{coind}_{\mathfrak{h},L}^{\mathfrak{g},K}(\xi_{\mathfrak{h},L} )$ and $\mathrm{C}_H^G(\xi)_{\mathfrak{g},K}$ are thus the same, which concludes the proof.
\end{proof}

\begin{proof}[Proof of Theorem \ref{comparisonInd}]
    Since $\xi$ is $L$-admissible, we have $\xi^{**}= \xi$, hence \[\mathrm{I}_H^G(\xi)_{\mathfrak{g},K} = (\mathrm{I}_H^G(\xi^*)^*)_{\mathfrak{g},K} = (\mathrm{I}_H^G(\xi^*)_{\mathfrak{g},K})^*.\]
    The dense injection $\mathrm{C}_H^G(\xi^* \otimes \delta^{1/2}) \hookrightarrow \mathrm{I}_H^G(\xi^*)$ is a Naimark equivalence of Banach representations of $G$, therefore the $(\mathfrak{g},K)$-module associated with $\mathrm{I}_H^G(\xi^*)$ identifies with that of $\mathrm{C}_H^G(\xi^* \otimes \delta^{1/2})$. Consequently, by Proposition \ref{intermediateComparison}, we have
    \[\mathrm{I}_H^G(\xi^*)_{\mathfrak{g},K} = \mathrm{C}_H^G(\xi^* \otimes \delta^{1/2})_{\mathfrak{g},K} \approx \mathrm{coind}^{\mathfrak{g},K}_{\mathfrak{h},L}\left[(\xi^*)_{\mathfrak{h},L} \otimes \delta^{1/2}_{\mathfrak{h},L}\right].\]
    On the other hand, by \cite{KnappVogan}*{Theorem 3.1}, we have an isomorphism
    \[\mathrm{coind}^{\mathfrak{g},K}_{\mathfrak{h},L}\left[(\xi^*)_{\mathfrak{h},L} \otimes \delta^{1/2}_{\mathfrak{h},L}\right]^* \approx \mathrm{ind}^{\mathfrak{g},K}_{\mathfrak{h},L}\left(\left[(\xi^*)_{\mathfrak{h},L} \otimes \delta^{1/2}_{\mathfrak{h},L}\right]^*\right).\]
    Since $\left[(\xi^*)_{\mathfrak{h},L} \otimes \delta^{1/2}_{\mathfrak{h},L}\right]^* \approx \xi_{\mathfrak{h},L} \otimes \delta^{-1/2}_{\mathfrak{h},L}$, we eventually get an isomorphism
    \[\mathrm{I}_H^G(\xi)_{\mathfrak{g},K}\approx \mathrm{ind}^{\mathfrak{g},K}_{\mathfrak{h},L}(\xi_{\mathfrak{h},L} \otimes \delta^{-1/2}_{\mathfrak{h},L}).\]
\end{proof}


\begin{thebibliography}{99}
\bib{AtiyahMacDonald}{book}{
    author = {Atiyah, M. F.},
    author = {MacDonald, I. G.},
    title = {Introduction To Commutative Algebra},
    date = {1969},
    publisher = {CRC Press}
}
\bib{Afgoustidis3}{article}{
        author = {Afgoustidis, A.},
        title = {On the analogy between real reductive groups and Cartan motion groups: A proof of the Connes-Kasparov isomorphism},
        journal = {Journal of Functional Analysis},
        volume = {277, no.7},
        date = {2019},
        pages = {2237--2258}
    }
\bib{Afgoustidis1}{article}{
        author = {Afgoustidis, A.},
        title = {On the analogy between real reductive groups and Cartan motion groups: The Mackey-Higson bijection},
        journal = {Cambridge Journal of Mathematics},
        volume = {9, no.3},
        date = {2021},
        pages = {551--575}
    }
    \bib{AfgoustidisAubert}{article}{
        author = {Afgoustidis, A.},
        author = {Aubert, A.},
        title = {Continuity of the Mackey-Higson bijection},
        journal = {Pacific Journal of Mathematics},
        volume = {310, no.2},
        date = {2021},
        pages = {257-273}
    }
    \bib{AfgoustidisClare}{article}{
        author = {Afgoustidis, A.},
        author = {Clare, P.},
        title = {The Mackey bijection as a stratified equivalence},
        journal = {arXiv:2602.22798},
        date = {2026}
    }
    \bib{AnhMarki}{article}{
    author = {\'Anh, P. N.},
    author = {M\'arki, L.},
    title = {Morita equivalence for rings without identity},
    journal = {Tsukuba Journal of Mathematics},
    volume = {11, no. 1},
    date = {1987},
    pages = {1-16}
    }
    \bib{BaumConnesHigson}{article}{
        title = {Classifying space for proper actions and K-theory of group C*-algebras},
        author = {Baum, P.},
        author = {Connes, A.},
        author = {Higson, N.},
        pages = {240--291},
        book = {
            title = {C*-algebras: 1943–1993 (San Antonio, TX, 1993)},
            series ={Contemporary Mathematics},
            volume = {167},
            publisher = {American Mathematical Society, Providence, RI},
            date = {1994}
        }
    }
    \bib{BaumNistor}{article}{
        author = {Baum, P.},
        author = {Nistor, V.},
        title = {Periodic cyclic homology of Iwahori–Hecke algebras},
        journal = {K-theory},
        volume = {27, no. 4},
        date = {2002},
        pages = {329-357}
    }
    \bib{BHS}{article}{
        author = {Bernstein, J.},
        author = {Higson, N.},
        author = {Subag, E.},
        title = {Algebraic Families of Harish-Chandra pairs},
        journal = {International Mathematics Research Notices},
        date = {2020},
        volume = {2020, no. 15},
        pages = {4776-4808}
        }
    \bib{Bourbaki}{book}{
        title = {Algebra},
        subtitle = {Chapter 8},
        author = {Bourbaki, N.},
        publisher = {Springer Cham},
        date = {2022}
    }
    \bib{GrAlLie1}{book}{
        title = {Groupes et algèbres de Lie},
        subtitle = {Chapitre 1},
        author = {Bourbaki, N.},
        publisher = {Springer Berlin, Heidelberg},
        date = {2007}
    }
    \bib{Bradd}{article}{
        title = {Compatible decomposition of the Casselman algebra and the reduced group C*-algebra of a real reductive group},
        author = {Bradd, J.},
        journal = {Cambridge Journal of Mathematics},
        volume = {13, no. 1},
        date = {2025},
        pages ={1-50}
    }
    \bib{BHY_operator}{article}{
        author = {Bradd, J.},
        author = {Higson, N.},
        author = {Yuncken, R.},
        title = {Operator K-Theory and Tempiric Representations},
        journal = {arXiv:2412.18924},
        date = {2024}
}
    \bib{CuntzMorita}{article}{
        title = {Morita Invariance in Cyclic Homology for Nonunital Algebras},
        author = {Cuntz, J.},
        journal = {K-theory},
        volume = {15},
        pages = {301-305},
        date = {1998}
    }
    \bib{CuntzQuillen}{article}{
        title = {Cyclic homology and singularity},
        author = {Cuntz, J.},
        author = {Quillen, D.},
        journal = {Journal of the American Mathematical Society},
        volume = {8, no. 2},
        date = {1995},
        pages = {373-442}
    }
    \bib{CuntzQuillenExcision}{article}{
        title = {Excision in bivariant periodic cyclic cohomology},
        author = {Cuntz, J.},
        author = {Quillen, D.},
        journal = {Inventiones Mathematicae},
        volume = {127},
        pages = {67-98},
        date = {1997}
    }
    \bib{DelormeCh}{article}{
        title = {Sur les représentations des groupes de déplacements de Cartan},
        author = {Champetier, C.},
        author = {Delorme, P.},
        journal = {Journal of Functional Analysis},
        volume = {43},
        date = {1981},
        pages = {258-279}
        }
    \bib{Membedding}{article}{
        title = {A Mackey embedding for reduced C*-algebras of real reductive groups},
        author = {Clare, P.},
        author = {Higson, N.},
        author= {Roman, A.},
        journal = {arXiv:2508.02504},
        date = {2025}
    }
    \bib{Delorme84}{article}{
        title = {Homomorphismes de Harish-Chandra liés aux $K$-types minimaux des séries principales généralisées des groupes de Lie réductifs connexes},
        author = {Delorme, P.},
        journal = {Annales Scientifiques de l'ENS},
        volume = {17, no.1},
        date = {1984}}
    \bib{DelormeAnnals}{article}{
        title = {Sur le théorème de Paley-Wiener d'Arthur},
        author = {Delorme, P.},
        journal = {Annals of Mathematics},
        volume = {162, no. 2},
        pages = {987-1029},
        date = {2005}
    }
    
    \bib{Dixmier}{book}{
        title = {C*-algebras},
        author = {Dixmier, J.},
        publisher = {Elsevier},
        series = {North-Holland Mathematical Library},
        volume ={15},
        date = {1977}
    }
    \bib{GaudillotEstrada}{article}{
        title = {Covariant representations of algebraic group actions and applications},
        author = {Gaudillot-Estrada, Y.},
        journal = {arXiv:2603.06259},
        date = {2026}
    }
    \bib{MAhigson}{article}{
        title = {The Mackey analogy and $K$-theory},
        author = {Higson, N.},
        pages = {149--172},
        book = {
            title = {Group representations, ergodic theory, and mathematical physics: a tribute to George W. Mackey},
            series = {Contemporary Mathematics},
            volume = {449},
            publisher = {American Mathematical Society},
            address = {Providence, RI},
            date = {2008}
        }
    }
    \bib{Jacobson}{article}{
        title = {A Topology for the Set of Primitive Ideals in an Arbitrary Ring},
        author = {Jacobson, N.},
        journal = {Proceedings of the National Academy of Sciences of the United States of America},
        volume = {31, no. 10},
        date = {1945},
        pages = {333-338}
    }
    \bib{Knapp}{book}{
        title = {Lie Groups Beyond an Introduction},
        author = {Knapp, A.},
        publisher = {Birkhäuser Boston, MA},
        series = {Progress in Mathematics},
        date = {1996}
    }
    \bib{KnappVogan}{book}{
        title = {Cohomological Induction and Unitary Representations},
        author = {Knapp, A.},
        author = {Vogan, D.},
        publisher = {Princeton University Press},
        series = {Princeton Mathematical Series},
        volume = {45},
        date = {1995}
    }
    \bib{Loday}{book}{
        title = {Cyclic homology},
        author = {Loday, J.},
        publisher = {Springer Berlin, Heidelberg},
        date = {1998}
    }
    \bib{Mukai}{book}{
        title = {An Introduction to Invariants and Moduli},
        author = {Mukai, S.},
        publisher = {Cambridge University Press},
        date = {2003}
    }
    \bib{MackeyAnalogy}{article}{
        title = {On the analogy between semisimple Lie groups and certain related semidirect product groups},
        author = {Mackey, G. W.},
        pages = {339--363},
        book = {
            title = {Lie groups and their representations (Proc. Summer School, Bolyai János Math. Soc., Budapest, 1971)},
            publisher = {Halsted},
            address = {New York},
            date = {1975}
        }
    }
\bib{Renard}{book}{
        title = {Représentations des groupes réductifs $p$-adiques},
        author = {Renard, D.},
        publisher = {Société mathématique de France},
        series = {Cours spécialisés},
        volume ={17},
        date = {2010}
    }
\bib{vandenBan}{article}{
    title = {Induced representations and the Langlands classification},
    author = {van den Ban, E. P.},
    pages = {123-155},
    book = {
        title = {Representation Theory and Automorphic Forms},
        publisher = {American Mathematical Society},
        series = {Proceedings of Symposia in Pure Mathematics},
        volume = {61},
        date = {1997},
        address = {Providence, RI}
    }
}
\bib{VoganAnnals}{article}{
    title = {The algebraic structure of the representation of semisimple Lie groups I},
    author = {Vogan, D.},
    journal = {Annals of Mathematics},
    volume = {109, no. 1},
    date = {1979},
    pages = {1-60}}
\bib{VoganGreen}{book}{
        title = {Representations of Real Reductive Lie Groups},
        author = {Vogan, D.},
        series = {Progress in Mathematics},
        volume = {15},
        publisher = {Birkhäuser},
        address = {Boston, MA},
        date = {1981}
    }
\bib{Wallach}{book}{
        title = {Real reductive groups I},
        author = {Wallach, N.},
        date = {1988},
        publisher = {Academic Press}
    }
\bib{Warner}{book}{
        title = {Harmonic Analysis on Semi-Simple Lie Groups I},
        author = {Warner, G.},
        publisher = {Springer Berlin, Heidelberg},
        date = {1972}
        }
\bib{Wassermann}{article}{
    title = {Une démonstration de la conjecture de Connes-Kasparov pour
les groupes de Lie linéaires connexes réductifs},
    author = {Wassermann, A.},
    journal = {Comptes Rendus de l'Académie des Sciences de Paris},
    volume = {304, no. 18},
    date = {1987},
    pages = {559-562}
}
\bib{Yu}{article}{
    title = {Mackey analogy as deformation of {$\mathcal{D}$}-modules},
    author = {Yu, S.},
    journal = {Mathematische Annalen},
    volume = {385},
    pages = {421-457},
    date = {2023}
}
\end{thebibliography}
\end{document}